\documentstyle[amscd]{amsart}
\newsymbol\boxtimes 1202

\newcommand{\nc}{\newcommand}

\nc{\ad}{{{\bf{ad}}}}
\nc{\AJ}{{\operatorname{aj}}}
\nc{\Aut}{{\operatorname{Aut}}}
\nc{\Bls}{{{\cal B}ls}}
\nc{\Boxtimes}{{\fbox{$\times$}}}
\nc{\blt}{{\bullet}}
\nc{\bSt}{{\mbox{\bf{St}}}}
\nc{\card}{{\operatorname{card}}}
\nc{\Cch}{{\check{C}}}
\nc{\cd}{{\operatorname{cd}}}
\nc{\Ch}{{\operatorname{Ch}}}
\nc{\chara}{{\operatorname{char}}}
\nc{\CHom}{{\cal{H}om}}
\nc{\Coker}{{\operatorname{Coker}}}
\nc{\codim}{{\operatorname{codim}}}
\nc{\Cone}{{\operatorname{Cone}}}
\nc{\cSgn}{{\cal{S}gn}}
\nc{\depth}{{\operatorname{depth}}}
\nc{\dirlim}{{\underset{\rightarrow}{\operatorname{lim}}}}
\nc{\dotbox}{{\overset{\bullet}{\boxtimes}}}
\nc{\dotimes}{{\overset{\bullet}{\otimes}}}
\nc{\Ed}{{\operatorname{Edge}}}
\nc{\ev}{{\operatorname{ev}}}
\nc{\emp}{{\emptyset}}
\nc{\Ext}{{\operatorname{Ext}}}
\nc{\Fac}{{\cal{F}ac}}
\nc{\Fun}{{\operatorname{F}}}
\nc{\FS}{{\cal{FS}}}
\nc{\Hom}{{\operatorname{Hom}}}
\nc{\had}{{{\hat{\mbox{\bf{ad}}}}}}
\nc{\hgt}{{\operatorname{ht}}}
\nc{\Id}{{\operatorname{Id}}}
\nc{\id}{{\operatorname{id}}}
\nc{\Ima}{{\operatorname{Im}}}
\nc{\ind}{{\operatorname{ind}}}
\nc{\Ind}{{\operatorname{Ind}}}
\nc{\infi}{{\operatorname{inf}}}
\nc{\infh}{{\frac{\infty}{2}}}
\nc{\invlim}{{\underset{\leftarrow}{\operatorname{lim}}}}
\nc{\Jac}{{{\cal J}ac}}
\nc{\Ker}{{\operatorname{Ker}}}
\nc{\lcm}{{\operatorname{lcm}}}
\nc{\Locsys}{{{\cal L}ocsys}}
\nc{\Map}{{{\cal M}ap}}
\nc{\modul}{{\operatorname{mod}}}
\nc{\Mor}{{\operatorname{Mor}}}
\nc{\MS}{{\cal{MS}}}
\nc{\Ob}{{\operatorname{Ob}}}
\nc{\opp}{{\operatorname{opp}}}
\nc{\Or}{{{\cal O}r}}
\nc{\Ord}{{{\cal O}rd}}
\nc{\Part}{{{\cal P}art}}
\nc{\PGL}{{\operatorname{PGL}}}
\nc{\Pic}{{\operatorname{Pic}}}
\nc{\pr}{{\operatorname{pr}}}
\nc{\Rep}{{{\cal{R}}ep}}
\nc{\rk}{{\operatorname{rk}}}
\nc{\Sets}{{{\cal{S}}ets}}
\nc{\Sew}{{{\cal{S}}ew}}
\nc{\sgn}{{\operatorname{sgn}}}
\nc{\Sh}{{{\cal S}h}}
\nc{\Sign}{{{\cal S}ign}}
\nc{\Spe}{{\mbox{\bf{Sp}}}}
\nc{\supr}{{\operatorname{sup}}}
\nc{\Supp}{{\operatorname{Supp}}}
\nc{\supp}{{\operatorname{supp}}}
\nc{\Teich}{{{\cal{T}}eich}}
\nc{\tFS}{{\widetilde{\cal{FS}}}}
\nc{\Tor}{{\operatorname{Tor}}}
\nc{\totimes}{{\tilde{\otimes}}}
\nc{\tr}{{\operatorname{tr}}}
\nc{\tRep}{{\widetilde{{\cal R}ep}}}
\nc{\tTeich}{{\widetilde{{\cal T}eich}}}
\nc{\Vect}{{{\cal V}ect}}
\nc{\Ve}{{\operatorname{Vert}}}
\nc{\wt}{{\widetilde}}

\nc{\bo}{{{\bf{0}}}}
\nc{\One}{{{\bf{1}}}}
\nc{\one}{{{\bf{1}}}}

\nc{\BA}{{\Bbb A}}
\nc{\bA}{{\bf{A}}}
\nc{\bbA}{{\overline{\bf{A}}}}
\nc{\ba}{{{\bf{a}}}}
\nc{\baB}{{\overline{B}}}
\nc{\baeta}{{\bar{\eta}}}
\nc{\baJ}{{\bar{J}}}
\nc{\BB}{{\Bbb B}}
\nc{\bB}{{{\bf{B}}}}
\nc{\bc}{{{\bf{c}}}}
\nc{\bC}{{\overline{C}}}
\nc{\BC}{{\Bbb{C}}}
\nc{\bCC}{{\overline{\cal{C}}}}
\nc{\bCM}{{\overline{\cal{M}}}}
\nc{\bD}{{\bar{D}}}
\nc{\BD}{{\overline{D}}}
\nc{\bd}{{{\bf{d}}}}
\nc{\BE}{{\overline{E}}}
\nc{\BF}{{\Bbb{F}}}
\nc{\bF}{{{\bf{F}}}}
\nc{\bg}{{{\bf{g}}}}
\nc{\bG}{{{\bf{G}}}}
\nc{\BG}{{\Bbb G}}
\nc{\bGamma}{{\overline{\Gamma}}}
\nc{\bbH}{{     {{\bf{H}}}_a       }}
\nc{\bH}{{{\bf{H}}}}
\nc{\bi}{{{\bf{i}}}}
\nc{\bI}{{{\bf{I}}}}
\nc{\bj}{{{\bf{j}}}}
\nc{\bK}{{{\bf{K}}}}
\nc{\bL}{{{\bf{L}}}}
\nc{\BL}{{\Bbb{L}}}
\nc{\blambda}{{\bar{\lambda}}}
\nc{\bM}{{{\bf{M}}}}
\nc{\bmu}{{\vec{\mu}}}
\nc{\bN}{{{\bf{N}}}}
\nc{\BN}{{\Bbb{N}}}
\nc{\bnu}{{{\boldmath{${\nu}$}}}}
\nc{\bof}{{{\bf{f}}}}
\nc{\bp}{{{\bf{p}}}}
\nc{\BP}{{\Bbb P}}
\nc{\bP}{{{\bf{P}}}}
\nc{\BPO}{{\overset{\circ}{\BP}}}
\nc{\BQ}{{\Bbb Q}}
\nc{\bQ}{{\bf Q}}
\nc{\bq}{{{\bf{q}}}}
\nc{\BR}{{\Bbb{R}}}
\nc{\bR}{{{\bf{R}}}}
\nc{\br}{{{\bf{r}}}}
\nc{\breta}{{\bar{\eta}}}
\nc{\bs}{{{\bf{s}}}}
\nc{\bS}{{{\bf{S}}}}
\nc{\bt}{{{\bf{t}}}}
\nc{\bU}{{{\bf{U}}}}
\nc{\bV}{{{\bf{V}}}}
\nc{\bu}{{{\bf{u}}}}
\nc{\BUpsilon}{{\bar{\Upsilon}}}
\nc{\bw}{{{\bf{w}}}}
\nc{\bx}{{{\bf{x}}}}
\nc{\bX}{{{\bf{X}}}}
\nc{\BZ}{{\Bbb{Z}}}
\nc{\bz}{{{\bf{z}}}}
\nc{\bZ}{{{\bf{Z}}}}
\nc{\bzero}{{\boldmath{$0$}}}

\nc{\CA}{{\cal A}}
\nc{\CAD}{{\overset{\bullet}{\cal{A}}}}
\nc{\CAO}{{\overset{\circ}{\cal{A}}}}
\nc{\CB}{{\cal B}}
\nc{\CC}{{\cal C}}
\nc{\CalD}{{\cal D}}
\nc{\CE}{{\cal E}}
\nc{\CF}{{\cal F}}
\nc{\CG}{{\cal G}}
\nc{\CH}{{\cal H}}
\nc{\CI}{{\cal I}}
\nc{\CID}{{\overset{\bullet}{\cal{I}}}}
\nc{\CJ}{{\cal J}}
\nc{\CK}{{\cal K}}
\nc{\CL}{{\cal L}}
\nc{\CM}{{\cal M}}
\nc{\CN}{{\cal N}}
\nc{\CO}{{\cal O}}
\nc{\CP}{{\cal P}}
\nc{\CPO}{{\overset{\circ}{\cal{P}}}}
\nc{\CQ}{{\cal Q}}
\nc{\CR}{{\cal R}}
\nc{\CS}{{\cal S}}
\nc{\CT}{{\cal T}}
\nc{\CTD}{{\overset{\bullet}{\cal{T}}}}
\nc{\CTPO}{{\overset{\circ}{\cal{T}\cal{P}}}}
\nc{\CU}{{\cal{U}}}
\nc{\CV}{{\cal V}}
\nc{\CW}{{\cal W}}
\nc{\CX}{{\cal X}}
\nc{\CY}{{\cal Y}}
\nc{\CZ}{{\cal Z}}

\nc{\dCL}{{\overset{\bullet}{\cal{L}}}}
\nc{\dd}{{\operatorname{d}}}
\nc{\ddelta}{{\overset{\bullet}{\delta}}}
\nc{\dfu}{{\overset{\bullet}{\frak{u}}}}
\nc{\dlambda}{{\overset{\bullet}{\lambda}}}
\nc{\DO}{{\overset{\circ}{D}}}
\nc{\dpar}{{\partial}}
\nc{\dS}{{\overset{\bullet}{S}}}
\nc{\dT}{{\overset{\bullet}{T}}}

\nc{\fa}{{\frak{a}}}
\nc{\fA}{{\frak{A}}}
\nc{\fb}{{\frak{b}}}
\nc{\fB}{{\frak{B}}}
\nc{\fC}{{\frak{C}}}
\nc{\fD}{{\frak{D}}}
\nc{\fe}{{\frak{e}}}
\nc{\fE}{{\frak{E}}}
\nc{\fF}{{\frak{F}}}
\nc{\ff}{{\frak{f}}}
\nc{\fg}{{\frak{g}}}
\nc{\fG}{{\frak{G}}}
\nc{\fH}{{\frak{H}}}
\nc{\fii}{{\frak{i}}}
\nc{\fj}{{\frak{j}}}
\nc{\fl}{{\frak{l}}}
\nc{\fL}{{\frak{L}}}
\nc{\fM}{{\frak{M}}}
\nc{\fN}{{\frak{N}}}
\nc{\fn}{{\frak{n}}}
\nc{\fp}{{\frak{p}}}
\nc{\fq}{{\frak{q}}}
\nc{\fQ}{{\frak{Q}}}
\nc{\fs}{{\frak{s}}}
\nc{\ft}{{\frak{t}}}
\nc{\fu}{{\frak{u}}}
\nc{\fU}{{\frak{U}}}
\nc{\fv}{{\frak{v}}}
\nc{\fV}{{\frak{V}}}
\nc{\fW}{{\frak{W}}}
\nc{\fx}{{\frak{x}}}
\nc{\fy}{{\frak{y}}}
\nc{\fZ}{{\frak{Z}}}

\nc{\hCH}{{\hat{\cal{H}}}}
\nc{\hCI}{{\hat{\cal{I}}}}
\nc{\hfC}{{\hat{\frak{C}}}}
\nc{\hfg}{{\hat{\frak{g}}}}
\nc{\hL}{{\hat{L}}}
\nc{\HO}{{\overset{\circ}{H}}}
\nc{\hpsi}{{\hat{\psi}}}
\nc{\hx}{{\hat{x}}}

\nc{\jo}{{\overset{\circ}{j}}}

\nc{\phid}{{\overset{\bullet}{\phi}}}

\nc{\tA}{{\tilde{A}}}
\nc{\ta}{{\tilde{a}}}
\nc{\tB}{{\tilde{B}}}
\nc{\tb}{{\tilde{b}}}
\nc{\tBP}{{\tilde{\BP}}}
\nc{\tC}{{\tilde{C}}}
\nc{\tc}{{\tilde{c}}}
\nc{\tCA}{{\tilde{\cal{A}}}}
\nc{\tCC}{{\tilde{\cal{C}}}}
\nc{\tCH}{{\tilde{\cal{H}}}}
\nc{\tCI}{{\tilde{\cal{I}}}}
\nc{\tCO}{{\tilde{\cal{O}}}}
\nc{\tCP}{{\tilde{\cal{P}}}}
\nc{\tCT}{{\tilde{\cal{T}}}}
\nc{\tD}{{\tilde{D}}}
\nc{\tDelta}{{\tilde{\Delta}}}
\nc{\tE}{{\tilde E}}
\nc{\tF}{{\tilde F}}
\nc{\tfD}{{\tilde{\frak{D}}}}
\nc{\tfF}{{\tilde{\frak{F}}}}
\nc{\tff}{{\tilde{\frak{f}}}}
\nc{\tfu}{{\tilde{\frak{u}}}}
\nc{\tJ}{{\tilde{J}}}
\nc{\tj}{{\tilde{j}}}
\nc{\tK}{{\tilde K}}
\nc{\tL}{{\tilde{L}}}
\nc{\tM}{{\tilde{M}}}
\nc{\tP}{{\tilde{P}}}
\nc{\tPhi}{{\tilde{\Phi}}}
\nc{\tpi}{\tilde{\pi}}
\nc{\TPO}{{\overset{\circ}{T\BP}}}
\nc{\tR}{{\tilde{R}}}
\nc{\tS}{{\tilde S}}
\nc{\tT}{{\tilde{T}}}
\nc{\ttau}{{\tilde{\tau}}}
\nc{\ttheta}{{\tilde{\theta}}}
\nc{\tU}{{\tilde{U}}}
\nc{\tUpsilon}{{\tilde{\Upsilon}}}
\nc{\tW}{{\tilde W}}
\nc{\ty}{{\tilde y}}
\nc{\tY}{{\tilde Y}}
\nc{\txi}{{\tilde{\xi}}}

\nc{\UD}{{\overset{\bullet}{U}}}
\nc{\UO}{{\overset{\circ}{U}}}

\nc{\vA}{{\vec{A}}}
\nc{\valpha}{{\vec{\alpha}}}
\nc{\vbeta}{{\vec{\beta}}}
\nc{\vc}{{\vec{c}}}
\nc{\vD}{{\vec{D}}}
\nc{\vd}{{\vec{d}}}
\nc{\vgamma}{{\vec{\gamma}}}
\nc{\vK}{{\vec{K}}}
\nc{\vlambda}{{\vec{\lambda}}}
\nc{\vmu}{{\vec{\mu}}}
\nc{\vnu}{{\vec{\nu}}}
\nc{\vo}{{\vec{0}}}
\nc{\vu}{{\vec{u}}}
\nc{\vx}{{\vec{x}}}
\nc{\vy}{\vec{y}}
\nc{\vzero}{\vec{0}}

\nc{\XO}{{\overset{\circ}{X}}}

\nc{\ya}{{\operatorname{aj}}}

\nc{\nen}{\newenvironment}
\nc{\ol}{\overline}
\nc{\ul}{\underline}
\nc{\ra}{\rightarrow}
\nc{\lra}{\longrightarrow}
\nc{\Lra}{\Longrightarrow}
\nc{\lla}{\longleftarrow}
\nc{\Llra}{\Longleftrightarrow}
\nc{\hra}{\hookrightarrow}
\nc{\iso}{\overset{\sim}{\lra}}
\nc{\rlh}{\rightleftharpoons}

\nc{\IC}{{\cal{IC}}}
\nc{\ic}{{\bf{IC}}}
\nc{\PS}{{\cal{PS}}}
\nc{\oCG}{{\overline{\cal G}}}
\nc{\osG}{{\overline{\sf G}}}
\nc{\oCQ}{{\overline{\cal Q}}}
\nc{\oCZ}{{\overline{\cal Z}}}
\nc{\dK}{{\overset{\bullet}{\cal K}}{}}
\nc{\dQ}{{\overset{\bullet}{\cal Q}}{}}
\nc{\dsQ}{{\overset{\bullet}{{\sf Q}}}{}}
\nc{\dZ}{{\overset{\bullet}{\cal Z}}{}}
\nc{\dCS}{{\overset{\bullet}{\cal S}}{}}
\nc{\dbQ}{{\overset{\bullet}{\bf Q}}{}}
\nc{\ddZ}{{\ddot{\cal Z}}{}}
\nc{\ddQ}{{\ddot{\cal Q}}{}}
\nc{\oZ}{{\overset{\circ}{\cal Z}}{}}
\nc{\oF}{{\overset{\circ}{\cal F}}{}}
\nc{\dP}{{\overset{\bullet}{\cal P}}{}}
\nc{\oP}{{\overset{\circ}{\cal P}}{}}
\nc{\oQ}{{\overset{\circ}{\cal Q}}{}}
\nc{\obp}{{\overset{\circ}{{\bf p}}}{}}
\nc{\tbj}{{\tilde{\bf j}}{}}
\nc{\tbp}{{\tilde{\bf p}}{}}
\nc{\tfC}{{\widetilde{\frak C}}{}}
\nc{\tfE}{{\widetilde{\frak E}}{}}
\nc{\tfj}{{\widetilde{\frak j}}{}}
\nc{\tmk}{{\widetilde M^a_K}}
\nc{\tbQ}{{\widetilde{\bf Q}}{}}
\nc{\hCQ}{{\widehat{\cal Q}}{}}
\nc{\tfQ}{{\widetilde{\frak Q}}{}}
\nc{\tfp}{{\widetilde{\frak p}}{}}
\nc{\ofQ}{{\overset{\circ}{{\frak Q}}}{}}
\nc{\osQ}{{\overset{\circ}{{\sf Q}}}{}}
\nc{\tGQ}{{\widetilde{\cal{GQ}}}{}}
\nc{\tCS}{{\widetilde{\cal S}}{}}
\nc{\oGQ}{{\overset{\circ}{\cal{GQ}}}{}}
\nc{\ooGQ}{{\overset{\circ\circ}{\cal{GQ}}}{}}
\nc{\oGZ}{{\overset{\circ}{\cal{GZ}}}{}}
\nc{\tGZ}{{\widetilde{\cal{GZ}}}{}}
\nc{\ufM}{{\underline{\frak M}}}
\nc{\ufQ}{{\underline{\frak Q}}}
\nc{\usQ}{{\underline{\sf Q}}}
\nc{\sG}{{\sf G}}
\nc{\sQ}{{\sf Q}}
\nc{\sM}{{\sf M}}
\nc{\usM}{{\underline{\sf M}}}
\nc{\Ue}{{U_\varepsilon}}
\nc{\Upe}{{\Upsilon_\varepsilon}}
\nc{\crho}{{\check{\rho}}}
\nc{\ctheta}{{\check{\theta}}}
\nc{\cR}{{\check{\cal R}}{}}
\nc{\del}{{\partial}}
\nc{\bzw}{{{\Bbb Z}[v,v^{-1}]}}
\nc{\hK}{{\hat K}}
\nc{\hP}{{\hat{\cal P}}}
\nc{\utP}{{\underline{\tilde P}{}}}

\nc{\Thm}[1]{Theorem~\ref{#1}}
\nc{\Prop}[1]{Proposition~\ref{#1}}
\nc{\Lem}[1]{Lemma~\ref{#1}}
\nc{\Cor}[1]{Corollary~\ref{#1}}
\nc{\Conj}[1]{Conjecture~\ref{#1}}
\nc{\Claim}[1]{Claim~\ref{#1}}
\nc{\Defn}[1]{Definition~\ref{#1}}
\nc{\Exa}[1]{Example~\ref{#1}}
\nc{\Rem}[1]{Remark~\ref{#1}}
\nc{\Note}[1]{Note~\ref{#1}}

\nen{thm}[1]{\label{#1}{\bf Theorem.\ } \em}{}
\nen{prop}[1]{\label{#1}{\bf Proposition.\ } \em}{}
\nen{lem}[1]{\label{#1}{\bf Lemma.\ } \em}{}
\nen{cor}[1]{\label{#1}{\bf Corollary.\ } \em}{}
\nen{conj}[1]{\label{#1}{\bf Conjecture.\ } \em}{}
\nen{claim}[1]{\label{#1}{\bf Claim.\ } \em}{}

\nen{defn}[1]{\label{#1}{\bf Definition.\ } }{}
\nen{exa}[1]{\label{#1}{\bf Example.\ } }{}

\nen{rem}[1]{\label{#1}{\em Remark.\ } }{}
\nen{note}[1]{\label{#1}{\em Note.\ } }{}
\nen{exer}[1]{\label{#1}{\em Exercise.\ } }{}

\begin{document}

%%%%%%%%%%%%%%%%%%%%%%%%%%%%%%%%%%MINE%%%%%%%%%%%%%%%%%%%%%%%%%%%%%%%%
%%%%%%%%%%%%%%%%%%%%%%%%%%%%%%%%%%%%%%%%%%%%%%%%FOUNDATIONS
%       1. Graphical Features
%       2. ARROWS:
%       3. VARIOUS SPECIAL SIGNS
%       4. LARGE OPERATORS
%       5. ENDINGS
%       6. Special ``Words'' in TEX:
%       7. FONT
%%%%%SYMBOLS:
%REDEFINE-the same format as define!: \redefine\gg{frak g}
%\predefine\greater{\gg}
%

%1. SHORTHAND for Graphical Features:
                                %Indentation &
\nc{\nn}{\newline}
\nc{\nnn}{\newpage}
\nc{\noi}{\noindent}
\nc{\nop}{\noindent {\bf Proof.} }
                                %BREAKS
\nc{\sbr}{\smallpagebreak}
\nc{\mbr}{\medpagebreak}
\nc{\bbr}{\bigpagebreak}

%2. ARROWS:

%%%\nc{\ra}{\rightarrow}
\nc{\raa}{\longrightarrow}
\nc{\lala}{\leftarrow}
\nc{\laa}{\longleftarrow}
\nc{\lrax}{\leftrightarrow}     %xxxx i.e. x-added since M uses this already

\nc{\Ra}{\Rightarrow}         %Implication.
\nc{\LRa}{\Rightarrow}        %Equivalence.

\nc{\inj}{\hookrightarrow}    %Injective map-right.
\nc{\injj}{\hookleftarrow}    %Injective map.
\nc{\sur}{\twoheadrightarrow} %Surjective map-right.
\nc{\surr}{\twoheadleftarrow} %Surjective map.
\nc{\mm}{\mapsto}             %Map on elements.
\nc{\va}{\uparrow}              %Up-arrow

%3. VARIOUS SPECIAL SIGNS

\nc{\bb}{\underset}           %Write Belllow
\nc{\aax}{\overset}     %%%%xxxxx          %and Above

\nc{\bsx}{\backslash}           %%%xxx
\nc{\bss}{\backslash}           %%%xxx
\nc{\barr}{\overline}         %Long bar accent
\nc{\sss}{\S}               %Section

\nc{\sub}{\subseteq}          %Inclusions
\nc{\suppp}{\supseteq}          %%%xxx

\nc{\ti}{\tilde}              %Tilde
\nc{\tii}{\widetilde}         %Tilde-wide
\nc{\ch}{\check}              %Check

\nc{\tim}{\times}             %Times
\nc{\btim}{\boxtimes}
\nc{\ten}{\otimes}            %Tensoring
\nc{\bten}{\boxtimes}         %Tensoring: outer
\nc{\pl}{\oplus}              %Direct sum
\nc{\con}{ @>\cong>> }  %Isomorphism with a right arrow
\nc{\conn}{     @<\cong<<  }    %Isomormphism with a left arrow

\nc{\half}{ \frac{1}{2} }     %%%xxxxxxxx  %Half
%%%\nc{\half12}{ \frac{1}{2} }       %Half

                                %DOTS:
\nc{\ci}{\circ}               %Circle dot
\nc{\cdx}{\cdot}             %%%xxx             %Dot(s)
\nc{\cdd}{\cdot}             %%%xxx             %Dot(s)
\nc{\cddd}{\cdot\cdot\cdot}

%4. LARGE OPERATORS                                %such as \sum

\nc{\cupp}{\bigcup}             %cup,cap
\nc{\capp}{\bigcap}
\nc{\tenn}{\bigotimes}          %times,plus
\nc{\pll}{\bigoplus}

\nc{\pii}{\prod}                %product
\nc{\ppii}{\bigprod}            %Big product

\nc{\cci}{\bigcoprod}
\nc{\wwe}{\bigwedge}            %wedge
\nc{\cce}{\bigcoprod}           %cowedge

%5. ENDINGS
\nc{\pp}{\endproclaim}        %Endproclaim
\nc{\hh}{\endheading}         %Endheading

                                %SHORTHAND for
%6. SPECIAL ``Words'' in TEX:
                                        %Special Symbols
\nc{\df}{ \overset{ \text{def}}= }
\nc{\inv}{ {}^{-1}      }
%%%\nc{\emp}{   \emptyset}      %Emptyset
\nc{\we}{\wedge}
\nc{\wee}{{     \overset{2}\wedge       }}
                                        %Particular Varieties
%%%%%%%%%%%%\redefine{\P}{\Bbb P}
\nc{\ppp}{{ \Bbb P^1 }}            %%%xxxx  %P1
%%%%%%\nc{\p1}{{ \P^1 }}              %P1
%%%%%%\nc{\a1}{{\A^1}}         %%%xxxx     %A1
\nc{\aaa}{{\Bbb A^1}}              %A1
                                        %Particular Fields:
\nc{\qlb}{ \barr{\Q_l} }      %Q-el-bar
\nc{\ffq}{ {\F_q} }           %Ef-q (finite field)
                                        %Actual ``Words''
\nc{\Spec}{{\text{Spec}\ {}     }}
\nc{\aand}{ \ \text{and}\ }
\nc{\hk}{       \text{hyperk\"ahler}    }
%%%%%%%\nc{\Hom}{\text{Hom}}
%%%%%%%\nc{\Ext}{\text{Ext}}
%%%%%%%\nc{\Ker}{\text{Ker}}
%%%%%%%\nc{\Coker}{\text{Coker}}
%%%%%%%\redefine{\Im}{{ \text{Im} }}
%%%%%%%\nc{\codim}{\text{codim}}
\nc{\rank}{{\ \text{rank}\ }}
                                        %Various
                                        %Relative Products, abelianizationa
                                        %and affinization
\nc{\timB} {{   \underset{B}\tim                }}
\nc{\timP}{{    \underset{P}\tim                }}
\nc{\timQ}{{    \underset{Q}\tim                }}
\nc{\ab}{       ^{ab}   }
\nc{\af}{       ^{aff}  }

%Already EXISTING shorthand ``words'':
%\nc{\mod}{\text{mod{\}}}
%\nc{\dim}{\text{dim}}
%\nc{\max}{\text{max}}
%\nc{\and
%\nc{\or}{\text{or}}

%7. FONT

%FONT: SCRIPT,      Cal by a double letter \??
\nc{\AAA}{\cal A}
\nc{\BBB}{\cal B}       %%%%xxxxx
%%%%\nc{\CC}{\cal C}
\nc{\DD}{\cal D}
\nc{\EE}{\cal E}
\nc{\FF}{\cal F}
\nc{\GG}{\cal G}
\nc{\HH}{\cal H}
\nc{\II}{\cal I}
\nc{\JJ}{\cal J}
\nc{\KK}{\cal K}
\nc{\LL}{\cal L}
\nc{\MM}{\cal M}
\nc{\NN}{\cal N}
\nc{\OO}{\cal O}
\nc{\PP}{\cal P}
\nc{\QQ}{\cal Q}
\nc{\RR}{\cal R}
\nc{\SSS}{\cal S}
%%%\redefine{\SS}{\cal S}
\nc{\TT}{\cal T}
\nc{\UU}{\cal U}
\nc{\VV}{\cal V}
\nc{\ZZ}{\cal Z}
\nc{\XX}{\cal X}
\nc{\YY}{\cal Y}

% FONT: BLACKBOARD         by a single letter \?
%%%%\redefine{\A}{\Bbb A }
%%%%\redefine{\B}{\Bbb B}
%%%%\redefine{\Cs}{\Bbb C^*}
%%%%\redefine{\C}{\Bbb C}
%%%%\redefine{\D}{\Bbb D}
%%%%\redefine{\E}{\Bbb E}
%%%%\redefine{\F}{\Bbb F}
\nc{\A}{\Bbb A }
\nc{\cs}{\Bbb C^*}
\nc{\ccs}{ \Bbb C^*}
\nc{\cc}{\Bbb C}
\nc{\f}{\Bbb F}
\nc{\g}{\Bbb G}
\nc{\h}{\Bbb H}
\nc{\I}{\Bbb I}
\nc{\J}{\Bbb J}
\nc{\K}{\Bbb K}
%\redefine{\L}{\Bbb L}
\nc{\M}{\Bbb M}
\nc{\N}{\Bbb N}
%\redefine{\O}{\Bbb O}
\nc{\p}{\Bbb P}
\nc{\Q}{\Bbb Q}
\nc{\R}{\Bbb R}
\nc{\s}{\Bbb S}
\nc{\T}{\Bbb T}
\nc{\U}{\Bbb U}
\nc{\V}{\Bbb V}
\nc{\Z}{\Bbb Z}
\nc{\X}{\Bbb X}
\nc{\Y}{\Bbb Y}

%FONT: FRAKTUR        by \f?
%\nc{\fA}{\frak A}
%\nc{\fB}{\frak B}
%\nc{\fC}{\frak C}
%\nc{\fD}{\frak D}
%\nc{\fE}{\frak E}
%\nc{\fF}{\frak F}
%\nc{\fG}{\frak G}
%\nc{\fH}{\frak H}
\nc{\fI}{\frak I}
\nc{\fJ}{\frak J}
\nc{\fK}{\frak K}
%\nc{\fL}{\frak L}
%\nc{\fM}{\frak M}
%\nc{\fN}{\frak N}
\nc{\fO}{\frak O}
\nc{\fP}{\frak P}
%\nc{\fQ}{\frak Q}
\nc{\fR}{\frak R}
\nc{\fS}{\frak S}
\nc{\fT}{\frak T}
%\nc{\fU}{{\frak U}}
%\nc{\fV}{\frak V}
%\nc{\fZ}{\frak Z}
\nc{\fX}{\frak X}
\nc{\fY}{\frak Y}
%\nc{\fa}{\frak a}
%\nc{\fb}{\frak b}
\nc{\fc}{\frak c}
\nc{\fd}{\frak d}
%\nc{\fe}{\frak e}
%\nc{\ff}{\frak f}
%\nc{\fg}{\frak g}
\nc{\fh}{\frak h}
%\nc{\fi}{\frak i}
%\nc{\fj}{\frak j}
\nc{\fk}{\frak k}
%\nc{\fl}{\frak l}
\nc{\fm}{\frak m}
%\nc{\fn}{\frak n}
\nc{\fo}{\frak o}
%\nc{\fp}{\frak p}
%\nc{\fq}{\frak q}
\nc{\fr}{\frak r}
%\nc{\fs}{\frak s}
%\nc{\ft}{\frak t}
%\nc{\fu}{\frak u}
%\nc{\fv}{\frak v}
\nc{\fz}{\frak z}
%\nc{\fx}{\frak x}
%\nc{\fy}{\frak y}

%FONT: GREEK
\nc{\al}{\alpha }
\nc{\be}{\beta }
\nc{\ga}{\gamma }
\nc{\de}{\delta }
%\nc{\del}{\partial }
\nc{\ep}{\varepsilon }
\nc{\vap}{\epsilon }

\nc{\ze}{\zeta }
\nc{\et}{\eta }
\nc{\th}{\theta}
%\nc{\theta}{\theta }   %%%%xxxxx%%%%
\nc{\vth}{\vartheta }

\nc{\io}{\iota }
\nc{\ka}{\kappa }
\nc{\la}{\lambda }
%mu
%nu
%xi
%pi
%rho
\nc{\vrho}{\varrho}
\nc{\si}{\sigma }
%\nc{\tau}{\tau }
\nc{\ups}{\upsilon }
%phi
\nc{\vphi}{\varphi }
%chi
%psi
\nc{\om}{\omega }

\nc{\Ga}{\Gamma }
\nc{\De}{\Delta }
\nc{\Th}{\Theta }
\nc{\La}{\Lambda }
%Xi
%Pi
\nc{\Si}{\Sigma }
\nc{\Ups}{\Upsilon }
%Phi
%Psi
\nc{\Om}{\Omega }

%%%%%%%%%%%%%%%%%%%%%%%%MINEendOF%%%%%%%%%%%%%%%%%%%%%%%%%%%%%%%%%%

%HOME
%LOCAL DEFINITIONS
\nc{\zp}{{\overset{\bullet}{\cal Z}}{}}
\nc{\zc}{{\overset{\circ}{\cal Z}}{}}
\nc{\qp}{{\overset{\bullet}{\cal Q }}{}}
\nc{\qc}{{\overset{\circ}{\cal Q}}{}}
\nc{\nii}{ ^{n\cdd i} }
\nc{\yy}{\infty}

\nc{\cz}{\cc[z]}
\nc{\czn}{{ \cc_{\le n}[z] }}
\nc{\pv}{{\p(V)}}
\nc{\qv}{Q(V)}
\nc{\qcv}{{\overset{\circ}Q(V)}}
\nc{\ppi}{\p(I)}
\nc{\pk}{\p(K)}
\nc{\zv}{\ZZ(V)}
\nc{\zcv}{\zc(V)}

\nc{\ii}{{i\in I}}
\nc{\all}{{ ^{(\alpha)} }}
\nc{\bee}{{ ^{(\beta)} }}
\nc{\gaa}{{ ^{(\gamma)} }}

%  top matter
\title[]{Semiinfinite flags. II. Local and global Intersection Cohomology of
Quasimaps' spaces.}
\author{Boris Feigin}
\address{Landau Institute of Theoretical Physics, ul. Kosygina 2, Moscow,
Russia}
\email{feigin@@landau.ac.ru}
\author{Michael Finkelberg}
\address{Independent Moscow University, Bolshoj Vlasjevskij pereulok, dom 11,
Moscow 121002 Russia}
\email{fnklberg@@mccme.ru}
\author{Alexander Kuznetsov}
\address{Independent Moscow University, Bolshoj Vlasjevskij pereulok, dom 11,
Moscow 121002 Russia}
\email{sasha@@ium.ips.ras.ru}
\author{Ivan Mirkovi\'c}
\address{Dept. of Mathematics and Statistics, University of Massachusetts
at Amherst, Amherst MA 01003-4515, USA}
\email{mirkovic@@math.umass.edu}
\thanks{M.F. and I.M. were partially supported by the NSF grant DMS 97-29992.}
%\date{June 1997}
\maketitle

\centerline{\em To Dmitri Borisovich Fuchs on the occasion of his 60th
birthday}

\section{Introduction}

\subsection{}
This paper is a sequel to ~\cite{fm}. We will make a free use of notations,
conventions and results of {\em loc. cit.}

One of the main results of the present work is a computation of local $\IC$
stalks of the Schubert strata closures in the spaces $\CZ^\alpha$. We prove
that the generating functions of these stalks are given by the
{\em generic} (or {\em periodic) Kazhdan-Lusztig polynomials}, see the
Theorem ~\ref{main}. We understand that this result was known to G.Lusztig
for a long time, cf. ~\cite{l2} ~\S11. His proof was never published though,
and as far as we understand, it differs from ours: for example we never
managed to find a direct geometric proof of the key property ~\cite{l1}
~11.1.(iv) of Lusztig's $R$-polynomials.

Our proof uses the standard convolution technique. The only nonstandard
feature is the check of pointwise purity of the $\IC$ sheaves involved
(Theorem ~\ref{finkel}).
Usually one proceeds by finding global transversal slices. We were not able
to find the good slices, and instead reduced the proof to the purity
properties of $\IC$ sheaves on the affine Grassmannian.

\subsection{}
The central result of this work is the geometric construction of the
universal enveloping algebra $U(\fn_+^L)$ of the nilpotent subalgebra of
the Langlands dual Lie algebra $\fg^L$ (Theorem ~\ref{!}). This construction
occupies the section 2. The geometric incarnation $\CA$ of $U(\fn_+^L)$
naturally acts on the global intersection cohomology
$\oplus_{\alpha\in\BN[I]}H^\bullet(\CQ^\alpha,\IC(\CQ^\alpha))$
of all Quasimaps' spaces, and this action extends to the geometrically defined
$\fg^L$-action (section 4).

In case $\bG=SL_n$ such action was constructed in ~\cite{fk}, and the present
paper grew out of attempts to generalize the results of {\em loc. cit.} to
the case of arbitrary simple $\bG$. 
For $\fg=\frak{sl}_n$ the geometric construction of
$U(\fn_+^L)$ given in {\em loc. cit.} used the Laumon resolution
$\CQ_L^\alpha$ of the Quasimaps' spaces $\CQ^\alpha$, and provided
$U(\fn_+^L)$ with the geometrically defined Poincar\'e-Birkhoff-Witt basis.
The present construction also provides $U(\fn_+^L)$ with a special basis
numbered by the irreducible components of 
the semiinfinite orbits' intersections
in the affine Grassmannian. We want to stress that the two bases are
{\em different}: the latter one looks more like a canonical basis of ~\cite{l}.

Another advantage of the new construction of $U(\fn_+^L)$ is its local nature.
It allows one to define a $\fg^L$-action on the global cohomology of
Quasimaps' spaces with coefficients in sheaves more general than just
$\IC(\CQ^\alpha)$. Namely, given a perverse sheaf $\CF\in\CP(\oCG_\eta,\bI)$
on the affine Grassmannian we defined in ~\cite{fm} its {\em convolution}
$\bc^\alpha_\CQ(\CF)$ --- a perverse sheaf on $\CQ^{\eta+\alpha}$.
In section 7 we construct the $\fg^L$-action on
$\oplus_{\alpha\in Y}H^\bullet(\CQ^{\eta+\alpha},\bc^\alpha_\CQ(\CF))$.
In case $\CF$ is $\bG[[z]]$-equivariant we conjecture that the resulting
$\fg^L$-module is {\em tilting}. As in ~\cite{fk}, this conjecture is
motivated by an analogy with the semiinfinite cohomology of quantum groups.

Recall that for a $\bG[[z]]$-equivariant sheaf $\CF$ on the affine Grassmannian
the action of $\fg^L$ on its global cohomology was constructed in ~\cite{mv}.
The relation between the various $\fg^L$-actions on global cohomology will
be discussed in a separate paper.

\subsection{}
Let us list the other points of interest in this paper.

In section 3 we compute the stalks of $\IC(\CZ^\alpha)$
(Theorem ~\ref{simple} and Corollary ~\ref{berezin}).
Their generating functions are expressed in terms
of Lusztig's $q$-analogue $\CK^\alpha(t)$
of Kostant partition function. In case
$\bG=SL_n$ this result was proved in ~\cite{ku} using the Laumon resolution
of the Quasimaps' space $\CQ^\alpha$. The proof in general case uses the
{\em Beilinson-Drinfeld} incarnation of $\CZ^\alpha$ (see ~\cite{fm} ~\S6).
Formally, Theorem ~\ref{simple} is just a particular case of the Theorem
~\ref{main} computing $\IC$-stalks of the general Schubert strata closures in
$\CZ^\alpha$. But the argument goes the other way around: we deduce ~\ref{main}
as a rather formal corollary of ~\ref{simple}. It is well known that
the generating function of $\IC$ stalks of $\bG[[z]]$-orbits' closures
in the affine Grassmannian is also given by $\CK^\alpha(t)$ in the stable
range (see ~\cite{lus}). This coincidence is explained in section 3: though the
local singularities of $\CZ^\alpha$ and $\oCG_\eta$ {\em are different},
their ``skeleta'' (intersections of semiinfinite orbits in $\oCG_\eta$,
and the {\em central fiber} in $\CZ^\alpha$) {\em are the same}.

In section 4 we prove that
$\oplus_{\alpha\in\BN[I]}H^\bullet(\CQ^\alpha,\IC(\CQ^\alpha))$
carries a natural pure Tate Hodge structure and compute its generating
function (Theorems ~\ref{ostrik}, ~\ref{feigin}).

In section 5 we prove that the adjacency order on the set of Schubert
strata in $\CZ$ is equivalent to Lusztig's order
(see ~\cite{l1}) on the set of {\em alcoves}. The proof uses a map
$\pi:\ \CQ^\alpha_K\lra\CQ^\alpha$ from Kontsevich's space of {\em stable maps}
to Drinfeld's Quasimaps' space. This map is constructed in the Appendix.
The construction is just an application of Givental's ``Main Lemma''
(see ~\cite{g}). His proof of the ``Main Lemma'' has not satisfied all of its
readers, so we include the complete proof into the Appendix.

In section 7 we collect various conjectures on the structure of 
$\fg^L$-modules of geometric origin. For a mysterious reason these conjectures
involve tilting modules --- either over $\fg^L$ itself, or over the related
quantum group. The conjecture ~\ref{roman} was recently proved in ~\cite{fkm}.
The conjecture ~\ref{denis} may be viewed as a description of an
``automorphic sheaf'' on the moduli space of $\bG$-torsors corresponding
to the trivial $\bG^L$-local system on $\BP^1$. 
Finally, in ~\ref{romka} we propose a direct geometric construction of
the $\fn^L_\pm$-action on $H^\bullet(\CG,\CF)$ for a $\bG[[z]]$-equivariant
perverse sheaf $\CF$ on the affine Grassmannian $\CG$. Surprisingly enough,
no direct construction of $\fn^L_\pm$-action has been found so far
(cf. ~\cite{mv} for the direct construction of the 
action of the dual Cartan ${\fh}^L\subset\fg^L$).

\subsection{}
We are deeply grateful to V.Lunts and L.Positselsky whose explanations helped
us at the moments of despair. We are very much obliged to R.Bezrukavnikov
who found a serious gap in the original version of 
geometric construction of $U(\fn^L_+)$.
It is a pleasure to thank S.Arkhipov, D.Gaitsgory, and V.Ostrik for the
inspiring discussions around the tilting conjectures in section 7. We are
very much obliged to D.Gaitsgory for the careful reading of this paper,
and for suggesting a lot of drastic simplifications of the arguments and proofs
of the conjectures therein. They will all appear in his forthcoming 
publication.

\section{Local Ext-algebra $\CA$}

\subsection{}
Throughout this paper $C$ will denote a genus zero curve $\BP^1$ with the
marked points $0,\infty$. The complement $C-\infty$ is the affine line $\BA^1$.
The {\em twisting} map $\sigma_{\beta,\gamma}:\ \CQ^\beta\times C^\gamma
\lra\CQ^{\beta+\gamma}$ defined in ~\cite{fm} ~3.4.1 restricts to the map
$\varsigma_{\beta,\gamma}:\ \CZ^\beta\times\BA^\gamma\lra\CZ^{\beta+\gamma}$.
We denote its image by $\del_\gamma\CZ^{\beta+\gamma}$. The map
$\varsigma_{\beta,\gamma}:\ \CZ^\beta\times\BA^\gamma\lra
\del_\gamma\CZ^{\beta+\gamma}$ is finite. Moreover, in case $\beta=0$,
the space $\CZ^\beta$ is just a point, and the map $\varsigma_{0,\gamma}:\
\BA^\gamma\lra\del_\gamma\CZ^\gamma$ is an isomorphism. We will identify
$\del_\gamma\CZ^\gamma$ with $\BA^\gamma$ via this map.

{\bf Definition.} $\CA_\alpha:=
\Ext^\bullet_{\CZ^\alpha}(\IC(\del_\alpha\CZ^\alpha),
\IC(\CZ^\alpha));\ \CA:=\oplus_{\alpha\in\BN[I]}\CA_\alpha$.

Here the $\Ext^\bullet_{\CZ^\alpha}$
is taken in the constructible derived category on
$\CZ^\alpha$. {\em A priori} $\CA_\alpha$ is a graded vector space.
In this section we will
show that it is concentrated in the degree $|\alpha|$, and we will define
a structure of cocommutative $\BN[I]$-graded bialgebra on $\CA$.

\subsection{}
\label{Weil}
To unburden the notations we will denote the $\IC$ sheaf $\IC(\CZ^\alpha)$
by $\IC^\alpha$ (see ~\cite{fm} 10.7.1).
We denote the closed embedding of $\BA^\alpha=\del_\alpha\CZ^\alpha$ into
$\CZ^\alpha$ by $\fs_\alpha$. Then $\CA_\alpha
=\Ext^\bullet_{\CZ^\alpha}(\IC(\del_\alpha\CZ^\alpha),\IC^\alpha)=
\Ext^\bullet_{\del_\alpha\CZ^\alpha}
(\IC(\del_\alpha\CZ^\alpha),\fs_\alpha^!\IC^\alpha)=
H^\bullet(\BA^\alpha,\fs_\alpha^!\IC^\alpha)[-|\alpha|]$
since $\IC(\del_\alpha\CZ^\alpha)=\ul\BC[|\alpha|]$.

{\bf Theorem.} a) $H^\bullet(\BA^\alpha,\fs_\alpha^!\IC^\alpha)$
is concentrated in degree 0;

b) $\dim H^0(\BA^\alpha,\fs_\alpha^!\IC^\alpha)=\CK(\alpha)$
where $\CK(\alpha)$ stands for the Kostant partition function.

The proof of the Theorem will occupy the subsections ~\ref{raz}--\ref{dva}.

\subsection{}
\label{raz}
We denote the closed embedding of the origin 0 into $\BA^\alpha$ by
$\iota_\alpha$. Since $\fs_\alpha^!\IC^\alpha$
is constructible with respect
to the diagonal stratification of $\BA^\alpha$, we have the canonical
isomorphism
$H^\bullet(\BA^\alpha,\fs_\alpha^!\IC^\alpha)=
\iota_\alpha^*\fs_\alpha^!\IC^\alpha$.

Recall the projection $\pi_\alpha:\ \CZ^\alpha\lra\BA^\alpha$
(see ~\cite{fm} ~7.3). We denote the {\em central fiber} $\pi_\alpha^{-1}(0)$
by $\CF^\alpha$, and we denote its closed embedding into $\CZ^\alpha$ by
$\iota_\alpha$.

The Cartan group $\bH$ acts on $\CZ^\alpha$ contracting it to the fixed point
set $(\CZ^\alpha)^\bH=\del_\alpha\CZ^\alpha$. The projection $\pi_\alpha$
is $\bH$-equivariant. Hence the canonical morphism
$\fs_\alpha^!\IC^\alpha\lra\pi_{\alpha!}\IC^\alpha$ of sheaves
on $\BA^\alpha$ is an isomorphism. By the proper base change we have the
canonical isomorphism
$\iota_\alpha^*\pi_{\alpha!}\IC^\alpha=\pi_{\alpha!}\iota_\alpha^*\IC^\alpha=
H^\bullet_c(\CF^\alpha,\iota_\alpha^*\IC^\alpha)$.

Combining all the above isomorphisms, we conclude that
$H^\bullet(\BA^\alpha,\fs_\alpha^!\IC^\alpha)=
H^\bullet_c(\CF^\alpha,\iota_\alpha^*\IC^\alpha)$.

\subsection{}
\label{leq0}
In this subsection we prove that
$H^\bullet_c(\CF^\alpha,\iota_\alpha^*\IC^\alpha)$ is concentrated in
nonpositive degrees. To this end we study the intersection of the
{\em fine stratification} of $\CZ^\alpha$ (see ~\cite{fm} ~8.4.1) with the
central fiber $\CF^\alpha$. Recall the description of the central fiber
given in {\em loc. cit.} ~6.4.1, ~6.4.2 in terms of semiinfinite orbits in
the affine Grassmannian: $\CF^\alpha=\ol{T}_{-\alpha}\cap S_0$.

\subsubsection{Lemma} The intersection of $\CF^\alpha$ with a fine stratum
$\oZ^\gamma\times(\BC^*)^{\beta-\gamma}_\Gamma,\ \gamma\leq\beta\leq\alpha,
\Gamma\in\fP(\beta-\gamma)$, is nonempty iff $\gamma=\beta$. In this case
$\CF^\alpha\cap\oZ^\beta\simeq T_{-\beta}\cap S_0$.

{\em Proof.} Evident. $\Box$

\subsubsection{Corollary}
\label{strat F}
Let $\CS$ be a stratum of the fine stratification
of $\CZ^\alpha$. Then either $\CS\cap\CF^\alpha=\emptyset$, or
$\dim(\CS\cap\CF^\alpha)=\frac{1}{2}\dim\CS$.

\subsubsection{} According to the above Lemma, we have
$\CF^\alpha=\sqcup_{\beta\leq\alpha}(\oZ^\beta\cap\CF^\alpha)$ --- a partition
into locally closed subsets. The restriction of $\IC^\alpha$ to
$\oZ^\beta\cap\CF^\alpha$ is concentrated in the degrees $\leq-\dim\oZ^\beta=
-2|\beta|$. Moreover, by the definition of $\IC$ sheaf, if $\beta<\alpha$,
the above inequality is strict.
Since $\dim(\oZ^\beta\cap\CF^\alpha)=|\beta|$ we conclude that
$H^\bullet_c(\oZ^\beta\cap\CF^\alpha,\IC^\alpha|_{\oZ^\beta\cap\CF^\alpha})$
is concentrated in nonpositive degrees, and if $\beta<\alpha$ it is
concentrated in strictly negative degrees.
Applying the Cousin spectral sequence associated with the partition
$\CF^\alpha=\sqcup_{\beta\leq\alpha}(\oZ^\beta\cap\CF^\alpha)$ we see that
$H^\bullet_c(\CF^\alpha,\iota_\alpha^*\IC^\alpha)$ is concentrated in
nonpositive degrees.

\subsubsection{} The above spectral sequence shows also that
$H^0_c(\CF^\alpha,\iota_\alpha^*\IC^\alpha)=H^0_c(\oZ^\alpha\cap\CF^\alpha,
\IC^\alpha|_{\oZ^\alpha\cap\CF^\alpha})=H^0_c(\oZ^\alpha\cap\CF^\alpha,
\ul\BC[2|\alpha|])=H^{2|\alpha|}_c(\oZ^\alpha\cap\CF^\alpha),\BC)=
H^{2|\alpha|}_c(T_{-\alpha}\cap S_0)$. The latter term has a canonical basis
of irreducible components of $T_{-\alpha}\cap S_0$. According to ~\cite{mv},
the number of irreducible components equals $\CK(\alpha)$. This completes
the proof of ~\ref{Weil} b).

\subsection{}
\label{dva}
It remains to show that
$H^\bullet(\BA^\alpha,\fs_\alpha^!\IC^\alpha)$ is concentrated in nonnegative
degrees, or dually, that
$H^\bullet_c(\BA^\alpha,\fs_\alpha^*\IC^\alpha)$ is concentrated in nonpositive
degrees. To this end we study the intersection of the fine stratification of
$\CZ^\alpha$ with $\BA^\alpha=\del_\alpha\CZ^\alpha$.

\subsubsection{Lemma} The intersection of $\del_\alpha\CZ^\alpha$ with a fine
stratum
$\oZ^\gamma\times(\BC^*)^{\beta-\gamma}_\Gamma,\ \gamma\leq\beta\leq\alpha,
\Gamma\in\fP(\beta-\gamma)$, is nonempty iff $\gamma=0$. In this case
$\del_\alpha\CZ^\alpha$ contains the stratum $\oZ^0\times(\BC^*)^\beta_\Gamma
\simeq(\BC^*)^\beta_\Gamma$.

{\em Proof.} Evident. $\Box$

\subsubsection{}
If $\Gamma=\{\{\gamma_1,\ldots,\gamma_k\}\}$ is a partition of $\beta$ then
$\dim(\BC^*)^\beta_\Gamma=k$, and due to the factorization property ~\cite{fm}
~7.3, the stalk of $\IC^\alpha$ at (any point in) the stratum
$(\BC^*)^\beta_\Gamma$ is isomorphic to $\IC^{\alpha-\beta}_\Gamma\simeq
\IC^0_{\{\{\alpha-\beta\}\}}\otimes\bigotimes_{r=1}^k\IC^0_{\{\{\gamma_r\}\}}$
(notations of {\em loc. cit.} ~ 10.7.1, ~10.7.2).

{\em Lemma.} The restriction of $\IC^\alpha$ to the stratum
$(\BC^*)^\beta_\Gamma$ is concentrated in the degrees $\leq-2k$ (here $k$
is the number of elements in the partition $\Gamma$).

{\em Proof.} In view of the above product formula, it suffices to check
that for any $r=1,\ldots,k$ the simple stalk $\IC^0_{\{\{\gamma_r\}\}}$
is concentrated in degree $\leq-2$. This is the stalk of $\IC^{\gamma_r}$
at the one-dimensional (diagonal) stratum
$(\BC^*)^{\gamma_r}_{\{\{\gamma_r\}\}}\subset\CZ^{\gamma_r}$. By the definition
of $\IC$ sheaf, its restriction to any nonopen $l$-dimensional stratum is
concentrated in degrees $<-l$. This completes the proof of the Lemma. $\Box$

\subsubsection{}
We consider the partition of $\del_\alpha\CZ^\alpha$ into locally closed
subsets: $\del_\alpha\CZ^\alpha=\sqcup(\BC^*)^\beta_\Gamma$. According to
the above Lemma, the restriction of $\IC^\alpha$ to a $k$-dimensional stratum
is concentrated in degrees $\leq-2k$.
Applying the Cousin spectral sequence associated with this partition we
conclude, exactly as in subsection ~\ref{leq0}, that
$H^\bullet_c(\BA^\alpha,\fs_\alpha^*\IC^\alpha)$ is concentrated in nonpositive
degrees. Dually,
$H^\bullet(\BA^\alpha,\fs_\alpha^!\IC^\alpha)$ is concentrated in nonnegative
degrees. This completes the proof of the Theorem ~\ref{Weil}. $\Box$

\subsubsection{Corollary} $\CA_\alpha=
\Ext^\bullet_{\CZ^\alpha}(\IC(\del_\alpha\CZ^\alpha),\IC^\alpha)$ is
concentrated in the degree $|\alpha|$.

{\em Proof.} See the beginning of ~\ref{Weil}. $\Box$

\subsection{}
\label{rest}
In the rest of this section we construct the multiplication map
$\CA_\beta\otimes\CA_\alpha\lra\CA_{\beta+\alpha}$. We start with the
following Lemma.

\subsubsection{Lemma}
\label{normal}
$\IC(\del_\alpha\CZ^{\alpha+\beta})=
(\varsigma_{\beta,\alpha})_*\IC(\CZ^\beta\times\BA^\alpha)=
(\varsigma_{\beta,\alpha})_*(\IC^\beta\boxtimes\ul\BC[|\alpha|])$.

{\em Proof.} The map $\varsigma_{\beta,\alpha}$ is finite and generically
one-to-one. $\Box$

\subsubsection{Remark} In fact, $\varsigma_{\beta,\alpha}$ is the normalization
map, but we do not need nor prove this fact.

\subsubsection{}
\label{ef}
We denote by $\IC_\alpha^\beta$ the sheaf
$\IC(\del_\alpha\CZ^{\alpha+\beta})$. In particular,
$\del_0\CZ^\beta=\CZ^\beta$, and $\IC_0^\beta=\IC^\beta$.

Consider the restriction of the map $\varsigma_{\beta+\alpha,\gamma}:\
\CZ^{\beta+\alpha}\times\BA^\gamma\lra\del_\gamma\CZ^{\gamma+\beta+\alpha}$
to $\del_\alpha\CZ^{\beta+\alpha}\times\BA^\gamma$. Evidently,
$\varsigma_{\beta+\alpha,\gamma}(\del_\alpha\CZ^{\beta+\alpha}\times\BA^\gamma)
=\del_{\alpha+\gamma}\CZ^{\gamma+\beta+\alpha}$.
The map $\varsigma_{\beta+\alpha,\gamma}$ restricted to
$\del_\alpha\CZ^{\beta+\alpha}\times\BA^\gamma$ is finite, so the semisimple
perverse sheaf $(\varsigma_{\beta+\alpha,\gamma})_*
\IC(\del_\alpha\CZ^{\beta+\alpha}\times\BA^\gamma)=
(\varsigma_{\beta+\alpha,\gamma})_*(\IC_\alpha^\beta\boxtimes\ul\BC[|\gamma|])$
contains the direct summand $\IC_{\alpha+\gamma}^\beta$ with multiplicity one.
Let $$\IC_{\alpha+\gamma}^\beta\stackrel{e}{\lra}
(\varsigma_{\beta+\alpha,\gamma})_*(\IC_\alpha^\beta\boxtimes\ul\BC[|\gamma|])
\stackrel{f}{\lra}\IC_{\alpha+\gamma}^\beta$$
denote the corresponding inclusion and projection.

\subsection{}
We include the spaces $\CA_\alpha$ into the wider framework.
For $\alpha,\beta,\gamma\in\BN[I]$ we define
$$\CA_\beta^{\alpha,\gamma}:=\Ext^\bullet_{\CZ^{\alpha+\beta+\gamma}}
(\IC^\gamma_{\alpha+\beta},\IC^{\beta+\gamma}_\alpha)$$
In particular, $\CA_\beta=\Ext^\bullet_{\CZ^\beta}(\IC^0_\beta,\IC^\beta_0)=
\CA_\beta^{0,0}$.

Now for $\delta\in\BN[I]$ we construct the {\em stabilization map}
$\tau^{\alpha,\gamma}_{\beta,\delta}:\
\CA_\beta^{\alpha,\gamma}\lra\CA_\beta^{\alpha,\gamma+\delta}$,
and the {\em costabilization map}
$\xi^{\alpha,\gamma}_{\beta,\delta}:\
\CA_\beta^{\alpha,\gamma}\lra\CA_\beta^{\alpha+\delta,\gamma}$.

\subsubsection{}
\label{bezruk}
We choose an open subset $U'\stackrel{j'}{\hookrightarrow}
\CZ^{\alpha+\beta+\gamma+\delta}$ (resp.
$U\stackrel{j}{\hookrightarrow}\CZ^{\alpha+\beta+\gamma}$) such that
$U'\cap\del_{\alpha+\beta}\CZ^{\alpha+\beta+\gamma+\delta}=W':=
\oZ^{\gamma+\delta}\times\BA^{\alpha+\beta}$ (resp.
$U\cap\del_{\alpha+\beta}\CZ^{\alpha+\beta+\gamma}=W:=
\oZ^\gamma\times\BA^{\alpha+\beta}$) (the open subset formed by all the
quasimaps of defect {\em exactly} $\alpha+\beta$).

We have $\Ext^\bullet_{U'}(\IC^{\gamma+\delta}_{\alpha+\beta},
\IC^{\beta+\gamma+\delta}_\alpha)=
\Ext^\bullet_{W'}(\IC^{\gamma+\delta}_{\alpha+\beta},
\varsigma_{\gamma+\delta,\alpha+\beta}^!\IC^{\beta+\gamma+\delta}_\alpha)$, 
and $\Ext^\bullet_U(\IC^\gamma_{\alpha+\beta},
\IC^{\beta+\gamma}_\alpha)=
\Ext^\bullet_W(\IC^\gamma_{\alpha+\beta},
\varsigma_{\gamma,\alpha+\beta}^!\IC^{\beta+\gamma}_\alpha)$, where
$\varsigma_{\gamma,\alpha+\beta}$ stands for the finite map
$\CZ^\gamma\times\BA^{\alpha+\beta}\to
\del_{\alpha+\beta}\CZ^{\alpha+\beta+\gamma}\hookrightarrow
\CZ^{\alpha+\beta+\gamma}$.

Now by the Lemma ~\ref{normal} and the factorization property, there exists
a unique complex of sheaves $\CM$ on $\BA^{\alpha+\beta}$ such that
$(\varsigma_{\gamma+\delta,\alpha+\beta}^!
\IC^{\beta+\gamma+\delta}_\alpha)|_{W'}
=\ul\BC[2|\gamma+\delta|]\boxtimes\CM$, and
$(\varsigma_{\gamma,\alpha+\beta}^!\IC^{\beta+\gamma}_\alpha)|_W=
\ul\BC[2|\gamma|]\boxtimes\CM$. 
As $\IC^\gamma_{\alpha+\beta}|_W=
\ul\BC[2|\gamma|]\boxtimes\ul\BC[|\alpha+\beta|]$, we 
get the restriction map $j^*:\ \Ext^\bullet_{\CZ^{\alpha+\beta+\gamma}}
(\IC^\gamma_{\alpha+\beta},\IC^{\beta+\gamma}_\alpha)\to
\Ext^\bullet_{\oZ^\gamma\times\BA^{\alpha+\beta}}
(\IC(\oZ^\gamma)\boxtimes\IC(\BA^{\alpha+\beta}),\IC(\oZ^\gamma)\boxtimes\CM)$.
We can project the right hand side to the summand 
$\Id\otimes\Ext^\bullet_{\BA^{\alpha+\beta}}(\IC(\BA^{\alpha+\beta}),\CM)$. 
The resulting map from $\Ext^\bullet_{\CZ^{\alpha+\beta+\gamma}}
(\IC^\gamma_{\alpha+\beta},\IC^{\beta+\gamma}_\alpha)$ to
$\Ext^\bullet_{\BA^{\alpha+\beta}}
(\IC(\BA^{\alpha+\beta}),\CM)$ will be also denoted by $j^*$.

Multiplying with $\Id\in\Ext^0_{\CZ^{\gamma+\delta}}(\IC^{\gamma+\delta},
\IC^{\gamma+\delta})$, we get an embedding from 
$\Ext^\bullet_{\BA^{\alpha+\beta}}(\IC(\BA^{\alpha+\beta}),\CM)$ to
$\Ext^\bullet_{\CZ^{\gamma+\delta}\times\BA^{\alpha+\beta}}
(\IC^{\gamma+\delta}\boxtimes\IC(\BA^{\alpha+\beta}),
\IC^{\gamma+\delta}\boxtimes\CM)$.

The key observation is that there exists a unique morphism
$c:\ \IC^{\gamma+\delta}\boxtimes\CM\to 
\varsigma_{\gamma+\delta,\alpha+\beta}^!
\IC^{\beta+\gamma+\delta}_\alpha$ extending the isomorphism
$(\varsigma_{\gamma+\delta,\alpha+\beta}^!
\IC^{\beta+\gamma+\delta}_\alpha)|_{W'}
=\ul\BC[2|\gamma+\delta|]\boxtimes\CM$ from $W'$ to 
$\CZ^{\gamma+\delta}\times\BA^{\alpha+\beta}$. 
In effect, let $\ic^{\beta+\gamma+\delta}_\alpha$ be the mixed Hodge
counterpart of $\IC^{\beta+\gamma+\delta}_\alpha$. It is a pure Hodge module
of weight $w=2|\beta+\gamma+\delta|+|\alpha|$. Then 
$(\varsigma_{\gamma+\delta,\alpha+\beta}^!
\ic^{\beta+\gamma+\delta}_\alpha)|_{W'}
=\ul\BQ[2|\gamma+\delta|]\boxtimes\bM$ for some Hodge module $\bM$ on 
$\BA^{\alpha+\beta}$. It follows from the Proposition ~\ref{stunt} below that
$\bM$ is a pure (hence semisimple) 
Hodge complex of weight $|\alpha+2\beta|$. Thus
$\ic^{\gamma+\delta}\boxtimes\bM$ is a subquotient of the direct sum of
cohomology of the mixed Hodge complex 
$\varsigma_{\gamma+\delta,\alpha+\beta}^!\ic^{\beta+\gamma+\delta}_\alpha$.
The latter complex has weight $\geq w$, while 
$\ic^{\gamma+\delta}\boxtimes\bM$ is pure of weight $w$. The desired morphism
$c$ is constructed by induction in the cohomology degree using the Proposition
5.1.15 of ~\cite{bbd}.

Finally, we define
$$\tau^{\alpha,\gamma}_{\beta,\delta}:\
\CA_\beta^{\alpha,\gamma}=\Ext^\bullet_{\CZ^{\alpha+\beta+\gamma}}
(\IC^\gamma_{\alpha+\beta},\IC^{\beta+\gamma}_\alpha)\lra
\Ext^\bullet_{\CZ^{\alpha+\beta+\gamma+\delta}}
(\IC^{\gamma+\delta}_{\alpha+\beta},\IC^{\beta+\gamma+\delta}_\alpha)=
\CA_\beta^{\alpha,\gamma+\delta}$$
as the following composition:
$\Ext^\bullet_{\CZ^{\alpha+\beta+\gamma}}
(\IC^\gamma_{\alpha+\beta},\IC^{\beta+\gamma}_\alpha)
\stackrel{j^*}{\lra}
\Ext^\bullet_{\BA^{\alpha+\beta}}
(\IC(\BA^{\alpha+\beta}),\CM)
\stackrel{\Id\otimes ?}{\lra}
\Ext^\bullet_{\CZ^{\gamma+\delta}\times\BA^{\alpha+\beta}}
(\IC^{\gamma+\delta}\boxtimes\IC(\BA^{\alpha+\beta}),
\IC^{\gamma+\delta}\boxtimes\CM)
\stackrel{c}{\lra}
\Ext^\bullet_{\CZ^{\gamma+\delta}\times\BA^{\alpha+\beta}}
(\IC^{\gamma+\delta}\boxtimes\IC(\BA^{\alpha+\beta}),
\varsigma_{\gamma+\delta,\alpha+\beta}^!
\IC^{\beta+\gamma+\delta}_\alpha)=
\Ext^\bullet_{\CZ^{\alpha+\beta+\gamma+\delta}}
(\IC^{\gamma+\delta}_{\alpha+\beta},\IC^{\beta+\gamma+\delta}_\alpha)$.
Quite obviously, the result does not depend on a choice of $U$ and $U'$.

\subsubsection{}
To define the costabilization map
$\xi^{\alpha,\gamma}_{\beta,\delta}:\
\CA_\beta^{\alpha,\gamma}\lra\CA_\beta^{\alpha+\delta,\gamma}$
we note that
$\CA_\beta^{\alpha,\gamma}:=\Ext^\bullet_{\CZ^{\alpha+\beta+\gamma}}
(\IC^\gamma_{\alpha+\beta},\IC^{\beta+\gamma}_\alpha)=
\Ext^\bullet_{\CZ^{\alpha+\beta+\gamma}\times\BA^\delta}
(\IC^\gamma_{\alpha+\beta}\boxtimes\ul\BC[|\delta|],
\IC^{\beta+\gamma}_\alpha\boxtimes\ul\BC[|\delta|])$ maps by
$(\varsigma_{\alpha+\beta+\gamma,\delta})_*$ to
$\Ext^\bullet_{\CZ^{\alpha+\beta+\gamma+\delta}}
((\varsigma_{\alpha+\beta+\gamma,\delta})_*
(\IC^\gamma_{\alpha+\beta}\boxtimes\ul\BC[|\delta|]),
(\varsigma_{\alpha+\beta+\gamma,\delta})_*
(\IC^{\beta+\gamma}_\alpha\boxtimes\ul\BC[|\delta|]))$.

According to ~\ref{ef}, the first argument of the latter Ext contains
the canonical direct summand $e\IC^\gamma_{\alpha+\delta+\beta}$, while
the second argument of the latter Ext canonically projects by $f$ to
$\IC^{\beta+\gamma}_{\alpha+\delta}$. Thus, the latter Ext canonically
projects to
$\Ext^\bullet_{\CZ^{\alpha+\beta+\gamma+\delta}}
(\IC^\gamma_{\alpha+\delta+\beta},
\IC^{\beta+\gamma}_{\alpha+\delta})=\CA_\beta^{\alpha+\delta,\gamma}$.
Composing this projection with the above map
$(\varsigma_{\alpha+\beta+\gamma,\delta})_*$ we obtain the desired map
$\xi^{\alpha,\gamma}_{\beta,\delta}:\
\CA_\beta^{\alpha,\gamma}\lra\CA_\beta^{\alpha+\delta,\gamma}$.

\subsection{Proposition}
\label{abs non}
For $\alpha,\beta,\gamma,\delta,\epsilon\in\BN[I]$ we have

a) $\tau_{\beta,\epsilon}^{\alpha,\gamma+\delta}\circ
\tau_{\beta,\delta}^{\alpha,\gamma}=
\tau_{\beta,\delta+\epsilon}^{\alpha,\gamma}:\
\CA_\beta^{\alpha,\gamma}\lra\CA_\beta^{\alpha,\gamma+\delta+\epsilon}$;

b) $\xi^{\alpha+\delta,\gamma}_{\beta,\epsilon}\circ
\xi^{\alpha,\gamma}_{\beta,\delta}=
\xi^{\alpha,\gamma}_{\beta,\delta+\epsilon}:\
\CA_\beta^{\alpha,\gamma}\lra\CA_\beta^{\alpha+\delta+\epsilon,\gamma}$;

c) $\xi^{\alpha,\gamma+\delta}_{\beta,\epsilon}\circ
\tau^{\alpha,\gamma}_{\beta,\delta}=
\tau^{\alpha+\epsilon,\gamma}_{\beta,\delta}\circ
\xi^{\alpha,\gamma}_{\beta,\epsilon}:\
\CA_\beta^{\alpha,\gamma}\lra\CA_\beta^{\alpha+\epsilon,\gamma+\delta}$.

{\em Proof.} Routine check. $\Box$

\subsection{Definition}
\label{multip}
Let $a\in\CA_\alpha=\CA_\alpha^{0,0},\
b\in\CA_\beta=\CA_\beta^{0,0}$. We define the product $a\cdot b$ as
the following composition:
$$a\cdot b:=\tau_{\alpha,\beta}^{0,0}(a)\circ\xi_{\beta,\alpha}^{0,0}(b)
\in\CA_{\alpha+\beta}^{0,0}=\CA_{\alpha+\beta}$$

The associativity of this multiplication follows immediately from the
Proposition ~\ref{abs non}.

\subsection{}
\label{stab}
To reward the patient reader who has fought his way through the above notation,
we repeat here the definition of $\tau_{\alpha,\beta}^{\gamma,\delta}$
in the particular case of $\gamma=\delta=0$. This is the only case needed
in the definition of multiplication (the general case is needed to prove the
associativity), and the definition simplifies somewhat in this case.

So we are going to define the map
$$\tau_{\alpha,\beta}^{0,0}:\ \CA_\alpha
=\Ext^\bullet_{\CZ^\alpha}(\IC(\del_\alpha\CZ^\alpha),\IC^\alpha)\lra
\Ext^\bullet_{\CZ^{\alpha+\beta}}(\IC(\del_\alpha\CZ^{\alpha+\beta}),
\IC^{\alpha+\beta})$$

First of all, tensoring with $\Id\in\Ext^0_{\CZ^\beta}(\IC^\beta,\IC^\beta)$,
we get the map from 
$\Ext^\bullet_{\CZ^\alpha}(\IC(\del_\alpha\CZ^\alpha),\IC^\alpha)=
\Ext^\bullet_{\BA^\alpha}(\IC(\BA^\alpha),\fs_\alpha^!\IC^\alpha)$
to $\Ext^\bullet_{\CZ^\beta\times\BA^\alpha}(\IC^\beta\boxtimes\IC(\BA^\alpha),
\IC^\beta\boxtimes\fs_\alpha^!\IC^\alpha)$. 

Now the same argument as in
~\ref{bezruk} shows that there is a unique morphism $c$ from 
$\IC^\beta\boxtimes\fs_\alpha^!\IC^\alpha$ to 
$\varsigma_{\beta,\alpha}^!\IC^{\alpha+\beta}$ extending the factorization
isomorphism from the open part $\oZ^\beta\times\BA^\alpha$ of 
$\CZ^\beta\times\BA^\alpha$. Thus $c$ induces the map from 
$\Ext^\bullet_{\CZ^\beta\times\BA^\alpha}(\IC^\beta\boxtimes\IC(\BA^\alpha),
\IC^\beta\boxtimes\fs_\alpha^!\IC^\alpha)$ to
$\Ext^\bullet_{\CZ^\beta\times\BA^\alpha}(\IC^\beta\boxtimes\IC(\BA^\alpha),
\varsigma_{\beta,\alpha}^!\IC^{\alpha+\beta})$.

The latter space equals $\Ext^\bullet_{\CZ^{\alpha+\beta}}
((\varsigma_{\beta,\alpha})_*(\IC^\beta\boxtimes\IC(\BA^\alpha)),
\IC^{\alpha+\beta})=
\Ext^\bullet_{\CZ^{\alpha+\beta}}(\IC(\del_\alpha\CZ^{\alpha+\beta}),
\IC^{\alpha+\beta})$.

Finally, we define $\tau_{\alpha,\beta}^{0,0}$ as the composition of above
maps.

\subsection{}
We close this section with a definition of comultiplication
$\Delta:\ \CA_{\alpha}\lra\oplus_{\beta+\gamma=\alpha}
\CA_\beta\otimes\CA_\gamma$.

\subsubsection{} We choose two disjoint open balls $\Omega,\Upsilon\subset\BC=
\BA^1$. Then $\Omega^\beta\times\Upsilon^\gamma\subset\BA^{\beta+\gamma}=
\BA^\alpha$. According to the factorization property of $\CZ^\alpha$,
we have $(\fs_\alpha^!\IC^\alpha)|_{\Omega^\beta\times\Upsilon^\gamma}=
\fs_\beta^!\IC^\beta|_{\Omega^\beta}\boxtimes
\fs_\gamma^!\IC^\gamma|_{\Upsilon^\gamma}$. We define the
comultiplication component $\Delta_{\beta,\gamma}:\ \CA_{\alpha}\lra
\CA_\beta\otimes\CA_\gamma$ as the following composition:
$\CA_\alpha=H^0(\BA^\alpha,\fs_\alpha^!\IC^\alpha)\lra
H^0(\Omega^\beta,\fs_\beta^!\IC^\beta)\otimes
H^0(\Upsilon^\gamma,\fs_\gamma^!\IC^\gamma)=
H^0(\BA^\beta,\fs_\beta^!\IC^\beta)\otimes
H^0(\BA^\gamma,\fs_\gamma^!\IC^\gamma)=\CA_\beta\otimes\CA_\gamma$.

Quite evidently, the result does not depend on a choice of $\Omega$ and
$\Upsilon$. This comultiplication is manifestly cocommutative and
coassociative. It also commutes with the multiplication defined in
~\ref{multip}, so it defines the structure of cocommutative bialgebra on $\CA$.

\subsubsection{}
\label{crit}
Consider the main diagonal stratum
$\BA^1\stackrel{u_\alpha}{\hookrightarrow}\BA^\alpha$. We will
denote by $U^\alpha\stackrel{j}{\hookrightarrow}\BA^\alpha$ the open inclusion
of the complement $U^\alpha=\BA^\alpha-\BA^1$. We say that a cohomology class
$a\in H^0(\BA^\alpha,\fs_\alpha^!\IC^\alpha)$ {\em is supported on the main
diagonal} if $0=j^*a\in H^0(U^\alpha,\fs_\alpha^!\IC^\alpha)$.
The following criterion will prove very useful in section ~4.

{\em Lemma.} $a\in\CA_\alpha$ is supported on the main diagonal iff
$a$ is primitive, i.e. $\Delta(a)=1\otimes a+a\otimes1$.

{\em Proof.} a) If $a\in\CA_\alpha$ is supported on the main diagonal
then, obviously, $\Delta_{\beta,\gamma}(a)=0$ when $\beta>0<\gamma$.

b) The converse follows by the Mayer-Vietoris type argument since $U^\alpha$
may be covered by the open subsets of the type
$\Omega^\beta\times\Upsilon^\gamma,\ \beta>0<\gamma. \quad \Box$

\section{Simple $\IC$-stalks}

\subsection{}
\label{simple}
Let $\check\CR{}^+\subset Y$ denote the set of positive coroots.

For $\alpha\in\BN[I]$ the following $q$-analogue of the Kostant partition
function $\CK(\alpha)$ was introduced in ~\cite{lus}:
$$\CK^\alpha(t):=\sum_\kappa t^{-|\kappa|}$$
Here the sum is taken over the set of partitions of $\alpha$ into a sum of
positive coroots $\ctheta\in\check\CR{}^+$, and for a partition $\kappa$ the
number of elements in $\kappa$ is denoted by $|\kappa|$. Finally, $t$ is a
formal variable of degree 2. We have $\CK(\alpha)=\CK^\alpha(1)$.

In this section we prove the following Theorem:

{\bf Theorem.} $\CK^\alpha(t)$ is the generating function of the simple stalk
$\IC^0_{\{\{\alpha\}\}}$ (see ~\cite{fm} ~10.7.1), that is (notations of
~\ref{raz}),

a) for odd $k$ we have $H^k\iota_\alpha^*\fs_\alpha^*\IC^\alpha=0$;

b) $\dim H^{2k}\iota_\alpha^*\fs_\alpha^*\IC^\alpha$ equals the coefficient
$\CK^\alpha_k$ of $t^k$ in $\CK^\alpha(t)$.

\subsubsection{Corollary}
\label{berezin}
The generating function of the stalk $\IC^0_\Gamma,\
\Gamma=\{\{\gamma_1,\ldots,\gamma_k\}\}$ (see ~\cite{fm} ~10.7.1) equals
$\prod_{r=1}^k\CK^{\gamma_r}(t)$.

The Corollary follows immediately from the above Theorem and
the factorization property. In its turn, the Corollary (along with the
Proposition ~10.7.3 of ~\cite{fm}) implies the Parity vanishing conjecture
~10.7.4 of {\em loc. cit.}

\subsection{}
We start with the proof of ~\ref{simple} ~a).
%The proof is a simple corollary
%of Localization Theorem for equivariant cohomology (see e.g. ~\cite{em} ~B.2).
Recall that the Cartan group $\bH$ acts on $\CZ^\alpha$ contracting
it to the fixed point set $(\CZ^\alpha)^\bH=\del_\alpha\CZ^\alpha=\BA^\alpha$.
Moreover, if $\BC^*$ is the group of automorphisms of $C=\BP^1$ preserving
$0,\infty$, then $\bH\times\BC^*$ acts on $\CZ^\alpha$ contracting it to the
fixed point set $(\CZ^\alpha)^{\bH\times\BC^*}=\fs_\alpha\circ\iota_\alpha(0)$.
Applying the Corollary ~14.3 of ~\cite{bl} we see that the equivariant
cohomology $H^\bullet_{\bH\times\BC^*}(\CZ^\alpha,\IC^\alpha)$ is
isomorphic to $H^\bullet_{\bH\times\BC^*}(\cdot)\otimes
H^\bullet(\CZ^\alpha,\IC^\alpha)$, and
$H^\bullet_\bH(\CZ^\alpha,\IC^\alpha)\cong
H^\bullet_\bH(\cdot)\otimes H^\bullet(\CZ^\alpha,\IC^\alpha)=
H^\bullet_\bH(\cdot)\otimes \iota_\alpha^*\fs_\alpha^*\IC^\alpha$.

On the other hand, according to the Theorem ~B.2 of ~\cite{em}, the natural map
$$H^\bullet_\bH(\cdot)\otimes H^\bullet(\BA^\alpha,\fs_\alpha^!\IC^\alpha)=
H^\bullet_\bH(\BA^\alpha,\fs_\alpha^!\IC^\alpha)\lra
H^\bullet_\bH(\CZ^\alpha,\IC^\alpha)$$
is an isomorphism after certain localization in $H^\bullet_\bH(\cdot)$.

Now $H^\bullet_\bH(\cdot)$ is concentrated in even degrees, and the above
localization is taken with respect to a homogeneous multiplicative subset.
By the Theorem ~\ref{Weil},
$H^\bullet(\BA^\alpha,\fs_\alpha^!\IC^\alpha)$ is concentrated in degree 0,
hence the (localized)
$H^\bullet_\bH(\cdot)\otimes H^\bullet(\BA^\alpha,\fs_\alpha^!\IC^\alpha)=
H^\bullet_\bH(\BA^\alpha,\fs_\alpha^!\IC^\alpha)$ is concentrated in even
degrees, hence the (localized)
$H^\bullet_\bH(\CZ^\alpha,\IC^\alpha)\cong
H^\bullet_\bH(\cdot)\otimes H^\bullet(\CZ^\alpha,\IC^\alpha)=
H^\bullet_\bH(\cdot)\otimes \iota_\alpha^*\fs_\alpha^*\IC^\alpha$
is concentrated in even degrees, hence
$\iota_\alpha^*\fs_\alpha^*\IC^\alpha$ is concentrated in even degrees.
This completes the proof of the Theorem ~\ref{simple} ~a).

\subsubsection{Corollary of the proof}
\label{cor}
$\dim(\iota_\alpha^*\fs_\alpha^*\IC^\alpha)=
\chi(\iota_\alpha^*\fs_\alpha^*\IC^\alpha)=\CK(\alpha)$.

\subsection{} To prove ~\ref{simple} ~b) we need to introduce and recall
some notation.

\subsubsection{}
As we have seen in ~\ref{strat F}, the fine stratification
of $\CZ^\alpha$ induces a partition of the central fiber $\CF^\alpha$ into
the locally closed subsets: $\CF^\alpha=\sqcup_{\beta\leq\alpha}
(\oZ^\beta\cap\CF^\alpha)=:\sqcup_{\beta\leq\alpha}\oF^\beta$. Moreover,
$\oF^\beta\simeq T_{-\beta}\cap S_0$. The closed embedding of the smallest
point stratum $\oF^0$ into $\CF^\alpha$ will be denoted by $\fs_\alpha$.
Recall that the closed embedding $\CF^\alpha\hookrightarrow\CZ^\alpha$ is
denoted by $\iota_\alpha$. Certainly, we have
$\iota_\alpha^*\fs_\alpha^*\IC^\alpha=\fs_\alpha^*\iota_\alpha^*\IC^\alpha$.

\subsubsection{}
\label{Hodge}
We will use the machinery of {\em mixed Hodge modules}. The sheaf $\IC^\alpha$
has its mixed Hodge counterpart $\ic^\alpha$ --- the irreducible pure Hodge
module on $\CZ^\alpha$ of weight $2|\alpha|$.
The following Theorem is a strong version of ~\ref{simple}.

{\bf Theorem.}
a) for odd $k$ we have $H^k\fs_\alpha^*\iota_\alpha^*\ic^\alpha=0$;

b) $H^{2k}\fs_\alpha^*\iota_\alpha^*\ic^\alpha$ is the direct sum of
$\CK^\alpha_k$ copies of the Tate Hodge structure $\BQ(k+|\alpha|)$.

\subsubsection{}
\label{Grass}
Recall that for a dominant $\eta\in Y^+$ we have a $\bG(\CO)$-orbit
$\CG_\eta$ and its closure $\oCG_\eta$ in the affine Grassmannian $\CG$.
Recall that $w_0$ is the longest element in the Weyl group $\CW_f$ of $\bG$.
For a sufficiently dominant $\eta\in Y^+$, the difference $\eta-w_0(\alpha)=:
\vartheta$ is also dominant, and $\CG_\eta\subset\oCG_\vartheta$.
Let $\ic(\oCG_\vartheta)_\eta$ denote the stalk of the Hodge module
$\ic(\oCG_\vartheta)$ at any point in $\CG_\eta$. G.Lusztig has proved the
following theorem in ~\cite{lus} (in fact, he used the language of Weil sheaves
instead of Hodge modules):

{\bf Theorem.} Let $\eta\in Y^+$ be sufficiently dominant. Then

a) for odd $k$ we have $H^k(\ic(\oCG_\vartheta)_\eta[-2|\eta|])=0$;

b) $H^{2k}(\ic(\oCG_\vartheta)_\eta[-2|\eta|])$ is the direct sum of
$\CK^\alpha_k$ copies of the Tate Hodge structure $\BQ(k+|\alpha|)$.

\subsubsection{}
\label{hard}
Comparing ~\ref{Grass} and ~\ref{Hodge} we see that the Theorems ~\ref{Hodge}
and ~\ref{simple} follow from the following

{\bf Theorem.}
$\ic(\oCG_\vartheta)_\eta[-2|\eta|]\cong\fs_\alpha^*\iota_\alpha^*\ic^\alpha$.

The rest of this section is devoted to the proof of this Theorem.

\subsection{}
We consider the intersection of semiinfinite orbits $\ol{T}_{w_0\eta-\alpha}
\cap S_{w_0\eta}$ in the affine Grassmannian $\CG$. It lies in the closure
$\oCG_\vartheta$, and the translation by the Cartan loop $w_0\eta(z)\in\bG((z))$
identifies it with $\ol{T}_{-\alpha}\cap S_0\simeq\CF^\alpha$. Moreover, the
partition of $\ol{T}_{w_0\eta-\alpha}\cap S_{w_0\eta}$ into intersections
with $\bG(\CO)$-orbits $\CG_{\eta-w_0\beta},\ 0\leq\beta\leq\alpha$,
corresponds under this identification to the partition
$\CF^\alpha=\sqcup_{0\leq\beta\leq\alpha}\oF^\beta$.
We denote the locally closed embedding $\CF^\alpha\simeq
\ol{T}_{w_0\eta-\alpha}\cap S_{w_0\eta}\hookrightarrow\oCG_\vartheta$ by
$\bi_\alpha$. Thus,
$\ic(\oCG_\vartheta)_\eta=\fs_\alpha^*\bi_\alpha^*\ic(\oCG_\vartheta)$.

\subsubsection{Conjecture} $\bi_\alpha^*\ic(\oCG_\vartheta)[-2|\eta|]\cong
\iota_\alpha^*\ic^\alpha$.

This conjecture would imply the Theorem ~\ref{hard}. Unfortunately, we cannot
prove this conjecture at the moment. Instead of it, we will prove its version
in the $K$-group of mixed Hodge modules.

\subsection{Proposition}
\label{stunt}
a) $\fs_\alpha^*\iota_\alpha^*\ic^\alpha$ is a pure Hodge complex of weight
$2|\alpha|$;

b) $\fs_\alpha^!\iota_\alpha^*\ic^\alpha$ is a pure Hodge structure of weight
$2|\alpha|$.

{\em Proof.} a) The natural map $H^\bullet(\CZ^\alpha,\ic^\alpha)\lra
\fs_\alpha^*\iota_\alpha^*\ic^\alpha$ is an isomorphism since $\ic^\alpha$
is smooth along the fine stratification. Now $\ic^\alpha$ is pure of weight
$2|\alpha|$, and $H^\bullet(?)$ increases weights, while both
$\iota_\alpha^*(?)$ and $\fs_\alpha^*(?)$ decrease weights.

b) We already know that
$\fs_\alpha^!\iota_\alpha^*\ic^\alpha=
H^\bullet(\BA^\alpha,\fs_\alpha^!\ic^\alpha)$, and both $H^\bullet(?)$ and
$\fs_\alpha^!(?)$ increase weights. On the other hand,
$\fs_\alpha^!\iota_\alpha^*\ic^\alpha=
\pi_{\alpha!}\iota_\alpha^*\ic^\alpha$, and both $\pi_{\alpha!}(?)$ and
$\iota_\alpha^*(?)$ decrease weights. $\Box$

\subsection{}
The above Proposition implies that to prove the Theorems ~\ref{hard} and
~\ref{Hodge} it suffices to identify the classes of
$[\fs_\alpha^*\iota_\alpha^*\ic^\alpha]$ and
$[\ic(\oCG_\vartheta)_\eta]$ in the $K$-group of mixed Hodge structures
$K(MHM(\cdot))$. We will proceed by induction in $\alpha$. The first step
$\alpha=0$ being trivial, we may suppose that our Theorems are established
for all $\beta<\alpha$.

Let us denote the complement $\CF^\alpha-\oF^0$ by $U$, and let us denote
its open embedding into $\CF^\alpha$ by $j$. By factorization,
$\iota_\alpha^*\ic^\alpha|_{\oF^\beta}$ is the constant Hodge module with the
stalk $\fs_{\alpha-\beta}^*\iota_{\alpha-\beta}^*\ic^{\alpha-\beta}[2|\beta|]$.
Thus, by induction, we know that the classes of
$[j^*\iota_\alpha^*\ic^\alpha]$ and
$[j^*\bi_\alpha^*\ic(\oCG_\vartheta)]$ in the $K$-group of mixed Hodge modules
$K(MHM(U))$ coincide.

Let us consider the following exact triangles of mixed Hodge complexes:
\begin{equation}
\label{1}
\fs_\alpha^!\iota_\alpha^*\ic^\alpha\lra
\fs_\alpha^*\iota_\alpha^*\ic^\alpha\lra
\fs_\alpha^*j_*j^*\iota_\alpha^*\ic^\alpha\lra
\fs_\alpha^!\iota_\alpha^*\ic^\alpha[1]\lra\ldots
\end{equation}

\begin{equation}
\label{2}
\fs_\alpha^!\bi_\alpha^*\ic(\oCG_\vartheta)\lra
\fs_\alpha^*\bi_\alpha^*\ic(\oCG_\vartheta)\lra
\fs_\alpha^*j_*j^*\bi_\alpha^*\ic(\oCG_\vartheta)\lra
\fs_\alpha^!\bi_\alpha^*\ic(\oCG_\vartheta)[1]\lra\ldots
\end{equation}

\subsubsection{}
We already know that
$\fs_\alpha^*\iota_\alpha^*\ic^\alpha$ lives in the even negative degrees, and
$\fs_\alpha^!\iota_\alpha^*\ic^\alpha$ lives in degree 0. Hence the differential
$\fs_\alpha^!\iota_\alpha^*\ic^\alpha\lra
\fs_\alpha^*\iota_\alpha^*\ic^\alpha$ vanishes (since the category of mixed
Hodge structures has the cohomological dimension 1), and
$\fs_\alpha^*j_*j^*\iota_\alpha^*\ic^\alpha=
\fs_\alpha^!\iota_\alpha^*\ic^\alpha[1]\oplus
\fs_\alpha^*\iota_\alpha^*\ic^\alpha$.

\subsubsection{}
By ~\cite{mv} we know that
$\fs_\alpha^!\bi_\alpha^*\ic(\oCG_\vartheta)$ lives in degree 0 (and has
dimension $\CK(\alpha)$). Now the same argument as above proves that
$\fs_\alpha^*j_*j^*\bi_\alpha^*\ic(\oCG_\vartheta)=
\fs_\alpha^!\bi_\alpha^*\ic(\oCG_\vartheta)[1]\oplus
\fs_\alpha^*\bi_\alpha^*\ic(\oCG_\vartheta)$.

Moreover, by ~\cite{mv} we know that
$\fs_\alpha^!\bi_\alpha^*\ic(\oCG_\vartheta)$ is a direct summand of
$H^\bullet(\oCG_\vartheta,\ic(\oCG_\vartheta))$, and hence it is a pure Tate Hodge
module of weight $2|\alpha|$. We deduce that
$\fs_\alpha^*j_*j^*\bi_\alpha^*\ic(\oCG_\vartheta)$ is a Tate Hodge module.

\subsubsection{}
By the induction assumption, the classes of
$[\fs_\alpha^*j_*j^*\bi_\alpha^*\ic(\oCG_\vartheta)]$ and
$[\fs_\alpha^*j_*j^*\iota_\alpha^*\ic^\alpha]$ in $K(MHM(\cdot))$ coincide.
So they both lie in the subgroup $K(TMHM(\cdot))$ of Tate Hodge modules.
We have
$[\fs_\alpha^*j_*j^*\iota_\alpha^*\ic^\alpha]=
[\fs_\alpha^*\iota_\alpha^*\ic^\alpha]-
[\fs_\alpha^!\iota_\alpha^*\ic^\alpha]$, and no cancellations
occur in the RHS. In effect, all graded parts of
$\fs_\alpha^*\iota_\alpha^*\ic^\alpha$ have weights $<2|\alpha|$ (since they
live in negative degrees), while
$\fs_\alpha^!\iota_\alpha^*\ic^\alpha$ is pure of weight $2|\alpha|$ (see
~\ref{stunt}). We deduce that both
$[\fs_\alpha^*\iota_\alpha^*\ic^\alpha]$ and
$[\fs_\alpha^!\iota_\alpha^*\ic^\alpha]$ lie in $K(TMHM(\cdot))$. Thus
$\fs_\alpha^!\iota_\alpha^*\ic^\alpha$ is a sum of $\CK(\alpha)$ copies of
$\BQ(|\alpha|)$, hence
$\fs_\alpha^!\iota_\alpha^*\ic^\alpha\cong
\fs_\alpha^!\bi_\alpha^*\ic(\oCG_\vartheta)$.
It follows that
$[\fs_\alpha^*\iota_\alpha^*\ic^\alpha]=
[\fs_\alpha^*\bi_\alpha^*\ic(\oCG_\vartheta)]$.
This completes the proof of the Theorem ~\ref{hard} along with the Theorems
~\ref{Hodge} and ~\ref{simple}. $\Box$

\subsection{Remark} It is easy to describe the semisimple complex
$\fs_\alpha^!\IC^\alpha$ on $\BA^\alpha$. For a Kostant partition 
$\kappa\in\fK(\alpha)$ let $\BA^\alpha_\kappa\subset\BA^\alpha$ denote the
closure of the corresponding diagonal stratum. It is known that the
normalization 
$\tilde\BA{}^\alpha_\kappa\stackrel{N^\kappa}{\lra}\BA^\alpha_\kappa$ is
isomorphic to an affine space. We have $\fs_\alpha^!\IC^\alpha\simeq
\bigoplus_{\kappa\in\fK(\alpha)}N^\kappa_*\ul\BC$.

\section{$\fg^L$-action on the global $\IC$-cohomology}

\subsection{}
\label{ostrik}
We compute the global Intersection Cohomology of $\CQ^\alpha$ following ~\S2 of
~\cite{fk}. First we introduce the open subset $\ddQ^\alpha\subset\CQ^\alpha$
formed by all the quasimaps $\phi$ such that the defect of $\phi$ does not
meet $\infty\in C$. Twisting by the multiples of $\infty$ we obtain the
partition into locally closed subsets $\CQ^\alpha=\sqcup_{\beta\leq\alpha}
\ddQ^\beta$. The evaluation at $\infty$ defines the projection
$p_\alpha:\ \ddQ^\alpha\lra\bX$ which is a locally trivial fibration with
a fiber isomorphic to $\CZ^\alpha$.

By the Poincar\'e duality and Theorem ~\ref{simple} the generating function
of $H^\bullet_c(\CZ^\beta,\IC^\beta)$ is given by $\CK^\beta(t^{-1})$.
By the parity vanishing, the spectral sequence of fibration
$p_\beta:\ \ddQ^\beta\lra\bX$ degenerates and provides us with the
isomorphism of graded vector spaces $H^\bullet_c(\ddQ^\beta,\IC(\CQ^\beta))
\cong H^\bullet(\bX,\ul\BC[\dim\bX])\otimes
H^\bullet_c(\CZ^\beta,\IC^\beta)$.
By the factorization property,
$\IC(\CQ^\alpha)|_{\ddQ^\beta}\simeq
\IC(\CQ^\beta)\otimes\IC^0_{\{\{\alpha-\beta\}\}}$. Finally, the parity
vanishing implies the degeneration of the Cousin spectral sequence associated
with the partition $\CQ^\alpha=\sqcup_{\beta\leq\alpha}\ddQ^\beta$.

Combining all the above equalities we arrive at the formula for the
generating function of $H^\bullet(\CQ^\alpha,\IC(\CQ^\alpha))$.
To write it down in a neat form we will consider the generating function
$P_\bG(t)$ of $\oplus_{\alpha\in\BN[I]}H^\bullet(\CQ^\alpha,\IC(\CQ^\alpha))$.
To record the information on $\alpha$ we will consider this generating
function as a formal cocharacter of $\bH$ with coefficients in the Laurent
polynomials in $t$. Formal cocharacters will be written multiplicatively,
so that the cocharacter corresponding to $\alpha$ will be denoted by $e^\alpha$.
Finally, for the reasons which will become clear later
(see Proposition ~\ref{h}),
we will make the following rescaling. We will attach to
$H^\bullet(\CQ^\alpha,\IC(\CQ^\alpha))$ the cocharacter
$e^{\alpha+2\crho}$.

With all this in mind, the Poincar\'e polynomial $P_\bG(t)$ of
$\oplus_{\alpha\in\BN[I]}H^\bullet(\CQ^\alpha,\IC(\CQ^\alpha))$
is calculated as follows:

{\bf Theorem.} $$P_\bG(t)=\frac{e^{2\check\rho}t^{-\frac{1}{2}\dim\bX}
\sum_{w\in\CW_f}t^{l(w)}}
{\prod_{\ctheta\in{\check\CR}{}^+}(1-te^\ctheta)(1-t^{-1}e^\ctheta)}$$
where $l(w)$ stands for the length of $w\in\CW_f$, and $2\check\rho$
stands for the sum of positive coroots $\ctheta\in{\check\CR}{}^+$.
$\Box$

\subsubsection{}
\label{feigin}
The above argument along with the Theorem ~\ref{Hodge} establishes also the
following Theorem:

{\bf Theorem.} a) For odd $k$ we have
$H^{k-2|\alpha|-\dim\bX}(\CQ^\alpha,\ic(\CQ^\alpha))=0$;

b) $H^{2k-2|\alpha|-\dim\bX}(\CQ^\alpha,\ic(\CQ^\alpha))$ is the direct sum
of a few copies of the Tate Hodge module $\BQ(k). \quad \Box$

\subsection{}
\label{prim}
In the rest of this section we equip
$\oplus_{\alpha\in\BN[I]}H^\bullet(\CQ^\alpha,\IC(\CQ^\alpha))$
with the structure of the Langlands dual Lie algebra $\fg^L$-module with the
character $\frac{|\CW_f|e^{2\check\rho}}
{\prod_{\ctheta\in{\check\CR}{}^+}(1-e^\ctheta)^2}$.

We start with the following Proposition describing the subspace Prim$(\CA)$
of the primitive elements of bialgebra $\CA$.

{\bf Proposition.} a) If $\alpha\in\BN[I]$ is not a positive coroot, that is
$\alpha\not\in\cR^+$, then $\CA_\alpha$ does not contain primitive elements;

b) dim(Prim$(\CA)\cap\CA_\alpha)\leq1$ for any $\alpha\in\BN[I]$.

{\em Proof.} We will use the Criterion ~\ref{crit}. Consider the long exact
sequence of cohomology (notations of ~\ref{crit}):
$$\ldots\lra H^\bullet(\BA^1,u_\alpha^!\fs_\alpha^!\IC^\alpha)\lra
H^\bullet(\BA^\alpha,\fs_\alpha^!\IC^\alpha)\lra
H^\bullet(U^\alpha,\fs_\alpha^!\IC^\alpha)\lra\ldots$$

Applying Poincar\'e duality to the Theorem ~\ref{simple} we see that
$u_\alpha^!\fs_\alpha^!\IC^\alpha$ is the direct sum of constant sheaves
living in nonnegative degrees, and $H^0(u_\alpha^!\fs_\alpha^!\IC^\alpha)$
is nonzero iff $\alpha\in\cR^+$. In this case
$H^0(u_\alpha^!\fs_\alpha^!\IC^\alpha)=\ul\BC$. We deduce that
$H^0(\BA^1,u_\alpha^!\fs_\alpha^!\IC^\alpha)$ is nonzero iff $\alpha\in\cR^+$,
and in this case
$\dim H^0(\BA^1,u_\alpha^!\fs_\alpha^!\IC^\alpha)=1$.
Now an element $a\in H^0(\BA^\alpha,\fs_\alpha^!\IC^\alpha)$ is supported on
the main diagonal iff it comes from
$H^0(\BA^1,u_\alpha^!\fs_\alpha^!\IC^\alpha)$. This completes the proof of
the Proposition. $\Box$

\subsubsection{Remark} Later on we will see that in fact the converse of
~\ref{prim} ~a) is also true (see ~\ref{!}).

\subsection{Definition}
\label{a}
a) $\ba$ is the Lie algebra formed by the primitive
elements Prim$(\CA)$;

b) $\fa\subset\ba$ is the Lie subalgebra generated by $\oplus_{i\in I}\CA_i
\subset$ Prim$(\CA)$.

\subsection{}
\label{Serre}
We are going to prove that $\fa=\ba\cong\fn^L_+$ ---
the nilpotent subalgebra of the Langlands dual Lie algebra $\fg^L$.
We choose generators $e_i\in\CA_i-0$, and start with the following Proposition:

{\bf Proposition.} The generators $e_i$ satisfy the Serre relations
of $\fn^L_+$.

{\em Proof.} The corresponding commutator has the weight lying out of
$\cR^+$, and hence vanishes by ~\ref{prim} ~a). $\Box$

\subsubsection{Remark}
\label{choice}
In fact, there is a preferred choice of generators
$e_i$. Note that $\CZ^i\cong\BA^2$, and $\del_i\CZ^i\cong\BA^1$.
So $\IC^i$ contains a subsheaf $\ul\BZ[2]$, and $H^0(\BA^1,\fs_i^!\IC^i)\supset
H^0(\BA^1,\fs_i^!\ul\BZ[2])=\BZ$. The element $1\in\BZ$ in the RHS corresponds
to the canonical element $e_i$ in the LHS.
We will use this choice of $e_i$ in what follows.

\subsection{}
\label{smirnoff}
We construct the action of $\CA$ on
$\oplus_{\alpha\in\BN[I]}H^\bullet(\CQ^\alpha,\IC(\CQ^\alpha))$ closely
following ~\ref{rest}--\ref{stab}.

\subsubsection{} First we construct the {\em stabilization map}
$\ft_{\alpha,\beta}:\ \CA_\alpha\lra\Ext^{|\alpha|}_{\CQ^{\alpha+\beta}}
(\IC(\del_\alpha\CQ^{\alpha+\beta}),\IC(\CQ^{\alpha+\beta}))$
where $\del_\alpha\CQ^{\alpha+\beta}$ stands for the image of the twisting
map $\sigma_{\beta,\alpha}:\ \CQ^\beta\times C^\alpha\lra\CQ^{\alpha+\beta}$
(see ~\cite{fm} ~3.4.1).

To this end we note that the sheaf $\fs_\alpha^!\IC^\alpha$ on
$\BA^\alpha$ extends smoothly to the sheaf $\fF_\alpha$ on $C^\alpha$,
and $H^0(\BA^\alpha,\fs_\alpha^!\IC^\alpha)=H^0(C^\alpha,\fF_\alpha)=
\Ext^{|\alpha|}_{C^\alpha}(\IC(C^\alpha),\fF_\alpha)$. 

Furthermore, tensoring with $\Id\in\Ext^0(\IC(\CQ^\beta),\IC(\CQ^\beta))$,
we get the map from $\CA_\alpha=\Ext^{|\alpha|}_{C^\alpha}
(\IC(C^\alpha),\fF_\alpha)$ to 
$\Ext^{|\alpha|}_{\CQ^\beta\times C^\alpha}
(\IC(\CQ^\beta)\boxtimes\IC(C^\alpha),\IC(\CQ^\beta)\boxtimes
\fF_\alpha)$.

Now the same argument as in ~\ref{bezruk} shows that there is a unique
morphism $c$ from $\IC(\CQ^\beta)\boxtimes\fF_\alpha$
to $\sigma_{\beta,\alpha}^!\IC(\CQ^{\alpha+\beta})$ extending the factorization
isomorphism from an open part $\oQ^\beta\times\BA^\alpha\subset
\CQ^\beta\times C^\alpha$. Thus $c$ induces the map from 
$\Ext^{|\alpha|}_{\CQ^\beta\times C^\alpha}
(\IC(\CQ^\beta)\boxtimes\IC(C^\alpha),\IC(\CQ^\beta)\boxtimes\fF_\alpha)$ to 
$\Ext^{|\alpha|}_{\CQ^\beta\times C^\alpha}
(\IC(\CQ^\beta)\boxtimes\IC(C^\alpha),
\sigma_{\beta,\alpha}^!\IC(\CQ^{\alpha+\beta}))$.

The latter space equals $\Ext^{|\alpha|}_{\CQ^{\alpha+\beta}}
((\sigma_{\beta,\alpha})_*(\IC(\CQ^\beta)\boxtimes\IC(C^\alpha)),
\IC(\CQ^{\alpha+\beta}))=\Ext^{|\alpha|}_{\CQ^{\alpha+\beta}}
(\IC(\del_\alpha\CQ^{\alpha+\beta}),\IC(\CQ^{\alpha+\beta}))$.

Finally, we define $\ft_{\alpha,\beta}$ as the composition of above maps.

\subsubsection{} 
\label{bezr}
Second, we construct the {\em costabilization map}
$$\fx_{\beta,\alpha}:\ H^\bullet(\CQ^\beta,\IC(\CQ^\beta))\lra
H^{\bullet-|\alpha|}(\del_\alpha\CQ^{\alpha+\beta},
\IC(\del_\alpha\CQ^{\alpha+\beta}))$$
To this end we note that exactly as in ~\ref{normal}, we have
$\IC(\del_\alpha\CQ^{\alpha+\beta})=(\sigma_{\beta,\alpha})_*
\IC(\CQ^\beta\times C^\alpha)=(\sigma_{\beta,\alpha})_*
(\IC(\CQ^\beta)\boxtimes\ul\BC[|\alpha|])$.
Let $[C^\alpha]\in H^{-|\alpha|}(C^\alpha,\ul\BC[|\alpha|])$ denote the
fundamental class of $C^\alpha$.

Now, for $h\in H^\bullet(\CQ^\beta,\IC(\CQ^\beta))$ we define
$\fx_{\beta,\alpha}(h)$ as $h\otimes[C^\alpha]\in
H^{\bullet-|\alpha|}(\CQ^\beta\times C^\alpha,\IC(\CQ^\beta)\boxtimes
\ul\BC[|\alpha|])=H^{\bullet-|\alpha|}(\del_\alpha\CQ^{\alpha+\beta},
(\sigma_{\beta,\alpha})_*(\IC(\CQ^\beta)\boxtimes\ul\BC[|\alpha|]))=
H^{\bullet-|\alpha|}(\del_\alpha\CQ^{\alpha+\beta},
\IC(\del_\alpha\CQ^{\alpha+\beta}))$.

\subsection{Definition}
Let $a\in\CA_\alpha,\ h\in H^\bullet(\CQ^\beta,\IC(\CQ^\beta))$. We define
the action $a(h)\in H^\bullet(\CQ^{\alpha+\beta},\IC(\CQ^{\alpha+\beta}))$
as the action of $\ft_{\alpha,\beta}(a)$ on the global cohomology applied
to $\fx_{\beta,\alpha}(h)$.

Let us stress that the action of $\CA_\alpha$ {\em preserves cohomological
degrees.}

\subsubsection{} For $a\in\CA_\alpha,\ b\in\CA_\beta,\
h\in H^\bullet(\CQ^\beta,\IC(\CQ^\beta))$ we have $a(b(h))=a\cdot b(h)$.
The proof is entirely similar to the proof of associativity of the
multiplication in $\CA$.

\subsection{}
\label{F}
For $\beta\in\BN[I],\ i\in I$ both
$H^\bullet(\CQ^\beta,\IC(\CQ^\beta))$ and
$H^\bullet(\CQ^{\beta+i},\IC(\CQ^{\beta+i}))$ are Poincar\'e selfdual.
We define the map
$$f_i:\ H^\bullet(\CQ^{\beta+i},\IC(\CQ^{\beta+i}))\lra
H^\bullet(\CQ^\beta,\IC(\CQ^\beta))$$ as the dual of the map
$$e_i:\ H^\bullet(\CQ^\beta,\IC(\CQ^\beta))\lra
H^\bullet(\CQ^{\beta+i},\IC(\CQ^{\beta+i}))$$
(see ~\ref{choice}).

It follows from ~\ref{Serre} that the maps $f_i,\ i\in I$, satisfy the Serre
relations of $\fn_-^L$.

\subsubsection{}
\label{f}
Let us sketch a more explicit construction of the operators $f_i$. 
Let $\sigma$ denote the closed embedding $\del_i\CQ^{\beta+i}\hookrightarrow
\CQ^{\beta+i}$, and let $j$ denote the embedding of the complementary open
subset $U$. We have an exact sequence of perverse sheaves
$$0\to\sigma^*j_{!*}\IC(U)[-1]\to j_!\IC(U)\to\IC(\CQ^{\beta+i})\to0$$
defining an element $\ff'_i$ in 
$\Ext^1(\IC(\CQ^{\beta+i}),\sigma^*j_{!*}\IC(U)[-1]$. 
The weight considerations as in ~\ref{bezruk} show that 
$\sigma^*j_{!*}\IC(U)[-1]$ canonically surjects to $\IC(\del_i\CQ^{\beta+i})$,
and thus $\ff'_i$ gives rise to $\ff_i\in\Ext^1(\IC(\CQ^{\beta+i}),
\IC(\del_i\CQ^{\beta+i})$.
The action of $\ff_i$ on the global cohomology defines the same named map
$\ff_i:\ H^\bullet(\CQ^{\beta+i},\IC(\CQ^{\beta+i}))\lra
H^{\bullet+1}(\del_i\CQ^{\beta+i},\IC(\del_i\CQ^{\beta+i}))$.
Now recall that
$H^\bullet(\del_i\CQ^{\beta+i},\IC(\del_i\CQ^{\beta+i}))=
H^\bullet(\CQ^\beta,\IC(\CQ^\beta))\otimes H^\bullet(C,\ul\BC[1])$.
So the projection of $H^\bullet(C,\ul\BC[1])=\BC[1]\oplus\BC[-1]$ to $\BC[-1]$
defines the projection $p:\
H^{\bullet+1}(\del_i\CQ^{\beta+i},\IC(\del_i\CQ^{\beta+i}))\lra
H^\bullet(\CQ^\beta,\IC(\CQ^\beta))$.

Finally, for $h\in H^\bullet(\CQ^{\beta+i},\IC(\CQ^{\beta+i}))$ we have
$f_i(h)=p\ff_i(h)\in H^\bullet(\CQ^\beta,\IC(\CQ^\beta))$.

\subsubsection{}
\label{e}
We leave to the reader the absolutely similar elementary
construction of the operators $e_i$. We only mention that the corresponding
local element $\fe_i\in
\Ext^1_{\CQ^{\beta+i}}(\IC(\del_i\CQ^{\beta+i}),\IC(\CQ^{\beta+i}))$
comes from the star extension of the constant sheaf on an open subset of
$\CQ^{\beta+i}$.

\subsection{Proposition}
\label{fe}
If $i,j\in I,\ i\not=j$, then $e_if_j=f_je_i:\
H^\bullet(\CQ^\beta,\IC(\CQ^\beta))\lra
H^\bullet(\CQ^{\beta+i-j},\IC(\CQ^{\beta+i-j}))$.

{\em Proof.} Using the elementary constructions ~\ref{e}, ~\ref{f} of
$e_i$ and $f_j$ we reduce the Proposition to the local calculation in a
smooth open subset $U$ of $\CQ^{\beta+i}$. This local calculation is nothing
else than the following fact. Let $D=D_1\cup D_2$ be a divisor with normal
crossings in $U$, consisting of two smooth irreducible components.
If we shriek extend the constant sheaf on $U-D$ through $D_1$ and then
star extend it through $D_2$ we get the same result if we first star extend
the constant sheaf on $U-D$ through $D_2$, and then shriek extend it through
$D_1. \quad \Box$

\subsection{}
\label{h}
For $i\in I$ we define the endomorphism $h_i$ of
$H^\bullet(\CQ^\beta,\IC(\CQ^\beta))$ as the scalar multiplication by
$\langle\beta+2\crho,i'\rangle$ where $i'\in X$ is the simple root.

{\bf Proposition.} $e_if_i-f_ie_i=h_i:\
H^\bullet(\CQ^\beta,\IC(\CQ^\beta))\lra H^\bullet(\CQ^\beta,\IC(\CQ^\beta))$.

{\em Proof.}
Let $U\subset\CQ^{\beta+i}$ be an open subset such that
$W:=U\cap\del_i\CQ^{\beta+i}$ consists of quasimaps of defect {\em exactly} $i$.
Then $W=\oQ^\beta\times C$ is a smooth divisor in $U$. Let $T_WU$ be the
normal bundle to $W$ in $U$. This is a line bundle on $W$.
For $\phi\in\oQ^\beta$ we can restrict $T_WU$ to $C=\phi\times C\subset
\oQ^\beta\times C$, and the degree of the restriction does not depend on a
choice of $\phi$. We denote this degree by $d_{\beta,i}$.

Using the elementary construction of $e_i,f_i$ we see that
the commutator $[e_i,f_i]$ acts as the scalar multiplication by $d_{\beta,i}$.

The equality $d_{\beta,i}=\langle\beta+2\crho,i'\rangle$ is proved
in ~\cite{fk} ~4.4. Let us make a few comments on this proof. As stated, it
works in the case $\bG=SL_n$, moreover, it uses the Laumon resolution
$\CQ^\alpha_L\lra\CQ^\alpha$. In fact, the proof uses not the whole
space $\CQ^\alpha_L$ but the open subspace $\CU_\alpha^L\subset\CQ^\alpha_L$
(see {\em loc. cit.} ~4.2.1). This open subset projects isomorphically into
$\CQ^\alpha$. More precisely, $\CU_\alpha^L\subset\CQ^\alpha$ consists of
all quasimaps with defect {\em at most a simple coroot}. Now one can see that
the calculations in $\CU_\alpha^L$ carried out in ~\S4 of {\em loc. cit.}
do not use any specifics of $SL_n$ and carry over to the case of arbitrary
$\bG$. This completes the proof of the Proposition. $\Box$

\subsection{}
\label{!}
Combining the results of ~\ref{Serre}, ~\ref{F}, ~\ref{fe},
~\ref{h} we see that $e_i,f_i,h_i,\ i\in I$, generate the action of the
Langlands dual Lie algebra $\fg^L$ on
$\oplus_{\alpha\in\BN[I]}H^\bullet(\CQ^\alpha,\IC(\CQ^\alpha))$.
Since $\fg^L$ is a {\em simple} Lie algebra, we obtain an {\em embedding}
$\fg^L\hookrightarrow$End
$(\oplus_{\alpha\in\BN[I]}H^\bullet(\CQ^\alpha,\IC(\CQ^\alpha)))$.
The image of $\fn_+^L\subset\fg^L$ under this embedding coincides with
the image of the Lie algebra $\fa\subset\CA$ (see ~\ref{a}) in
End$(\oplus_{\alpha\in\BN[I]}H^\bullet(\CQ^\alpha,\IC(\CQ^\alpha)))$.
We deduce that $\fa\cong\fn_+^L$. Hence we obtain the embedding of
the enveloping algebras $U(\fn_+^L)\cong U(\fa)\hookrightarrow\CA$.
Comparing the graded dimensions we see that this embedding is in fact an
isomorphism. In particular, Prim$(\CA)=\fa\cong\fn_+^L$.

We have proved the following Theorem.

{\bf Theorem.} The bialgebra $\CA$ is isomorphic to the universal enveloping
algebra $U(\fn_+^L)$.

\section{Closures of Schubert strata and combinatorics of alcoves}

\subsection{}
\label{fuck}
We recall some combinatorics from ~\cite{l1} and ~\cite{l4}.
Let $E$ be an $\BR$-vector space $Y\otimes_\BZ\BR$. We equip it with a
scalar product $(|)$ extending it by linearity from the basis $I:\
(i|j)\df\ d\frac{i\cdot j}{d_id_j}$. Here $i\cdot j$ is a part of Cartan datum
(see ~\cite{l}, 1.1), $d_i\df\ \frac{i\cdot i}{2}$, and
$d\df\ \operatorname{max}_{i\in I}d_i$.

The Weyl group $\CW_f$ acts on $E$ by orthogonal reflections; it is generated
by the reflections $s_i(y)\df\ y-\langle y,i'\rangle i$ (notations of
{\em loc. cit.}, 2.2; we have extended the pairing $\langle,\rangle:\
Y\times X\lra\BZ$ by linearity to $\langle,\rangle:\ E\times X\lra\BR$.)
We shall regard $\CW_f$ as acting on $E$ {\em on the right}.
The set $I\CW_f$ is the set $\check\CR$ of {\em coroots}, and the set
$I\CW_f\cap\BN[I]$ is the set $\check\CR^+$ of {\em positive coroots}.
Recall that $\crho=\frac{1}{2}\sum_{\ctheta\in\check\CR^+}\ctheta\in E$.

We consider the following collection $$\fF\df\ \{H_{\theta,n},\ \theta\in\CR,\
n\in\BZ\}$$ of affine hyperplanes in $E$:
$$H_{\theta,n}=\{y\in E | \langle y,\theta\rangle=n\}$$
Each $H\in\fF$ defines an orthogonal reflection $y\mapsto y\sigma_H$ in $E$
with fixed point set $H$. Let $\Omega$ be the group of affine motions
generated by the $\sigma_H\ (H\in\fF)$. We shall regard $\Omega$ as acting
on the right on $E$.

The connected components of $E-\cup_{H\in\fF}$ are called {\em alcoves}.
The group $\Omega$ acts simply transitively on the set $\fA$ of alcoves.
We shall denote by $A^-_0$ the following alcove:
$$A^-_0\df\ \{y\in E | \langle y,i'\rangle<0\ \forall i\in I;\
\langle y,\theta_0\rangle>-1\}$$
where $\theta_0\in\CR^+$ is the highest root.
Then $\Omega$ is generated by the reflections in the walls of $A^-_0$.
The subgroup of $\Omega$ generated by the reflections in the walls $H_{i',0}$
is just the Weyl group $\CW_f$.

Note that in ~\cite{l4}, \S\S1,8 the group $\Omega$ is called the
{\em affine Weyl group} and denoted by $W$ (not to be confused with the
{\em dual affine Weyl group} $W^\sharp$). We will follow the notations of
~\cite{l1} instead. In particular, we will call $\CW$ the group defined as
in {\em loc. cit.}, 1.1. It acts on $\fA$ simply transitively {\em on the left},
commuting with the action of $\Omega$. As a Coxeter group it is canonically
isomorphic to the affine Weyl group $\Omega$. In the notations of
{\em loc. cit.} the subgroup $\CW_f\subset\Omega$ is nothing else then
$\Omega_0$, and $w_0=\omega_0$.

Recall that the set $S$ of Coxeter generators of $\CW$ can be represented as
the faces (walls) of $A^-_0$. The generator corresponding to the wall
$H_{i',0},\ i\in I,$ will be denoted $s_i$, and the generator corresponding
to the wall $H_{\theta_0,-1}$ will be denoted by $s_0$.

It is easy to see that the intersection of $\Omega$ with the group of
translations of $E$ is (the group of translations by vectors in) $Y$.
Thus we obtain a normal subgroup $Y\subset\Omega$. It is known that $\Omega$
is a semidirect product of $\CW_f$ and $Y$ (see ~\cite{l4}, 1.5).
In particular, each element $\omega\in\Omega$ can be uniquely written in the
form $\omega=w\chi,\ w\in \CW_f,\ \chi\in Y$. Combining this observation with
the action of $\Omega$ on $\fA$ we obtain a bijection
$$\CW_f\times Y\iso\fA:\ (w,\chi)\mapsto A^-_0w\chi$$
We will use this bijection freely, so we will write $(w,\chi)$ for an alcove
often.

We will denote by $\beta_0\in\check\CR^+$ the coroot dual to the highest
root $\theta_0$, and by $s_{\beta_0}\in \CW_f$  the
reflection in $\CW_f$ taking $y\in E$ to $y-\langle y,\theta_0\rangle\beta_0$.

\subsubsection{Lemma}
\label{s0}
For any $i\in I$ we have $s_i(w,\chi)=(s_iw,\chi)$, and
$s_0(w,\chi)=(s_{\beta_0}w,\chi-\beta_0w)$.

{\em Proof.} Clear. $\Box$

\subsubsection{}
\label{dim}
Recall that for a pair of alcoves $A,B\in\fA$ the {\em distance} $d(A,B)$
was defined in ~\cite{l1}, 1.4.
For $\chi=\sum_{i\in I}a_ii\in Y$ we define $|\chi|\in\BZ$ as follows:
$$|\chi|\df\ \sum_{i\in I}a_i$$
For $w\in\CW_f$ let $l(w)$ denote the usual length function:
$$l(w)=\sharp(\check\CR^+w\cap-\check\CR^+)=\dim\bX_w$$

{\em Lemma.} Let $B=(w,\chi)$. Then $d(A^-_0,B)=2|\chi|+l(w)$.

{\em Proof.} For $A=(w_1,\chi_1),B=(w_2,\chi_2)$ let us define
$d'(A,B)\df\ 2|\chi_2|-2|\chi_1|+l(w_2)-l(w_1)$. We want to prove
$d'(A,B)=d(A,B)$. To this end it suffices to check the properties
~\cite{l1}(1.4.1),(1.4.2) for $d'$ instead of $d$. This in turn follows
easily from the Lemma ~\ref{s0}. $\Box$

\subsubsection{}
\label{leq}
Recall a few properties of the partial order $\leq$ on $\fA$
introduced in ~\cite{l1}, 1.5.

\begin{equation}
\label{01}
(w_1,\chi)\leq(w_2,\chi)\Leftrightarrow w_1\leq w_2,
\end{equation}
where $w_1\leq w_2$ stands for the usual Bruhat order on $\CW_f$.

\begin{equation}
\label{02}
(w_1,\chi_1)\leq(w_2,\chi_2)\Leftrightarrow
(w_1,\chi_1+\chi)\leq(w_2,\chi_2+\chi)\ \forall\chi\in Y,
\end{equation}

\begin{equation}
\label{03}
(w,0)\leq(ws_{\check\theta},\check{\theta})\
\forall \check{\theta}\in\check\CR^+,
\end{equation}
where $s_{\check\theta}\in\CW_f$ is the reflection in $\CW_f$ taking
$y\in E$ to $y-\langle y,\theta\rangle\check{\theta},\ \theta$ being a root
dual to the coroot $\check{\theta}$.

\begin{equation}
\label{04}
s_0(w,\chi)\leq(w,\chi)\Leftrightarrow\beta_0w>0.
\end{equation}

The proofs are easy if not contained explicitly in ~\cite{l1}.
For the converse we have the following

{\em Lemma.} Let $\preceq$ be the minimal partial order on $\fA$ enjoying
the properties ~(\ref{01}), ~(\ref{02}), ~(\ref{03}).
Then $A\preceq B\Leftrightarrow A\leq B$.

{\em Proof.} Easy. $\Box$

\subsection{}
\label{closure}
Recall the notations of ~\cite{fm} ~8.4, ~9.1.
The Schubert stratification of $\CZ^\alpha_\chi$ reads as follows:
$$\CZ^\alpha_\chi=\bigsqcup_{w\in\CW_f}^{0\leq\beta\leq\alpha}
\dZ^\beta_{w,\chi-\alpha+\beta}$$
or, after a change of variable $\eta=\chi-\alpha+\beta$,
$$\CZ^\alpha_\chi=\bigsqcup_{w\in\CW_f}^{\chi-\alpha\leq\eta\leq\chi}
\dZ^{\eta-\chi+\alpha}_{w,\eta}$$
We want to describe the closure of a stratum
$\dZ^{\eta-\chi+\alpha}_{w,\eta}$.

{\bf Theorem.} Fix a pair of alcoves $A=(w,\eta),B=(y,\xi)$.
For any $\chi\geq\eta,\xi$ and
$\alpha\in\BN[I]$ sufficiently dominant
(such that $\eta-\chi+\alpha\geq10\crho\leq\xi-\chi+\alpha$) we have

a) $\dim\dZ^{\xi-\chi+\alpha}_{y,\xi}-\dim\dZ^{\eta-\chi+\alpha}_{w,\eta}=
d(A,B)$;

b) $\dZ^{\eta-\chi+\alpha}_{w,\eta}$ lies in the closure of
$\dZ^{\xi-\chi+\alpha}_{y,\xi}$ iff $A\leq B$.

{\em Proof.} a) follows immediately comparing the Lemmas ~\ref{dim}
and ~\cite{fm} ~8.5.2.

The proof of b) occupies the rest of this section.

\subsection{}
Let $\th\in\CR^+$ be a positive root. Let $\fg_\th$ denote
the corresponding ${\frak{sl}}_2$ Lie subalgebra in
$\fg$, and let $\bG_\th$ denote the corresponding
$SL_2$-subgroup in $\bG$. Take any $w\in\CW_f$ and let $y=s_\th w$.
We will view $w,y$ as $\bH$-fixed points in $\bX$. The $\bG_\theta$-orbit
$\bG_\th w=\bG_\th y$ is a smooth
rational curve of degree $\ctheta$ in $\bX$
(here $\check\th$ stands for the dual coroot of the root $\th$).
We will view it as a closed subset of $\bX$ rather than a parametrized curve.
Let us denote this curve by $L_{w,y}$. It is fixed by the $\bH$ action on $\bX$
(as a subset, not pointwise).

\subsubsection{Lemma}
Let $L\subset\bX$ be an irreducible rational curve fixed by the
Cartan action. Then $L=L_{w,y}$ for some $w,y=s_\th w\in\CW_f$.

{\em Proof.}
Since $L$ is a fixed curve it must contain a fixed point $w\in\CW_f$.
Let $f:\hat L\to L\subset\bX$ be the normalization of $L$. The $\bH$-action
on $L$ extends to the $\bH$-action on $\hat L$. We choose an $\bH$-fixed
point in $f^{-1}(w)$ and preserve the name $w$ for this point.
Let $t$ be a formal coordinate at $w\in\hat L$.
The map $f$ gives rise to the
homomorphism $f^*:\widehat\CO_{\bX,w}\to\BC[[t]]=\widehat\CO_{\hat L,w}$
of $X$-graded rings from the completion of the local ring $\CO_{\bX,w}$
to the ring of formal power series in $t$. We have
$\widehat\CO_{\bX,w}=\BC[[T^*_w\bX]]$ --- the ring of formal power
series on the tangent space of $\bX$ at the point $w$, and the grading
is induced by the grading of the cotangent space
$$
T^*_w\bX=\bigoplus\limits_{\vth\in w\CR^+}\BC x_\vth,
$$
where $x_\vth\in\fg^*$ is an $\bH$-eigenvector with the eigenvalue $\vartheta$.
Since the roots in
$w\CR^+$ are pairwise linearly independent, it follows that there
exists a root $\th\in w\CR^+$ and a positive integer $n$ such that the
map $f^*$ is given by
$$
f^*(x_\vth)=\begin{cases}
t^n,    &\text{ if }\vth=\th\\
0,      &\text{ otherwise.}
\end{cases}
$$
Therefore the formal jet of the map $f:\hat L\to\bX$ coincides with the
formal jet of the map $\vphi_{w,y}^n:\BP^1\to\bX$ given by the composition
of the $n$-fold covering $\BP^1\to L_{w,y}$ ramified over the points $w$
and $y$, and of the embedding $L_{w,y}\to\bX$. Now irreducibility of $L$
implies $L=L_{w,y}. \quad \Box$

\subsubsection{Remark}
The fixed curves $L_{w,y}$ and $L_{w',y'}$ intersect nontrivially
iff $w=w'$ or $w=y'$ or $y=w'$ or $y=y'$. In effect, a point of intersection
of two fixed curves has to be a fixed point.

\subsection{Definition}
Let $f:\CC\to\bX$ be a stable map (see ~\cite{ko})
from a genus 0 curve $\CC$ into the flag variety $\bX$.
Let $w,y\in\CW_f$.
If $f(\CC)\cap\bX_w\ne\emptyset$ and $f(\CC)\cap\overline\bX_y\ne\emptyset$
we will say that the pair $(w,y)$ is $(f,\CC)$-connected.
If $\deg\CC=\ga$ we will say also that $(w,y)$ is $\ga$-connected.

\subsubsection{Lemma}
\label{conn}
A pair $(w,y)$ is $\ga$-connected iff there exists a collection
$(\check\th_1,\dots,\check\th_k)$ of positive coroots
$\check\th_r\in\check\CR^+$ such that $\check\th_1+\ldots+\check\th_k\le\ga$
and $s_{\th_k}\dots s_{\th_1}w\le y$.

{\em Proof.}
Suppose the pair $(w,y)$ is $(f,\CC)$-connected, where $f:\CC\to\bX$
is a stable map of degree $\gamma$. Then we have $f(\CC)\cap\bX_w\ne\emptyset$.
Acting on $f$ by an element of the Borel subgroup we can assume that
$w\in f(\CC)$. Acting by the Cartan subgroup $\bH$
we can degenerate $(f,\CC)$ into
an $\bH$-fixed stable map $f':\CC'\to\bX$. Since $w\in f(\CC)$ and $w$ is
an $\bH$-fixed point, we will have $w\in f'(\CC')$. Similarly, since
$f(\CC)\cap\overline\bX_y\ne\emptyset$ and $\overline\bX_y$ is a closed
$\bH$-invariant subspace, we will have $f'(\CC')\cap\overline\bX_y\ne\emptyset$,
hence $f'(\CC')\ni y'$ for some $y'\in\CW_f$ such that $y'\le y$. The image of
an $\bH$-fixed stable map is a connected union of $\bH$-fixed curves, hence
there exists a sequence $L_1,\dots,L_k$ of $\bH$-fixed curves such that
$w\in L_1$, $y'\in L_k$ and $L_r\cap L_{r+1}\ne\emptyset$.
We can assume that $L_r=L_{w_r,w_{r+1}}$
($r=1,\dots,k$), where $w_1=w$, $w_{r+1}=s_{\th_r}w_r$ and $w_{k+1}=y'$.
%$L_1=L_{w,s_{\th_1}w}$, $L_2=L_{s_{\th_1}w,s_{\th_2}s_{\th_1}w},\dots,
%L_k=L_{s_{\th_{k-1}}\dots s_{\th_1}w,s_{\th_k}\dots s_{\th_1}w}$ and
%$s_{\th_k}\dots s_{\th_1}w=y'$.
Since
$$
\ga=\deg f=\deg f'\ge\deg L_1+\ldots+\deg L_k=\check\th_1+\ldots+\check\th_k
$$
the Lemma follows. $\Box$

\subsubsection{Corollary}
\label{codim1}
If $w\not\le y$ and $(w,y)$ is $\ga$-connected then
$$
l(w)\le l(y)+2|\ga|-1
$$
and equality holds only if $\ga=\check\th$, $l(s_\th)=2|\ctheta|-1$ and
$s_\th w=y$ for some $\ctheta\in\cR^+$.

{\em Proof.}
By the Lemma \ref{conn} we have
\begin{multline*}
l(y)\ge l(s_{\th_k}\dots s_{\th_1}w)\ge
l(w)-l(s_{\th_1})-\ldots-l(s_{\th_k})\ge\\
\ge l(w)-(2|\ctheta_1|-1)-\ldots-(2|\ctheta_k|-1)\ge
l(w)-2|\ga|+k\ge l(w)-2|\ga|+1.
\end{multline*}
The Corollary follows. $\Box$

\subsection{Proposition}
\label{criter}
A Schubert stratum $\dZ^{\eta-\chi+\al}_{w,\eta}$ lies in the closure of
$\dZ^{\xi-\chi+\al}_{y,\xi}$ iff $\eta\le\xi$ and the pair
$(w,y)$ is $(\xi-\eta)$-connected.

{\em Proof.}
We start with the ``only if'' part.
Both strata lie in the closed subspace
$\CZ^{\xi-\chi+\al}_\xi\subset\CZ^\al_\chi$, hence we can consider
them as subspaces in $\CZ^{\xi-\chi+\al}_\xi$.

Let $\CQ^{\xi-\chi+\al}_K=
\overline M_{0,0}(\BP^1\times\bX,(1,\xi-\chi+\al))$
be the Kontsevich space of stable maps from the genus zero curves without
marked points to $\BP^1\times\bX$ of bidegree $(1,\xi-\chi+\alpha)$
(see ~\cite{ko}). The map $\pi:\ \CQ^{\xi-\chi+\al}_K\lra\CQ^{\xi-\chi+\al}$
is constructed in the Appendix.
For a stable map $f:\ \CC\lra\BP^1\times\bX$ in $\CQ^{\xi-\chi+\al}_K$ we
denote by $f':\ \CC\lra\BP^1$ (resp. $f'':\ \CC\lra\bX$)
its composition with the first (resp. second) projection.

Let $\CK^{\xi-\chi+\al}_\xi\subset\CQ^{\xi-\chi+\al}_K$
be the preimage of the subspace
$\CZ^{\xi-\chi+\al}_\xi\subset\CQ^{\xi-\chi+\al}$ under the map
$\pi:\CQ_K^{\xi-\chi+\al}\to\CQ^{\xi-\chi+\al}$. It is a locally closed
subspace consisting of all stable maps $f:\CC\to\BP^1\times\bX$
such that ${f'}^{-1}(\infty)$ is a point
(i.e. $\CC$ has no vertical component over $\infty\in\BP^1$), and
$f''({f'}^{-1}(\infty))=w_0\in\bX$. Consider the preimage
of the Schubert stratification of $\CZ^{\xi-\chi+\al}_\xi$:
$$
\CK^{\xi-\chi+\al}_\xi=\bigsqcup_{w\in\CW_f}^{\chi-\al\le\eta\le\chi}
\dK^{\eta-\chi+\al}_{w,\eta}.
$$
Here $\dK^{\eta-\chi+\al}_{w,\eta}$ is just the subspace of all stable
maps $f:\CC\to\BP^1\times\bX$ such that
$$
\deg\CC_1=\xi-\eta\quad\text{and}\quad f''(P)\in\bX_w,
$$
where $\CC_1={f'}^{-1}(0)$
is the vertical component of $\CC$ over the point $0\in\BP^1$, and
$\CC_0$ is the main component of $\CC$ (the one which projects
isomorphically onto $\BP^1$ under $f'$); the point $P$ is the intersection of
these components: $P=\CC_0\cap\CC_1$;
and as always, $f''({f'}^{-1}(\infty))=w_0$.
It follows that $\dZ^{\eta-\chi+\al}_{w,\eta}$
lies in the closure of $\dZ^{\xi-\chi+\al}_{y,\xi}$ iff the map
$$
\pi:\dK^{\eta-\chi+\al}_{w,\eta}\bigcap
\overline{\dK^{\xi-\chi+\al}_{y,\xi}}
\to\dZ^{\eta-\chi+\al}_{w,\eta}
$$
is surjective. In particular, the above intersection must be non-empty.

Consider the subspace
$\CK^{\xi-\chi+\al}_{y,\xi}\subset\CK^{\xi-\chi+\alpha}_\xi$ formed by the
stable maps $f:\CC\to\BP^1\times\bX$ such that
$f(\CC)\cap(\{0\}\times\overline\bX_y)\ne\emptyset$. In other words,
$\CK^{\xi-\chi+\al}_{y,\xi}\subset\CK^{\xi-\chi+\alpha}_\xi$ is formed by the
stable maps such that
$f''(\CC_1)\cap\overline\bX_y\ne\emptyset$. The subspace
$\CK^{\xi-\chi+\al}_{y,\xi}$ is obviously closed and irreducible.
It contains $\dK^{\xi-\chi+\al}_{y,\xi}$ as an open subspace, hence
$$
\CK^{\xi-\chi+\al}_{y,\xi}=\overline{\dK^{\xi-\chi+\al}_{y,\xi}}.
$$
The intersection $\dK^{\eta-\chi+\al}_{w,\eta}\bigcap
\overline{\dK^{\xi-\chi+\al}_{y,\xi}}$ consists of all stable maps
such that
$$
\deg\CC_1=\xi-\eta,\qquad f''(P)\in\bX_w\quad\text{and}\quad
f''(\CC_1)\cap\overline\bX_y\ne\emptyset.
$$
Therefore the intersection is non-empty only if the pair $(w,y)$
is $(\xi-\eta)$-connected (recall that $P\in\CC_1$). This completes the proof
of the ``only if'' part.

On the other hand, if $(w,y)$ is $(\xi-\eta)$-connected, then
for any stable map $(f,\CC)$ we can replace the component
$\CC_1$ by a stable curve of degree $(\xi-\eta)$, connecting the point
$f''(P)\in\bX_w$ with $\overline\bX_y$. Such replacement will not affect
$\pi(f,\CC)$, but the new curve will lie in the intersection
$\dK^{\eta-\chi+\al}_{w,\eta}\bigcap\overline{\dK^{\xi-\chi+\al}_{y,\xi}}$,
hence the above intersection maps surjectively onto
$\dZ^{\eta-\chi+\al}_{w,\eta}$ and the Proposition follows. $\Box$

\subsection{}
Now we can prove the Theorem ~\ref{closure}.

We define the {\em adjacency order} $\preceq$ on $\fA$ as follows.
We say that $A=(w,\eta)\preceq(y,\xi)=B$ if for any $\chi\geq\eta,\xi$
and sufficiently dominant $\alpha\in\BN[I]$ the stratum
$\dZ^{\eta-\chi+\alpha}_{w,\eta}$ lies in the closure of
$\dZ^{\xi-\chi+\alpha}_{y,\xi}$.

First we check that $\preceq$ satisfies the relations ~(\ref{01}--\ref{03})
of ~\ref{leq}.
To this end let us rephrase these properties via criterion \ref{criter}.
The property ~(\ref{01}) means that $(w_1,w_2)$ is $0$-connected iff
$w_1\le w_2$. The property ~(\ref{02}) expresses the fact that the adjacency
of strata depends on difference of degrees only, which is evident from the
criterion \ref{criter}. The property ~(\ref{03})
means that the pair $(w,s_\th w)$
is $\check\th$-connected. This is indeed so since the curve $L_{w,s_\th w}$
has degree $\check\th$ and connects $w$ with $s_\th w$.

Now the adjacency order $\preceq$ is clearly generated by the adjacencies
in codimension 1. So it remains to check that the only adjacencies in
codimension 1 are the ones of type ~(\ref{03}) or type ~(\ref{01})
with $l(w_1)=l(w_2)-1$.

Assume that $\dZ^{\eta-\chi+\al}_{w,\eta}$
lies in the closure of $\dZ^{\xi-\chi+\al}_{y,\xi}$ and has codimension 1.
Since
$$
\dim\dZ^{\eta-\chi+\al}_{w,\eta}=2|\eta-\chi+\al|+l(w)-\dim\bX,\quad
\dim\dZ^{\xi-\chi+\al}_{y,\xi}=2|\xi-\chi+\al|+l(y)-\dim\bX
$$
it follows that $l(w)=l(y)+2|\xi-\eta|-1$. On the other hand, by
the criterion \ref{criter} the pair $(w,y)$ is $(\xi-\eta)$-connected.
If $\xi=\eta$ then $l(w)=l(y)-1$ and we are in the situation of type
~(\ref{01}). If $\xi>\eta$ then $w\not\le y$ and by the Corollary
\ref{codim1} it follows that $\xi-\eta=\check\th$ and $y=s_\th w$
for some $\th\in\CR^+$, i.e. we are in the situation of type ~(\ref{03}).

The Theorem follows. $\Box$

\subsection{}
\label{who}
Recall that $\dQ^\gamma\subset\CQ^\gamma$ is the open subset formed by all
the quasimaps without defect at $0\in C$ (i.e. defined at $0\in C$, see
~\cite{fm} ~8.1). For $w\in\CW_f$ we define the locally closed subset
$\dQ^\gamma_w\subset\dQ^\gamma$ formed by all the quasimaps $\phi$ such
that $\phi(0)\in\bX_w\subset\bX$ (cf. {\em loc. cit.} ~8.4). Its closure
will be denoted by $\CQ^\gamma_w$; it coincides with the closure of the
fine Schubert stratum $\oQ^\gamma_w$ defined in {\em loc. cit.} ~8.4.1.

{\bf Corollary of the Proof.} Let $\gamma\leq\beta\in\BN[I],\ w,y\in\CW_f$.
Then $\CQ^\beta_y\supset\CQ^\gamma_w$ iff $(y,\beta)=:B\geq A:=(w,\gamma)$.
Also, $\dim\CQ^\beta_y-\dim\CQ^\gamma_w=d(A,B). \quad \Box$

\subsection{}
\label{t'ma}
We will denote the closure of a fine Schubert stratum $\oZ^\beta_y$
by $\CZ^\beta_y$ (note that it was denoted by $\oCZ^\beta_y$ in
~\cite{fm}. We hope this will cause no confusion.)
Similarly, the closure of $\oZ^\beta_{y,\xi}$ (see {\em loc. cit.} ~9.4)
will be denoted by $\CZ^\beta_{y,\xi}$.

We promote ``alcovic'' notations. Suppose we are in the situation of the
above theorem, i.e. $A=(w,\eta)\leq(y,\xi)=B;\ \chi\geq\eta,\xi$, and
$\alpha\in\BN[I]$ is sufficiently dominant. We introduce new variables
$\gamma=\eta-\chi+\alpha$ and $\beta=\xi-\chi+\alpha$
(note that both $\gamma$ and $\beta$ necessarily lie in $\BN[I]$).
Then the above Theorem claims that $\dZ^\gamma_{w,\eta}$ lies in the
closure $\CZ^\beta_{y,\xi}$ of $\dZ^\beta_{y,\xi}$. We will denote this closure
by $\CZ^\beta_B$; thus we have the closed embedding
$\CZ^\gamma_A\hra\CZ^\beta_B$.

With our new notations, it is only natural to rename the snop $\CL(w,\eta)$
(see {\em loc. cit.} ~9.4) into $\CL(A)$, and $\CM(w,\eta)$ into
$\CM(A)$, and $\CalD\CM(w,\eta)$ into $\CalD\CM(A)$.

\subsubsection{}
\label{mrak}
Recall (see {\em loc. cit.} ~9.3) that we call $\alpha\in\BN[I]$ sufficiently
dominant, and write $\alpha\gg0$, iff $\alpha\geq10\crho$. For such
$\alpha$ we defined the open subvariety $\ddZ^\alpha\subset\CZ^\alpha$.
For $\beta\leq\alpha$ the closed subvariety $\CZ^\beta_y\cap\ddZ^\alpha\subset
\ddZ^\alpha$ is nonempty iff $\beta\gg0$. In this case it
will be denoted by the same symbol $\CZ^\beta_y$; we hope this
will cause no confusion.
Again, for $0\ll\gamma\leq\beta$ we have
$\CZ^\beta_y\supset\CZ^\gamma_w$ iff $(y,\beta)=:B\geq A:=(w,\gamma)$.
Also, $\dim\CZ^\beta_y-\dim\CZ^\gamma_w=d(A,B).$

\section{Mixed snops and convolution}

\subsection{}
In this section we compute the stalks of irreducible snops. First we prove
that they carry a natural Tate Hodge structure.
We start with preliminary results about affine Grassmannian.
Recall the stratification of the affine Grassmannian $\CG$ by the Iwahori
orbits $\CG_{w,\eta}$ described e.g. in ~\cite{fm} ~10.4.
The orbit closure $\oCG_{w,\eta}$ is partitioned into its intersections
with semiinfinite orbits: $\oCG_{w,\eta}=
\sqcup_{\alpha\in Y}(\oCG_{w,\eta}\cap T_\alpha)$.
The irreducible perverse sheaf $\IC(\oCG_{w,\eta})$ has the mixed Hodge
module counterpart $\ic(\oCG_{w,\eta})$.

{\bf Proposition.} a) $H^\bullet_c(\oCG_{w,\eta}\cap T_\alpha,
\ic(\oCG_{w,\eta}))$ is a pure Hodge complex of weight $\dim\oCG_{w,\eta}$.
It is a direct sum of Tate Hodge structures;

b) $H^\bullet_c(\oCG_{w,\eta}\cap\ol{T}_\alpha,
\ic(\oCG_{w,\eta}))$ is a pure Hodge complex of weight $\dim\oCG_{w,\eta}$.
It is a direct sum of Tate Hodge structures.

{\em Proof.} a) Let us denote the locally closed embedding of
$\oCG_{w,\eta}\cap T_\alpha$ (resp.
$\oCG_{w,\eta}\cap S_\alpha$) into $\oCG_{w,\eta}$ by $i_T$ (resp. $i_S$).
One can construct a natural isomorphism
$H^\bullet(\oCG_{w,\eta}\cap S_\alpha,i_S^!\ic(\oCG_{w,\eta}))\iso
H^\bullet_c(\oCG_{w,\eta}\cap T_\alpha,i_T^*\ic(\oCG_{w,\eta}))$ 
(see ~\cite{mv}, Opposite parabolic restrictions).
On the other hand,
$\ic(\oCG_{w,\eta})$ is a pure Hodge module of weight $\dim\oCG_{w,\eta}$,
while $H^\bullet_c(?)$ and $i_T^*$ decrease weights, and
$H^\bullet(?)$ and $i_S^!$ increase weights.

It remains to prove that $H^\bullet_c(\oCG_{w,\eta}\cap T_\alpha,
\ic(\oCG_{w,\eta}))$ is a direct sum of Tate Hodge structures. To this end
consider the Cousin spectral sequence associated with the partition into
locally closed subsets $\oCG_{w,\eta}=
\sqcup_{\alpha\in Y}(\oCG_{w,\eta}\cap T_\alpha)$. We have
$E_2^{p,q}\Longrightarrow H^{p+q}(\oCG_{w,\eta},\ic(\oCG_{w,\eta}))$, and
$E_2^{p,q}=\oplus
H^{p+q}_c(\oCG_{w,\eta}\cap T_\alpha,i_T^*\ic(\oCG_{w,\eta}))$; the sum
is taken over $\alpha$ such that the codimension of
$\oCG_{w,\eta}\cap T_\alpha$ in $\oCG_{w,\eta}$ equals $q$. For the weight
reasons the differentials in this spectral sequence vanish, and it collapses
at the second term. Hence
$H^\bullet(\oCG_{w,\eta},\ic(\oCG_{w,\eta}))\cong\oplus_{\alpha\in Y}
H^\bullet_c(\oCG_{w,\eta}\cap T_\alpha,i_T^*\ic(\oCG_{w,\eta}))$. On the
other hand, the LHS is well known to be a direct sum of Tate Hodge structures.

So a) is proved. Now b) follows by the application of Cousin spectral
sequence associated with the partition into the locally closed subsets
$\oCG_{w,\eta}\cap\ol{T}_\alpha=\sqcup_{\beta\geq\alpha}
(\oCG_{w,\eta}\cap T_\beta). \quad \Box$

\subsubsection{Remark}
\label{zhopa}
Let $g\in\bG$ lie in the normalizer $N(\bH)$ of the Cartan subgroup $\bH$,
and let $T_\alpha^g$ denote $g(T_\alpha)$: the action of $g$
on the semiinfinite orbit in $\CG$. The same argument as above proves that
$H^\bullet_c(\oCG_{w,\eta}\cap T_\alpha^g,
\ic(\oCG_{w,\eta}))$ is a pure Hodge complex of weight $\dim\oCG_{w,\eta}$;
it is a direct sum of Tate Hodge structures.

\subsection{}
\label{purity}
According to the Theorem ~\ref{Hodge}, the simple stalk
$\fs_\alpha^*\iota_\alpha^*\ic^\alpha=\ic^0_{\{\{\alpha\}\}}$ is pure of
weight $2|\alpha|$, and it is a direct sum of Tate Hodge structures.
By factorization, the Hodge module $\ic^\alpha$ on $\CZ^\alpha$
is pointwise pure, and all its stalks are direct sums of Tate Hodge structures.
Hence the Hodge module $\ic(\CQ^\alpha)$ on $\CQ^\alpha$
is pointwise pure, and all its stalks are direct sums of Tate Hodge structures.
In what follows, such Hodge modules will be called {\em pointwise pure Tate
Hodge modules}.

Recall the proalgebraic variety $\fQ^\alpha$ introduced in ~\cite{fm} ~10.6.
The theory of mixed Hodge modules on such varieties is developed in ~\cite{kt1}.
In particular, we will be interested in the irreducible Hodge module
$\ic(\fQ^\alpha)$, cf. ~\cite{fm} ~10.7.3. The same argument as in
{\em loc. cit.} ~12.7 proves that $\ic(\fQ^\alpha)$ is a pointwise pure Tate
Hodge module.

For the future reference, let us collect the above facts into a Theorem.

{\bf Theorem}. The Hodge modules $\ic^\alpha,\ic(\CQ^\alpha),\ic(\fQ^\alpha)$
are the pointwise pure Tate Hodge modules.

\subsection{Theorem}
\label{finkel}
$\ic(\CQ^\alpha_w)$ (see ~\ref{who})
is a pointwise pure Tate Hodge module.

{\em Proof.}
The usual argument with factorization reduces the proof to the study of
the stalk of $\ic(\CQ^\alpha_w)$ at the fine Schubert stratum
$\oQ^0_y\times(\BP^1-0)^\beta_\Gamma\times0^{\alpha-\beta}
\subset\CQ^\alpha,\ \beta\leq\alpha,\ y\in\CW_f$
(see {\em loc. cit.} ~8.4.1). Moreover, by the same factorization argument,
this stalk will not change if we add $\gamma\in\BN[I]$ to both $\alpha$ and
$\beta$. We will choose $\gamma$ sufficiently dominant to ensure that
$\alpha+\gamma\in Y^+$. So we may and will suppose that $\alpha$ is
dominant.
Furthermore, we can choose any point in the stratum
$\oQ^0_y\times(\BP^1-0)^\beta_\Gamma\times0^{\alpha-\beta}$,
and we will choose a point $\phi=(\fL_\lambda)_{\lambda\in X^+}\in
\oQ^\beta_y$ where $\fL_\lambda=\CV_\lambda^{\bN_+^y}(\langle\beta-\alpha,
\lambda\rangle\cdot0-\langle D,\lambda\rangle)$, and
$\bN_+^y$ is the nilpotent subgroup conjugated to $\bN_+$ by any representative
$\dot{y}$ of $y$ in the normalizer of $\bH$ (see {\em loc. cit.} ~3.3);
$D\in(\BP^1-0)^\beta_\Gamma$.

Recall the setup of the Theorem ~12.12 of {\em loc. cit.}
We consider the Hodge module $\CF=\ic(\oCG_{w,\alpha})$ on $\oCG_\alpha$
(recall that we assumed $\alpha\in Y^+$), and
the convolution $\bc^0_\CQ(\CF)=\bq_*(\ic(\fQ^0)\otimes\bp^*\CF)
[-\dim\ul\fM^\alpha]$ which is a mixed Hodge module on $\CQ^\alpha$ smooth
along the fine Schubert stratification.
By the Proposition ~12.3b) of {\em loc. cit.} $\bc^0_\CQ(\ic(\oCG_{w,\alpha}))=
\bq_*\ic(\CG\CQ^0_{w,\alpha})$.
By the Decomposition Theorem, $\bq_*\ic(\CG\CQ^0_{w,\alpha})$ is a direct
sum of irreducible Hodge modules on $\CQ^\alpha$.
By ~10.4.2 of {\em loc. cit.} $\bq(\CG\CQ^0_{w,\alpha})=\CQ^\alpha_w$,
and hence $\bq_*\ic(\CG\CQ^0_{w,\alpha})$ contains a direct summand
$\ic(\CQ^\alpha_w)$. Thus it suffices to prove that
$\bc^0_\CQ(\CF)$ is a pointwise pure Tate Hodge module.

More precisely, we are interested in the stalk of $\bc^0_\CQ(\CF)$ at the
point $\phi\in\oQ^\beta_y$. As in the proof of ~12.9.1 of {\em loc. cit.},
we have $\bq^{-1}(\phi)=\oCG_{w,\alpha}\cap\ol{T}{}_{\beta-\alpha}^y=
\sqcup_{\gamma\geq0}(\oCG_{w,\alpha}\cap T_{\beta-\alpha+\gamma}^y)$
where $T_{\beta-\alpha+\gamma}^y$ is $\dot{y}(T_{\beta-\alpha+\gamma})$:
the action of $\dot{y}\in\bG$ on the semiinfinite orbit in $\CG$.
According to the Lemma ~13.1 of {\em loc. cit.}, we have
$\ic(\CG\CQ^0_{w,\alpha})|_{\oCG_{w,\alpha}\cap T_{\beta-\alpha+\gamma}^y}=
\ic(\oCG_{w,\alpha})|_{\oCG_{w,\alpha}\cap T_{\beta-\alpha+\gamma}^y}\otimes
\ic^\gamma_\Gamma$. By ~\ref{zhopa},
$H^\bullet_c(\oCG_{w,\alpha}\cap T_{\beta-\alpha+\gamma}^y,
\ic(\CG\CQ^0_{w,\alpha}))$ is a pure direct sum of Tate Hodge structures,
and by ~\ref{purity}, $\ic^\gamma_\Gamma$ is also a pure direct sum of
Tate Hodge structures; thus
$H^\bullet_c(\oCG_{w,\alpha}\cap T_{\beta-\alpha+\gamma}^y,
\ic(\CG\CQ^0_{w,\alpha})\otimes\ic^\gamma_\Gamma)$
is also a pure direct sum of Tate Hodge structures.

Hence the Cousin spectral sequence associated to
the partition of $\bq^{-1}(\phi)=\oCG_{w,\alpha}\cap\ol{T}{}_{\beta-\alpha}^y$
into the locally closed subsets $\oCG_{w,\alpha}\cap T_{\beta-\alpha+\gamma}^y$
collapses at the second term for the weight reasons.
We deduce that $\bc^0_\CQ(\CF)_{(\phi)}=
H^\bullet(\bq^{-1}(\phi),\ic(\CG\CQ^0_{w,\alpha}))\cong
\oplus_{\gamma\geq0}
H^\bullet_c(\oCG_{w,\alpha}\cap T_{\beta-\alpha+\gamma}^y,
\ic(\CG\CQ^0_{w,\alpha})\otimes\ic^\gamma_\Gamma)$
is also a pure direct sum of Tate Hodge structures. This completes the
proof of the Theorem. $\Box$

\subsection{}
\label{this}
The subsections ~\ref{this}--\ref{long} are quite parallel to \S\S10--12
of ~\cite{fm}.

\subsubsection{}
Let $\sG$ be the usual affine flag manifold $\bG((z))/\bI$. It is the
ind-scheme representing the functor of isomorphism classes of triples
$(\CT,\tau,\ft)$ where $\CT$ is a $\bG$-torsor on $\BP^1$, and $\tau$ is its
section (trivialization) defined off 0, while $\ft$ is its $\bB$-reduction at 0.
It is equipped with a natural action of proalgebraic Iwahori group $\bI$,
and the orbits of this action are naturally numbered by the affine Weyl group
$\CW:\ \sG=\sqcup_{x\in\CW}\sG_x$. An orbit $\sG_y$ lies in the closure
$\osG_x$ iff $y\leq x$ in the usual Bruhat order on $\CW$. The orbit closure
$\osG_x$ has a natural structure of projective variety. The natural projection
$\pr:\ \sG\lra\CG$ is a fibration with the fiber $\bX$. All these facts are
very well known, see e.g. ~\cite{kt}.

\subsubsection{}
Let $\sM$ be the affine flag scheme representing the functor of isomorphism
classes of $\bG$-torsors on $\BP^1$ equipped with trivialization in the
formal neighbourhood of $\infty$ and with $\bB$-reduction at 0
(see ~\cite{ka} and ~\cite{kt1}). It is equipped with a natural action of
the opposite Iwahori group $\bI_-\subset\bG[[z^{-1}]]$, and the orbits of
this action are naturally numbered by the affine Weyl group
$\CW:\ \sM=\sqcup_{x\in\CW}\sM_x$. An orbit $\sM_y$ lies in the closure
$\usM_x$ iff $y\geq x$ in the usual Bruhat order on $\CW$.
For any $x\in\CW$ the union of orbits $\sM^x:=\sqcup_{\CW\ni y\leq x}\sM_y$
forms an open subscheme of $\sM$. This subscheme is a projective limit of
schemes of finite type, all the maps in projective system being fibrations
with affine fibers. Moreover, $\sM^x$ is equipped with a free action of a
prounipotent group $\bG^x$ (a congruence subgroup in $\bG[[z^{-1}]]$) such
that the quotient $\usM^x$ is a smooth scheme of finite type.

The natural projection $\pr:\ \sM\lra\fM$ (see ~\cite{fm} ~10.6) forgetting
a $\bB$-reduction at 0 is a fibration with the fiber $\bX$.

Restricting a trivialization of a $\bG$-torsor from $\BP^1-0$ to the formal
neighbourhood of $\infty$ we obtain the closed embedding $\bi:\ \sG\lra\sM$.
The intersection of $\sM_x$ and $\sG_y$ is nonempty iff $y\leq x$, and
then it is transversal. Thus, $\osG_x\subset\sM^x$. According to ~\cite{kt},
the composition $\osG_x\hookrightarrow\sM^x\lra\usM^x$ is a closed embedding.

\subsubsection{}
Following ~10.6.3 of ~\cite{fm} we define
for {\em arbitrary} $\alpha\in Y$ the scheme $\osQ^\alpha$
(resp. $\sQ^\alpha$) representing the functor of isomorphism classes of triples
$(\CT,\ft,(\fL_\lambda)_{\lambda\in X^+})$
where $\CT$ is a $\bG$-torsor trivialized in the formal neighbourhood of
$\infty\in\BP^1$, and $\ft$ is its $\bB$-reduction at 0, while
$\fL_\lambda\subset\CV_\lambda^\CT,\
\lambda\in X^+$, is a collection of line subbundles (resp. invertible
subsheaves) of degree $\langle-\alpha,\lambda\rangle$ satisfying the Pl\"ucker
conditions (cf. {\em loc. cit.}). The evident projection $\osQ^\alpha\lra\sM$
(resp. $\sQ^\alpha\lra\sM$) will be denoted by $\obp$ (resp. $\bp$).
The open embedding $\osQ^\alpha\hra\sQ^\alpha$ will be denoted by $\bj$.
Clearly, $\bp$ is projective, and $\obp=\bp\circ\bj$.

The free action of prounipotent group $\bG^x$ on $\sM^x$ lifts to
the free action of $\bG^x$ on the open subscheme
$\bp^{-1}(\sM^x)\subset\sQ^\alpha$.
The quotient is a scheme of finite type $\usQ^{\alpha,x}$ equipped with
the projective morphism $\bp$ to $\usM^x$. There exists a
$\bI_-$-invariant stratification $\fS$ of $\sQ^\alpha$ such that
$\bp$ is stratified with respect to $\fS$ and the stratification
$\sM=\sqcup_{x\in\CW}\sM_x$. One can define perverse sheaves and mixed Hodge
modules on $\sQ^\alpha$ smooth along $\fS$ following the lines of ~\cite{kt1}.
In particular, we have the irreducible Hodge module
$\ic(\sQ^\alpha)$.

\subsubsection{}
\label{following}
Following ~8.4.1 of ~\cite{fm}
we introduce the {\em fine stratification} of $\sQ^\alpha$ according to
defects of invertible subsheaves:
$$\sQ^\alpha=\bigsqcup^{\alpha\geq\beta\geq\gamma}_{\Gamma\in\fP(\beta-\gamma)}
\osQ^\gamma\times(C-0)^{\beta-\gamma}_\Gamma\times0^{\alpha-\beta}$$
and the {\em fine Schubert stratification} of $\sQ^\alpha$:
$$\sQ^\alpha=
\bigsqcup^{\alpha\geq\beta\geq\gamma}_{w\in\CW_f,\ \Gamma\in\fP(\beta-\gamma)}
\osQ^\gamma_w\times(C-0)^{\beta-\gamma}_\Gamma\times0^{\alpha-\beta}$$
where $\osQ^\gamma_w\subset\osQ^\gamma$ is defined as follows. A collection
of line subbundles $\fL_\lambda\subset\CV^\CT_\lambda$ defines a collection
of lines $L_\lambda\subset(\CV^\CT_\lambda)_{(0)}$ in the fiber over $0\in C$.
This collection of lines defines a $\bB$-reduction of $\CT$ at $0\in C$,
and if the relative position of this reduction and $\ft$ is $w\in\CW_f$, we
say that the triple $(\CT,\ft,(\fL_\lambda)_{\lambda\in X^+})$ lies in
$\osQ^\gamma_w$.

The closure of $\osQ^\gamma_w$ is denoted by $\sQ^\gamma_w\subset
\sQ^\gamma$. Also, for an alcove $A=(w,\gamma)$ we will often write
$\osQ_A$ for $\osQ^\gamma_w$, and $\sQ_A$ for $\sQ^\gamma_w$.
An open subset of $\sQ^\gamma_w$ formed by the triples
$(\CT,\ft,(\fL_\lambda)_{\lambda\in X^+})$ such that the collection
$(\fL_\lambda)_{\lambda\in X^+}$ has no defect at $0\in C$ will be denoted
by $\dsQ^\gamma_w$.

The {\em refined stratification} of $\sQ^\alpha$ is a refinement of the
fine Schubert stratification: it additionally subdivides $\osQ^\gamma_w$
into the strata $\osQ^{\gamma,x}_w:=\bp^{-1}(\sM_x)$.
The map $\bp:\ \sQ^\alpha\lra\sM$ is clearly stratified with respect
to the refined stratification of $\sQ^\alpha$ and the stratification
$\sM=\sqcup_{x\in\CW}\sM_x$.

\subsubsection{}
\label{parallel}
For $\beta\geq\alpha$ we have a closed embedding $\sQ^\alpha\hookrightarrow
\sQ^\beta$, sending a triple $(\CT,\tau,(\fL_\lambda)_{\lambda\in X^+})$
to $(\CT,\tau,
(fL_\lambda(\langle\beta-\alpha,\lambda\rangle\cdot0))_{\lambda\in X^+})$.
These embeddings form an inductive system, and we will denote its union
by $\sQ$. A {\em perverse sheaf} or a {\em mixed Hodge module on $\sQ$
supported on $\sQ^\alpha$} is just the same as a perverse sheaf or a mixed
Hodge module on $\sQ^\alpha$. Let $\aleph(\sQ)$ be the additive category formed
by the (finite) direct sums of irreducible Hodge modules $\ic(\sQ^\gamma_w)$
and their Tate twists on $\sQ$.
Let $\mho(\sQ)$ be the additive category formed by the direct sums
of Goresky-MacPherson sheaves $\IC(\sQ^\gamma_w)$ and their shifts in the
derived category. The subcategories of sheaves supported on $\sQ^\alpha\subset
\sQ$ will be denoted by $\aleph(\sQ^\alpha)$ and $\mho(\sQ^\alpha)$.

\subsubsection{Remark}
\label{brav} 
The category of perverse sheaves on $\sQ$ of finite length,
with all the irreducible constituents of the form $\IC(\sQ^\gamma_w)$, 
is equivalent
to the category $\PS$ defined in ~\cite{fm}. This is the good working 
definition of $\PS$ promised in ~\cite{fm} instead of the ugly provisional
definition given there.

\subsubsection{}
\label{similar}
Similarly, let
$\aleph(\sG)$ be the additive category formed
by the (finite) direct sums of irreducible Hodge modules $\ic(\osG_x)$
and their Tate twists on $\sG$.
Let $\mho(\sG)$ be the additive category formed by the direct sums
of Goresky-MacPherson sheaves $\IC(\osG_x)$ and their shifts in the
derived category. The subcategories of sheaves supported on $\osG_x\subset\sG$
will be denoted by $\aleph(\osG_x)$ and $\mho(\osG_x)$.

For $\alpha\in\BN[I]$
let $\aleph(\CQ^\alpha)$ be the additive category formed
by the (finite) direct sums of irreducible Hodge modules $\ic(\CQ^\gamma_w),\
\gamma\leq\alpha$, and their Tate twists on $\CQ^\alpha$.
Let $\mho(\CQ^\alpha)$ be the additive category formed by the direct sums
of Goresky-MacPherson sheaves $\IC(\CQ^\gamma_w),\ \gamma\leq\alpha$,
and their shifts in the derived category.

For $\alpha\in\BN[I]$ sufficiently dominant
let $\aleph(\ddZ^\alpha)$ be the additive category formed
by the (finite) direct sums of irreducible Hodge modules $\ic(\CZ^\gamma_w),\
0\ll\gamma\leq\alpha$, and their Tate twists on $\ddZ^\alpha$ (see ~\ref{mrak}).
Let $\mho(\ddZ^\alpha)$ be the additive category formed by the direct sums
of Goresky-MacPherson sheaves $\IC(\CZ^\gamma_w),\ 0\ll\gamma\leq\alpha$,
and their shifts in the derived category.

Finally, let $\aleph(\PS)$ be the additive category formed by the collections
$\fH^\alpha_\chi\in\aleph(\ddZ^\alpha_\chi)$ together with factorization
isomorphisms as in ~\cite{fm} ~9.3. Let $\mho(\PS)$ be the additive category
formed by the direct sums of irreducible snops $\CL(A)\in\PS$
(see ~\ref{t'ma}) and their shifts in the derived category.

\subsection{}
\label{that}
For $\fF,\fG\in\mho(\sG)$ their convolution $\fF\star\fG\in\mho(\sG)$
was studied e.g. in ~\cite{kt}. We will define the convolution
$\mho(\sG)\times\mho(\sQ)\lra\mho(\sQ)$. Its construction occupies the
subsections ~\ref{that}--\ref{long}.

\subsubsection{}
\label{degree 0}
Let $x\in\CW,\ w\in\CW_f,\ \alpha\in Y$.
We define the {\em convolution diagram} $\sG\sQ^\alpha_{w,x}$ as the cartesian
product of $\osG_x$ and $\sQ^\alpha_w$ over $\sM$. Thus we have the following
cartesian diagram (cf. ~\cite{fm} ~12.2):
$$
\begin{CD}
\sG\sQ^\alpha_{w,x}       @>\bi>>    \sQ^\alpha_w        \\
@V{\bp}VV               @V{\bp}VV               \\
\osG_x    @>{\bi}>>     \sM
\end{CD}
$$
The same argument as in {\em loc. cit.} ~12.2 proves that
$\ic(\sQ^\alpha_w)\otimes\bp^*\bi_*\ic(\osG_x)[-\dim\usM^x]\cong
\ic(\sG\sQ^\alpha_{w,x})$, and moreover, for any $\fF\in\aleph(\osG_x)$
the complex $\ic(\sQ^\alpha_w)\otimes\bp^*\bi_*\fF[-\dim\usM^x]$
is in fact a semisimple Hodge module (living in cohomological degree 0).
Furthermore, the Verdier duality $\CalD$ takes
$\ic(\sQ^\alpha_w)\otimes\bp^*\bi_*\fF[-\dim\usM^x]$ to
$\ic(\sQ^\alpha_w)\otimes\bp^*\bi_*\CalD\fF[-\dim\usM^x]$
(see {\em loc. cit.} ~12.1.b).

\subsubsection{}
\label{convolution}
It is well known that $\CW_f\backslash\CW/\CW_f=Y^+$. Let us denote the
double coset of $x\in\CW$ by $\eta\in Y^+$. Then the image
$\pr(\osG_x)$ in the affine Grassmannian $\CG$ lies in
$\oCG_\eta$. Suppose
$\eta+\alpha\in\BN[I]$. Comparing the above definition of
$\sG\sQ^\alpha_{w,x}$ with {\em loc. cit.} ~12.2 we obtain the natural
projection $\pr:\ \sG\sQ^\alpha_{w,x}\lra\CG\CQ^\alpha_\eta$.
Recall the map $\bq:\ \CG\CQ^\alpha_\eta\lra\CQ^{\eta+\alpha}$ defined in
{\em loc. cit.} ~11.2.

For $\fF\in\mho(\oCG_x)$ we define the {\em convolution}
$$\fF\star\IC(\sQ^\alpha_w):=\bq_*\pr_*
(\ic(\sQ^\alpha_w)\otimes\bp^*\bi_*\fF[-\dim\usM^x])$$ By the decomposition
theorem $\fF\star\IC(\sQ^\alpha_w)\in\mho(\CQ^{\eta+\alpha})$.
By additivity, the convolution extends to the functor
$$\star:\ \mho(\osG_x)\times\mho(\sQ^\alpha)\lra\mho(\CQ^{\eta+\alpha})$$
By the last sentence of ~\ref{degree 0} this functor commutes with Verdier
duality:
$$\CalD\fF\star\CalD\fH\iso\CalD(\fF\star\fH)$$

\subsection{}
\label{label}
Let $\alpha\in Y,\ w\in\CW_f$. Suppose a fine Schubert stratum
$\osQ^\gamma_w\times(C-0)^{\beta-\gamma}_\Gamma\times0^{\alpha-\beta}$
lies in $\sQ^\alpha_y$ (see ~\ref{following}). The stalk of $\ic(\sQ^\alpha_y)$
at (any point $\phi$ in) the stratum
$\osQ^\gamma_w\times(C-0)^{\beta-\gamma}_\Gamma\times0^{\alpha-\beta}$
will be denoted by $\ic(\sQ^\alpha_y)_{\beta,\gamma,\Gamma,w}$.
Let $\xi\in Y^+$ be dominant enough so that $\xi+\gamma\in\BN[I]$
(hence $\xi+\alpha\in\BN[I]$). Then we may consider the stalk
$\ic(\CQ^{\xi+\alpha}_y)_{\xi+\beta,\xi+\gamma,\Gamma,w}$
of $\ic(\CQ^{\xi+\alpha}_y)$ at (any point in) the stratum
$\oQ^{\xi+\gamma}_w\times(C-0)^{\beta-\gamma}_\Gamma\times0^{\alpha-\beta}
\subset\CQ^{\xi+\alpha}$.

{\bf Theorem.}
$\ic(\sQ^\alpha_y)_{\beta,\gamma,\Gamma,w}$ is isomorphic, up to a shift, to
$\ic(\CQ^{\xi+\alpha}_y)_{\xi+\beta,\xi+\gamma,\Gamma,w}$.

{\em Proof.} The argument is parallel to the one in ~\cite{fm} ~12.6, ~12.7.
Let us consider another copy $\CG^1$ of affine Grassmannian: the ind-scheme
representing the functor of isomorphism classes of pairs $(\CT,\tau')$ where
$\CT$ is a $\bG$-torsor on $C$, and $\tau'$ is its section (trivialization)
{\em defined off $1\in C$}. Restricting a trivialization to the formal
neighbourhood of $\infty\in C$ we obtain the closed embedding $\bi^1:\
\CG^1\hookrightarrow\fM$. Furthermore, since $\tau'$ identifies the fiber of
$\CT$ over $0\in C$ with $\bG$, the set of $\bB$-reductions of $\CT$ at 0
is canonically isomorphic to $\bX$. Choosing the reduction $\ft=e=\bB\in\bX$
we obtain the same named closed embedding $\bi^1:\ \CG^1\hookrightarrow\sM$.

We consider the cartesian product $\CG^1\sQ^\alpha_{y,\eta}$ of
$\oCG_\eta$ and $\sQ^\alpha_y$ over $\sM$. Thus we have the following
cartesian diagram (cf. ~\ref{degree 0}):
$$
\begin{CD}
\CG^1\sQ^\alpha_{y,\eta}       @>\bi^1>>    \sQ^\alpha_y        \\
@V{\bp}VV               @V{\bp}VV               \\
\oCG^1_\eta    @>{\bi^1}>>     \sM
\end{CD}
$$
The same argument as in ~\ref{degree 0} proves that up to a shift
$\ic(\sQ^\alpha_y)\otimes\bp^*\bi^1_*\ic(\oCG^1_\eta)\cong
\ic(\CG^1\sQ^\alpha_{y,\eta})$.

Recall the convolution diagram $\CG\CQ^\alpha_\eta$
(see ~\cite{fm} ~11.1, ~12.2).
The point $0\in C$ played a special role in the definition of this diagram.
Let us replace 0 by 1 in the definition of $\CG\CQ^\alpha_\eta$.
Let us call the result $\CG^1\CQ^\alpha_\eta$. One defines the map
$\bq^1:\ \CG^1\CQ^\alpha_\eta\lra\CQ^{\eta+\alpha}$ as in {\em loc. cit.} ~11.2.
Then the Proposition ~12.6 of {\em loc. cit.} states that the restriction
of $\bq^1$ to $\dQ^{1,\eta+\alpha}$ is an isomorphism, where
$\dQ^{1,\eta+\alpha}\subset\CQ^{\eta+\alpha}$ is an open subset formed by
all the quasimaps defined at $1\in C$ (i.e. without defect at 1).
A moment of reflection shows that there is no difference
between the definitions of
$\CG^1\sQ^\alpha_{w_0,\eta}$ and $\CG^1\CQ^\alpha_\eta$. Identifying them
and imbedding $\CG^1\sQ^\alpha_{y,\eta}$ into $\CG^1\sQ^\alpha_{w_0,\eta}$
we obtain the map $\bq^1:\ \CG^1\sQ^\alpha_{y,\eta}\lra\CQ^{\eta+\alpha}$
which is a closed embedding over $\dQ^{1,\eta+\alpha}$. Clearly, the image
of this embedding coincides with $\CQ^{1,\eta+\alpha}_y:=\CQ^{\eta+\alpha}_y
\cap\dQ^{1,\eta+\alpha}$. Evidently, $\CQ^{1,\eta+\alpha}_y$ is open in
$\CQ^{\eta+\alpha}_y$.

Now let us consider the stalk of $\ic(\sQ^\alpha_y)$ at a point
$\phi=(\CT,\ft,(\fL_\lambda)_{\lambda\in X^+})$ in the stratum
$\osQ^\gamma_w\times(C-0)^{\beta-\gamma}_\Gamma\times0^{\alpha-\beta}$.
Suppose that the isomorphism class of $\bG$-torsor $\CT$ equals $\eta\in Y^+$,
i.e. $\CT\in\fM_\eta$. The stalk in question does not depend on a choice
of $\CT\in\fM_\eta$ and the defect $D\in(C-0)^{\beta-\gamma}_\Gamma$.
In particular, we may (and will) suppose that $\CT\in\bi^1(\CG^1_\eta)$,
and $D\in(C-0-1)^{\beta-\gamma}$. We have seen that up to a shift
$\ic(\sQ^\alpha_y)_\phi\otimes\bp^*\bi^1_*\ic(\oCG_\eta)_\CT\cong
\ic(\CG^1\sQ^\alpha_{y,\eta})_\phi$. But the stalk of $\ic(\oCG^1_\eta)$ at
any point $\CT\in\CG^1_\eta$ is isomorphic, up to a shift, to the trivial
Tate Hodge structure $\BQ(0)$. We deduce that up to a shift
$\ic(\sQ^\alpha_y)_\phi\cong\ic(\CG^1\sQ^\alpha_{y,\eta})_\phi$.

On the other hand, the latter stalk is isomorphic to the stalk of
$\ic(\CQ^{1,\eta+\alpha}_y)$ at the point $\bq^1(\phi)$. This point lies in
$\oQ^{\eta+\gamma}_w\times(C-0)^{\beta-\gamma}_\Gamma\times0^{\alpha-\beta}$.

Thus we have proved that up to a shift
$\ic(\sQ^\alpha_y)_{\beta,\gamma,\Gamma,w}$ is isomorphic to
$\ic(\CQ^{\eta+\alpha}_y)_{\eta+\beta,\eta+\gamma,\Gamma,w}$.

The standard use of factorization shows that the latter stalk is independent
of the shift $\eta\mapsto\xi$. This completes the proof of the Theorem. $\Box$

\subsubsection{Corollary}
\label{pure Tate}
The Hodge module $\ic(\sQ^\alpha_y)$ is pointwise pure Tate. $\Box$

\subsubsection{Corollary}
Let $\beta,\gamma\in Y,\ w,y\in\CW_f$. Then $\sQ^\beta_y\supset\sQ^\gamma_w$
iff $(y,\beta)=:B\geq A:=(w,\gamma)$.

{\em Proof.} $\sQ^\beta_y\supset\sQ^\gamma_w$ iff the stalk of
$\ic(\sQ^\beta_y)$ at the generic point of $\sQ^\gamma_w$ does not vanish.
Also, $\CQ^\beta_y\supset\CQ^\gamma_w$ iff the stalk of
$\ic(\CQ^\beta_y)$ at the generic point of $\CQ^\gamma_w$ does not vanish.
Now recall that the order on alcoves is translation invariant, and apply
~\ref{who}. $\Box$

\subsection{}
\label{short}
Before we finish the construction of convolution
$\star:\ \mho(\sG)\times\mho(\sQ)\lra\mho(\sQ)$ we have to define
the numerous functors between the categories
$\aleph$ and $\mho$ of ~\ref{parallel} and ~\ref{similar}.

First of all, we choose an additive equivalence $\fE^\PS_\sQ:\
\aleph(\sQ)\lra\aleph(\PS)$ sending $\ic(\sQ^\eta_w),\ \eta\in Y,\ w\in\CW_f$,
to the collection $\ic(\CZ^\alpha_{w,\eta})\in\aleph(\ddZ^\alpha_\eta)$
(see ~\ref{mrak} and ~\cite{fm} ~9.3). Let $\fE_\PS^\sQ:\
\aleph(\PS)\lra\aleph(\sQ)$ be an inverse equivalence.

We preserve the name $\fE^\PS_\sQ$ for an additive equivalence
$\mho(\sQ)\lra\mho(\PS)$ sending $\IC(\sQ^\eta_w),\ \eta\in Y,\ w\in\CW_f$,
to $\CL(w,\eta)\in\mho(\PS)$ and commuting with shifts. Let $\fE_\PS^\sQ:\
\mho(\PS)\lra\mho(\sQ)$ be an inverse equivalence.

For $\alpha\in\BN[I]$ let $\aleph_{\gg0}(\CQ^\alpha)\subset\aleph(\CQ^\alpha)$
be the direct summand subcategory formed by the direct sums of irreducible
Hodge modules $\ic(\CQ^\gamma_w),\ 0\ll\gamma\leq\alpha$, and their Tate
twists. Let $\mho_{\gg0}(\CQ^\alpha)\subset\mho(\CQ^\alpha)$
be the direct summand subcategory formed by the direct sums of
Goresky-MacPherson sheaves $\IC(\CQ^\gamma_w),\ 0\ll\gamma\leq\alpha$, and
their shifts in derived category. Let
$\fE_{\CQ^\alpha}^{\gg0}:\ \aleph(\CQ^\alpha)\lra\aleph_{\gg0}(\CQ^\alpha),\
\mho(\CQ^\alpha)\lra\mho_{\gg0}(\CQ^\alpha)$ denote the projections to the
direct summands, and let
$\fE^{\CQ^\alpha}_{\gg0}:\ \aleph_{\gg0}(\CQ^\alpha)\lra\aleph(\CQ^\alpha),\
\mho_{\gg0}(\CQ^\alpha)\lra\mho(\CQ^\alpha)$ denote the embeddings of direct
summands.

Let $\ss_\alpha$ denote the locally closed embedding
$\ddZ^\alpha\hookrightarrow\CQ^\alpha$. One can easily see that for
$\gamma\gg0$ we have
$\ss_\alpha^*\ic(\CQ^\gamma_w)[-\dim\bX]=\ic(\CZ^\gamma_w)$.
Hence the functor $\ss_\alpha^*[-\dim\bX]$ defines an equivalence
$\fE^{\ddZ^\alpha}_{\CQ^\alpha}:\
\aleph_{\gg0}(\CQ^\alpha)\lra\aleph(\ddZ^\alpha),\
\mho_{\gg0}(\CQ^\alpha)\lra\mho(\ddZ^\alpha)$. Let
$\fE_{\ddZ^\alpha}^{\CQ^\alpha}:\
\aleph(\ddZ^\alpha)\lra\aleph_{\gg0}(\CQ^\alpha),\
\mho(\ddZ^\alpha)\lra\mho_{\gg0}(\CQ^\alpha)$ be an inverse equivalence.

Finally, for $\xi\in Y$ let $\xi_*:\ \aleph(\sQ)\lra\aleph(\sQ)$
(resp. $\mho(\sQ)\lra\mho(\sQ)$) be an additive equivalence taking
$\ic(\sQ^\alpha_w)$ to $\ic(\sQ^{\xi+\alpha}_w)$
(resp. $\IC(\sQ^\alpha_w)$ to $\IC(\sQ^{\xi+\alpha}_w)$ and commuting with
shifts).

\subsection{}
\label{long}
Finally, we are in a position to define the convolution
$\star:\ \mho(\sG)\times\mho(\sQ)\lra\mho(\sQ)$. To this end recall the
convolution $\mho(\osG_x)\times\mho(\sQ^\alpha)\lra\mho(\CQ^{\eta+\alpha})$
defined in ~\ref{that}. To stress the dependence on $x\in\CW$ and $\alpha\in Y$
we will call this convolution $\star^\alpha_x$ (recall that $\eta\in Y^+=
\CW_f\backslash\CW/\CW_f$ is just the double coset of $x$).

Let $\fF\in\mho(\osG_x),\ \fH\in\mho(\sQ^\alpha)$.
The collection $\fE^{\ddZ^\beta}_{\CQ^\beta}\circ\fE^{\gg0}_{\CQ^\beta}
(\fF\star^{\beta-\eta}_x((\beta-\alpha-\eta)_*\fH))
\in\mho(\ddZ^\beta_{\alpha+\eta}),\ \beta\in\BN[I]$, defines an object in
$\mho(\PS)$. Applying $\fE_\PS^\sQ$ to this object we obtain the desired
convolution $\fF\star^\alpha_x\fH\in\mho(\sQ)$.
Note that if $y\leq x,\ \gamma\leq\alpha$, and $\fF$ is actually supported
on $\osG_y$, while $\fH$ is actually supported on $\sQ^\gamma$, then we
have the canonical isomorphism $\fF\star^\gamma_y\fH\equiv\fF\star^\alpha_x\fH$.
Hence the convolutions $\star^\alpha_x:\ \mho(\osG_x)\times\mho(\sQ^\alpha)\lra
\mho(\sQ)$ extends to the convolution
$\star:\ \mho(\sG)\times\mho(\sQ)\lra\mho(\sQ)$.

As we have seen in ~\ref{that}, this convolution commutes with Verdier duality.

Moreover, for $\fF,\fG\in\mho(\sG),\ \fH\in\mho(\sQ)$ one can easily construct
a canonical isomorphism $\fF\star(\fG\star\fH)\cong(\fF\star\fG)\star\fH$.

Finally, composing this convolution with the mutually inverse equivalences
$\fE^\PS_\sQ:\ \mho(\sQ)\rightleftharpoons\mho(\PS)\ :\fE_\PS^\sQ$ we
obtain the same named convolution functor
$\star:\ \mho(\sG)\times\mho(\PS)\lra\mho(\PS)$. It is also associative
with respect to convolution in $\mho(\sG)$, and commutes with Verdier duality.

\subsection{}
\label{main}
Let $K(\sQ)$ (resp. $K(\PS)$) denote the $K$-group of the category
$\mho(\sQ)$ (resp. $\mho(\PS)$): it is generated by the isomorphism classes
of objects in $\mho(\sQ)$ with relations $[K_1\oplus K_2]=[K_1]+[K_2]$.
Both $K(\sQ)$ and $K(\PS)$ are equipped with
the structure of $\BZ[v,v^{-1}]$-modules: $-v$ acts as the shift $[-1]$.

$K(\PS)$ is a free $\bzw$-module with a basis $[\CL(A)],\ A\in\fA$
(see ~\ref{t'ma}). Similarly, $K(\sQ)$ is a free $\bzw$-module with a basis
$[\CL'(A)],\ A\in\fA$ where for $A=(w,\eta)$ we denote by $\CL'(A)$ the
Goresky-MacPherson sheaf $\IC(\sQ^\eta_w)=\IC(\sQ_A)$ (see ~\ref{following}).

We denote by $\hK(\PS)\supset K(\PS)$ (resp. $\hK(\sQ)\supset K(\sQ)$)
the completed module formed by the possibly infinite sums
$\sum_{A\in\fA}a_A[\CL(A)]$ (resp. $\sum_{A\in\fA}a_A[\CL'(A)]$)
$a_A\in\bzw$ such that $\exists B\in\fA$ with the property
$(a_A\ne0\ \Rightarrow\ A\leq B)$.

We introduce another (topological) basis $\{[\CM'(A)],\ A\in\fA\}$
of $\hK(\sQ)$. It is characterized by the following property. The stalk
at a fine Schubert stratum $\osQ_B\subset\sQ_B$ defines the functional
$f_B:\ \hK(\sQ)\lra\bzw$. We require $f_A([\CM'(A)])=f_A([\CL'(A)])$ and
$f_B([\CM'(A)])=0$ for $B\ne A$.
In a similar way one defines a (topological) basis $\{[\CM(A)],\ A\in\fA\}$
of $\hK(\PS)$.

The main goal of this section is a computation of the transition matrix
between the bases $[\CL(A)]$ and $[\CM(A)]$. We will make a free use of
the notations and results of the excellent exposition ~\cite{s}. In particular,
the polynomials $\ol{q}_{B,A}$ (generic Kazhdan-Lusztig polynomials)
are introduced in {\em loc. cit.}, Theorems ~6.1, ~6.4.
We will prove the following

{\bf Theorem.} a) $[\CL(A)]=\sum_{B\leq A}(-1)^{d(A,B)}\ol{q}_{B,A}[\CM(B)]$;

b) $[\CL'(A)]=\sum_{B\leq A}(-1)^{d(A,B)}\ol{q}_{B,A}[\CM'(B)]$.

The proof of the Theorem occupies the subsections ~\ref{start}--\ref{finish}.

\subsection{}
\label{start}
For the time being let us define $\fq_{B,A}\in\bzw$ (resp. $\fq'_{B,A}\in\bzw$)
by the requirement
$[\CL(A)]=\sum_{B\leq A}\fq_{B,A}[\CM(B)]$ (resp.
$[\CL'(A)]=\sum_{B\leq A}\fq'_{B,A}[\CM'(B)]$). Thus we have to prove
$\fq_{B,A}=\fq'_{B,A}=(-1)^{d(A,B)}\ol{q}_{B,A}$.

By the Corollary ~\ref{pure Tate} the Hodge module $\ic(\sQ_A)$ is pointwise
pure Tate. This allows us to reconstruct the stalk $\ic(\sQ_A)_B$
of $\ic(\sQ_A)$ at
$\osQ_B\subset\sQ_A$ from the polynomials $\fq'_{B,A}$. Namely, suppose
that the stalk of $\ic(\sQ_A)$ at $\osQ_A$ is the Hodge structure $\BQ(0)$
living in cohomological degree $\delta_A$ (this is the definition of
$\delta_A$). Then for an odd integer $k$ we have
$H^{\delta_A+k}\ic(\sQ_A)_B=0$. For arbitrary integer $k$ the cohomology
$H^{\delta_A+2k}\ic(\sQ_A)_B$ is a sum of a few copies of $\BQ(k)$. Finally,
$\dim H^{\delta_A+k}\ic(\sQ_A)_B$ equals the coefficient of
$(-v)^{-d(B,A)+k}$ in $\fq'_{B,A}$ (for arbitrary $k\in\BZ$).
The same discussion applies to the stalks $\ic(\CZ_A)_B$ of $\ic(\CZ_A)$.

But the stalks $\ic(\CZ_A)_B$ and $\ic(\sQ_A)_B$ are isomorphic. In effect,
if $A=(w,\alpha),\ B=(y,\beta)$, and $\xi\in Y^+$ is big enough, then by
the Theorem ~\ref{label} $\ic(\sQ_A)_B$ is isomorphic to the stalk of
$\ic(\CQ^{\xi+\alpha}_w)$ at $\oQ^{\xi+\beta}_y$ (up to a shift; but we may
disregard shifts: as we have just seen the cohomological degrees contain
exactly the same information as the Hodge structures by the pointwise Tate
purity). And the latter stalk is isomorphic to $\ic(\CZ_A)_B$ since
$\ic(\CZ^{\xi+\alpha}_w)=\ss_{\xi+\alpha}^*\ic(\CQ^{\xi+\alpha}_w)[-\dim\bX]$
where $\ss_{\xi+\alpha}$ stands for the locally closed embedding of
$\ddZ^{\xi+\alpha}$ into $\CQ^{\xi+\alpha}$.

This implies that the functor $\fE_\PS^\sQ:\ \mho(\PS)\lra\mho(\sQ)$
inducing an isomorphism of $\hK(\PS)$ and $\hK(\sQ)$ sending
$[\CL(A)]$ to $[\CL'(A)]$ {\em also sends $[\CM(A)]$ to $[\CM'(A)]$}.
In other words, $\fE_\PS^\sQ$ preserves stalks.
In particular, $\fq_{B,A}=\fq'_{B,A}$ for all $B\leq A$.
From now on we will identify $\hK(\PS)$ and $\hK(\sQ)$, we will write
$[\CL(A)]$ for $[\CL'(A)]$, and $[\CM(A)]$ for $[\CM'(A)]$.

\subsubsection{Remark} As we have seen, the stalk of $\ic(\CQ^\alpha_w)$
at a fine Schubert stratum $\oQ^\beta_y$ can be reconstructed in terms of
polynomials $\fq_{B,A}$ for $A=(w,\alpha),B=(y,\beta)$. Since we know the
simple stalks $\ic^0_{\{\{\gamma\}\}}$, the factorization property allows
us to reconstruct all the other stalks of $\ic(\CQ^\alpha_w)$.

\subsection{}
\label{kato}
We know the simple stalks $\ic^0_{\{\{\gamma\}\}}$ by the Theorems
~\ref{simple}, ~\ref{Hodge}. They coincide with the stalks of
$\ic(\sQ^\alpha)=\ic(\sQ^\alpha_{w_0})$. Thus we deduce the following formula:
$$[\CL(w_0,\alpha)]=\sum_{w\in\CW_f}^{\beta\in\BN[I]}
(-v)^{l(w)-l(w_0)}\CK^\beta(v^2)[\CM(w,\alpha-\beta)]$$
Comparing with Kato's Theorem (see ~\cite{kat} or ~\cite{s} ~6.3) we see
that for any $B\in\fA$ and $A=(w_0,\alpha)$ we have
$\fq_{B,A}=(-1)^{d(A,B)}\ol{q}_{B,A}$.

\subsection{}
\label{last}
The $K$-group $K(\sG):=K(\mho(\sG))$ is also a $\bzw$-module in a natural way
($-v$ acts as a shift $[-1]$). The convolution $\star:\ \mho(\sG)\times
\mho(\sG)\lra\mho(\sG)$ makes $K(\sG)$ into a $\bzw$-algebra. This is the
affine Hecke algebra $\CH$ as described e.g. in ~\cite{s} ~\S2. The basis
$\{H_x,\ x\in\CW\}$ corresponds to perverse shriek extensions of constant
sheaves on orbits $\sG_x$. The involution $H\mapsto\ol{H}$ of {\em loc. cit.}
is induced by the Verdier duality $\CalD:\ \mho(\sG)\lra\mho(\sG)$.

Let $s_i\in\CW,\ i\in I\sqcup0$ be a simple reflection (see ~\ref{fuck}).
Then $\osG_{s_i}\subset\sG$ is a projective line, to be denoted by $\osG_i$.
The class $[\IC(\osG_i)]\in\CH$ equals $H_{s_i}-v^{-1}$
(we identify $H_e$ with $1\in\CH$). Following W.Soergel, we will denote
$[\IC(\osG_i)]$ by $\tC_i=\tC_{s_i}$.

The convolution $\star:\ \mho(\sG)\times\mho(\sQ)\lra\mho(\sQ)$ gives rise to
the structure of $\CH$-module on $\hK(\sQ)$.

{\bf Proposition.} Let $i\in I\sqcup0,\ A\in\fA$. If $s_iA>A$ then
$\tC_i[\CM(A)]=[\CM(s_iA)]-v^{-1}[\CM(A)]$, and if $s_iA<A$ then
$\tC_i[\CM(A)]=[\CM(s_iA)]-v[\CM(A)]$.

{\em Proof.} As everything is invariant with respect to translations, we may
(and will) assume that $A=(w,\alpha)$ for $\alpha\in\BN[I]$.
Consider the convolution diagram $\sG\sQ^\alpha_{w,s_i}\stackrel{\bi}
{\hookrightarrow}\sQ^\alpha_w$. We will denote by $\sG\dsQ^\alpha_{w,s_i}$
the intersection of $\sG\sQ^\alpha_{w,s_i}$ with $\dsQ^\alpha_w\subset
\sQ^\alpha_w$ (see ~\ref{following}). To compute $\tC_i[\CM(A)]$ we have
to compute the stalks of $(\bq\circ\pr)_!\IC(\sG\dsQ^\alpha_{w,s_i})$.
Recall that $\bq\circ\pr$ maps $\sG\sQ^\alpha_{w,s_i}$ to
$\CQ^{\eta+\alpha}$ where $\eta\in Y^+=\CW_f\backslash\CW/\CW_f$ is the
double coset of $s_i$. Note that the double coset of $s_0$ equals $\beta_0$
(see ~\ref{fuck}) --- the coroot dual to the highest root. And the double
coset of $s_i,\ i\in I$, equals $0\in Y^+$. So in any case we have the map
$\bq\circ\pr:\ \sG\dsQ^\alpha_{w,s_i}\lra\CQ^{\beta_0+\alpha}$.

According to ~\ref{label} the sheaf $\IC(\sG\dsQ^\alpha_{w,s_i})$ is constant
along the fibers of $\bq\circ\pr$. So the Proposition follows from the
following Claim.

\subsubsection{Claim}
\label{claim}
Let $\phi\in\dQ^\beta_y\subset\CQ^{\beta_0+\alpha}$,
and let $B=(y,\beta)\in\fA$. Then

a) if $B\ne A,s_iA$, then $(\bq\circ\pr)^{-1}(\phi)=\emptyset$;

b) if $B$ equals $A$ or $s_iA$, and $s_iA>A$, then $(\bq\circ\pr)^{-1}(\phi)$
consists of exactly one point;

c) if $B$ equals $A$ or $s_iA$, and $s_iA<A$, then $(\bq\circ\pr)^{-1}(\phi)$
is isomorphic to the affine line $\BA^1$.

\subsubsection{Remark} If $i\in I$ then the double coset of $s_i$ equals 0, and
the above map $\bq$ is just an isomorphism. In this case the convolution boils
down to the usual convolution on $\bX$. The order relation between $s_iA$
and $A$ also reduces to the usual Bruhat order on $\CW_f$
(see ~\ref{s0} and ~\ref{leq}). So ~\ref{claim} follows in this case from
the standard facts about the geometry of $\bX$.

\subsubsection{} In the general case note first that the (isomorphism type of
the) fiber $(\bq\circ\pr)^{-1}(\phi)$ is independent of $\phi\in\dQ^\beta_y$.
The proof is absolutely similar to the one in ~\cite{fm} ~12.6.b). It makes
use of the {\em local convolution diagram} (see {\em loc. cit.} ~11.5).
We spare the reader the bulk of notation needed to carry this notion over to
our situation.

So we may (and will) choose $\phi$ lying in the fine Schubert stratum
$\oQ^0_y\times(C-0)^\beta_\Gamma\times0^{\alpha+\beta_0-\beta}
\subset\dQ^\beta_y$. Moreover, we will assume $\phi=y\times D\times
0^{\alpha+\beta_0-\beta}$ where $D\in(C-0)^\beta_\Gamma$,
and $y\in\bX_y=\oQ^0_y$ is the $\bH$-fixed point.
We may (and will) view the fiber $(\bq\circ\pr)^{-1}(\phi)$ as a subset of
$\osG_i$. As such, it is isomorphic to the intersection of $\osG_i$ with the
orbit $T^y_{y,\beta-\alpha-\beta_0}$ of $\bN_-^y((z))$ in $\sG$ passing through
the $\bH$-fixed point $(y,\beta-\alpha-\beta_0)\in\Omega$
(we use the identification $\CW_f\times Y\iso\Omega$ of ~\ref{fuck} along
with the fact that the set of $\bH$-fixed points of $\sG$ is naturally
identified with $\Omega$), cf. the proof of ~12.6.b) in {\em loc. cit.}
Here as always $\bN_-^y\subset\bG$ denotes the
subgroup $\dot{y}\bN_-\dot{y}{}^{-1}$ for a representative $\dot{y}\in N(\bH)$
of $y$. Now the Claim follows by the standard Bruhat-Tits theory.

This completes the proof of the Proposition. $\Box$

\subsection{}
\label{finish}
Comparing the above Proposition with ~\cite{l1} or ~\cite{s} ~4.1
we see that the $\CH$-module
$\hK(\sQ)$ is isomorphic to Lusztig's completed {\em periodic Hecke module}
$\hP$ (see ~\cite{s} ~\S6), and this isomorphism $\varphi$ takes
$[\CM(A)]\in\hK(\sQ)$ to a basic element $A\in\hP$.
Moreover, comparing ~\ref{kato} with the Theorem 6.4.2 (due to Kato) of
{\em loc. cit.} we see that for $A=(w_0,\alpha)$ we have
$\varphi([\CL(A)])=\utP_A$ (notations of {\em loc. cit.} ~\S6).

The Verdier duality $\CalD:\ \mho(\sQ)\lra\mho(\sQ)$ gives rise to an
involution $N\mapsto\CalD N$ of $\hK(\sQ)$. As Verdier duality commutes with
the convolution, this involution is compatible with the involution
$H\mapsto\ol{H}$ of $\CH$, that is, $\CalD(H(N))=\ol{H}(\CalD N)$.
Certainly, $\CalD[\CL(A)]=[\CL(A)]$ for any $A\in\fA$.

Recall the Lusztig's involution $P\mapsto\ol{P}$ of the periodic module
$\hP$ (see e.g. {\em loc. cit.} ~4.3). It is also compatible with the
involution on the Hecke algebra, i.e. $\ol{H(P)}=\ol{H}(\ol{P})$. It also
preserves the elements $\utP_A,\ A\in\fA$ (see {\em loc. cit.} ~6.4).
Since the elements $\utP_A,\ A=(w_0,\alpha)\in\fA$ clearly generate the
$\CH$-module $\hP$, we conclude that $\varphi(\CalD N)=\ol{\varphi(N)}$
for any $N\in\hK(\sQ)$.

In other words, we may identify $\CH$-modules $\hP$ and $\hK(\sQ)$ with
their involutions and standard bases. Now the Theorem ~6.4.2 of {\em loc. cit.}
carries over to $\hK(\sQ)$ and claims that for any $A\in\fA$ the
element $\varphi^{-1}(\utP_A)=\sum_B(-1)^{d(A,B)}\ol{q}_{B,A}[\CM(B)]$ is the
only $\CalD$-invariant element of $\hK(\sQ)$ lying in
$[\CM(A)]+\sum_{B\leq A}v^{-1}\BZ[v^{-1}][\CM(B)]$.
But $[\CL(A)]$ is also $\CalD$-invariant and it also lies in
$[\CM(A)]+\sum_{B\leq A}v^{-1}\BZ[v^{-1}][\CM(B)]$ by the definition of
Goresky-MacPherson extension.
This completes the proof of the Theorem ~\ref{main}. $\Box$

\section{Tilting Conjectures}

\subsection{} Recall the notations of ~\cite{fm} ~12.12, ~12.13. So let
$\eta\in Y^+,\ \beta\in Y,\ \eta+\beta\in\BN[I]$. For a sheaf
$\CF\in\CP(\oCG_\eta,\bI)$ we consider the complex $\bc^\beta_\CQ(\CF)=
\bq_*(\IC(\fQ^\beta)\otimes\bp^*\CF)[-\dim\ufM^\eta]$ on $\CQ^{\eta+\beta}$.
If $\eta+\beta\not\in\BN[I]$ we set $\CQ^{\eta+\beta}=\emptyset$, and
$\bc^\beta_\CQ(\CF)=0$. We will define a natural action of $\fg^L$ on
$\oplus_{\beta\in Y}H^\bullet(\CQ^{\eta+\beta},\bc^\beta_\CQ(\CF))=
\oplus_{\beta\in Y}H^\bullet(\CG\CQ^\beta_\eta,\IC(\fQ^\beta)\otimes
\bp^*\CF)$ closely following ~\ref{smirnoff}--\ref{h}.

\subsubsection{} For $\beta\in Y,\ \alpha\in\BN[I]$ we have the usual
{\em twisting map} $\sigma_{\beta,\alpha}:\ \CG\CQ^\beta_\eta\times C^\alpha
\lra\CG\CQ^{\alpha+\beta}_\eta$ (resp. $\fQ^\beta\times C^\alpha\lra
\fQ^{\alpha+\beta}$). We define
$\del_\alpha\CG\CQ^{\alpha+\beta}_\eta$ (resp. $\del_\alpha\fQ^{\alpha+\beta}$)
as the image of this map. We will construct the {\em stabilization map}
$\fv_{\alpha,\beta}:\ \CA_\alpha\lra
\Ext_{\CG\CQ^{\alpha+\beta}_\eta}^{|\alpha|}
(\IC(\del_\alpha\fQ^{\alpha+\beta})\otimes\bp^*\CF,
\IC(\fQ^{\alpha+\beta})\otimes\bp^*\CF)$.

We will follow the construction of ~\ref{bezr}. First, tensoring with 
$\Id\in\Ext^0_{\CG\CQ^\beta_\eta}(\IC(\fQ^\beta)\otimes\bp^*\CF,
\IC(\fQ^\beta)\otimes\bp^*\CF)$, we obtain the map 
from $\CA_\alpha=\Ext^{|\alpha|}_{C^\alpha}(\IC(C^\alpha),\fF_\alpha)$
to $\Ext^{|\alpha|}_{\CG\CQ^\beta_\eta\times C^\alpha}
(\IC(\fQ^\beta)\otimes\bp^*\CF\boxtimes\IC(C^\alpha),
\IC(\fQ^\beta)\otimes\bp^*\CF\boxtimes\fF_\alpha)$ (notations of ~\ref{bezr}).

As in ~\ref{bezr}, there is a canonical map $c$ from 
$\IC(\fQ^\beta)\otimes\bp^*\CF\boxtimes\fF_\alpha$ to $\sigma_{\beta,\alpha}^!
(\IC(\fQ^{\alpha+\beta})\otimes\bp^*\CF)$ extending the factorization 
isomorphism from an open part $\oGQ^\beta_\eta\times\BA^\alpha\subset
\CG\CQ^\beta_\eta\times C^\alpha$. Thus $c$ induces the map from 
$\Ext^{|\alpha|}_{\CG\CQ^\beta_\eta\times C^\alpha}
(\IC(\fQ^\beta)\otimes\bp^*\CF\boxtimes\IC(C^\alpha),
\IC(\fQ^\beta)\otimes\bp^*\CF\boxtimes\fF_\alpha)$ to
$\Ext^{|\alpha|}_{\CG\CQ^\beta_\eta\times C^\alpha}
(\IC(\fQ^\beta)\otimes\bp^*\CF\boxtimes\IC(C^\alpha),
\sigma_{\beta,\alpha}^!
(\IC(\fQ^{\alpha+\beta})\otimes\bp^*\CF))$.

The latter space equals 
$\Ext^{|\alpha|}_{\CG\CQ^{\alpha+\beta}_\eta}
((\sigma_{\beta,\alpha})_*
(\IC(\fQ^\beta)\otimes\bp^*\CF\boxtimes\IC(C^\alpha)),
\IC(\fQ^{\alpha+\beta})\otimes\bp^*\CF)=
\Ext_{\CG\CQ^{\alpha+\beta}_\eta}^{|\alpha|}
(\IC(\del_\alpha\fQ^{\alpha+\beta})\otimes\bp^*\CF,
\IC(\fQ^{\alpha+\beta})\otimes\bp^*\CF)$.

Finally, we define $\fv_{\alpha,\beta}$ as the composition of above maps.

\subsubsection{} We construct the {\em costabilization map}
$$\fy_{\beta,\alpha}:\ H^\bullet(\CG\CQ_\eta^\beta,
\IC(\fQ^\beta)\otimes\bp^*\CF)\lra
H^{\bullet-|\alpha|}(\del_\alpha\CG\CQ_\eta^{\alpha+\beta},
\IC(\del_\alpha\fQ^{\alpha+\beta})\otimes\bp^*\CF)$$
To this end we note that exactly as in ~\ref{normal}, we have
$\IC(\del_\alpha\fQ^{\alpha+\beta})=(\sigma_{\beta,\alpha})_*
\IC(\fQ^\beta\times C^\alpha)=(\sigma_{\beta,\alpha})_*
(\IC(\fQ^\beta)\boxtimes\ul\BC[|\alpha|])$.
Let $[C^\alpha]\in H^{-|\alpha|}(C^\alpha,\ul\BC[|\alpha|])$ denote the
fundamental class of $C^\alpha$.

Now, for $h\in H^\bullet(\CG\CQ_\eta^\beta,
\IC(\fQ^\beta)\otimes\bp^*\CF)$ we define
$\fy_{\beta,\alpha}(h)$ as $h\otimes[C^\alpha]\in
H^{\bullet-|\alpha|}(\CG\CQ_\eta^\beta\times C^\alpha,
\IC(\fQ^\beta)\otimes\bp^*\CF\boxtimes
\ul\BC[|\alpha|])=H^{\bullet-|\alpha|}(\del_\alpha\CG\CQ_\eta^{\alpha+\beta},
(\sigma_{\beta,\alpha})_*(\IC(\fQ^\beta)\otimes\bp^*\CF
\boxtimes\ul\BC[|\alpha|]))=
H^{\bullet-|\alpha|}(\del_\alpha\CG\CQ_\eta^{\alpha+\beta},
\IC(\del_\alpha\fQ^{\alpha+\beta})\otimes\bp^*\CF)$.

\subsection{Definition}
\label{Roma}
Let $a\in\CA_\alpha,\ h\in H^\bullet(\CG\CQ_\eta^\beta,
\IC(\fQ^\beta)\otimes\bp^*\CF)$. We define
the action $a(h)\in H^\bullet(\CG\CQ_\eta^{\alpha+\beta},
\IC(\fQ^{\alpha+\beta})\otimes\bp^*\CF)$
as the action of $\fv_{\alpha,\beta}(a)$ on the global cohomology applied
to $\fy_{\beta,\alpha}(h)$.

Let us stress that the action of $\CA_\alpha$ {\em preserves cohomological
degrees.}

\subsubsection{} For $a\in\CA_\alpha,\ b\in\CA_\beta,\
h\in H^\bullet(\CG\CQ_\eta^\beta,\IC(\fQ^\beta)\otimes\bp^*\CF)$
we have $a(b(h))=a\cdot b(h)$.
The proof is entirely similar to the proof of associativity of the
multiplication in $\CA$.

\subsection{}
\label{FF}
For $\beta\in Y$ the graded space
$H^\bullet(\CG\CQ_\eta^\beta,\IC(\fQ^\beta)\otimes\bp^*\CF)$ is
up to a shift Poincar\'e dual to
$H^\bullet(\CG\CQ_\eta^\beta,\IC(\fQ^\beta)\otimes\bp^*\CF)$
(cf. ~\cite{fm} ~12.1b). We define the map
$$f_i:\ H^\bullet(\CG\CQ_\eta^{\beta+i},
\IC(\fQ^{\beta+i})\otimes\bp^*\CF)\lra
H^\bullet(\CG\CQ_\eta^\beta,\IC(\fQ^\beta)\otimes\bp^*\CF)$$
as the dual of the map
$$e_i:\ H^\bullet(\CG\CQ_\eta^\beta,\IC(\fQ^\beta)\otimes\bp^*\CF)\lra
H^\bullet(\CG\CQ_\eta^{\beta+i},\IC(\fQ^{\beta+i})\otimes\bp^*\CF)$$
(see ~\ref{choice}).

It follows from ~\ref{Serre} that the maps $f_i,\ i\in I$, satisfy the Serre
relations of $\fn_-^L$.

\subsubsection{}
\label{ff}
We safely leave to the reader an elementary construction of the operators
$f_i,e_i$ absolutely similar to the one in ~\ref{f},~\ref{e}.

\subsection{}
\label{hh}
For $i\in I$ we define the endomorphism $h_i$ of
$H^\bullet(\CG\CQ_\eta^\beta,\IC(\fQ^\beta)\otimes\bp^*\CF)$
as the scalar multiplication by
$\langle\beta+2\crho,i'\rangle$ where $i'\in X$ is the simple root.
Exactly as in ~\ref{fe},~\ref{h}, one proves that the operators
$e_i,f_i,h_i,\ i\in I$, generate the action of the Langlands dual Lie algebra
$\fg^L$ on
$\oplus_{\beta\in Y}H^\bullet(\CG\CQ^\beta_\eta,\IC(\fQ^\beta)\otimes\bp^*\CF)=
\oplus_{\beta\in Y}H^\bullet(\CQ^{\eta+\beta},\bc^\beta_\CQ(\CF))$.

In the particular case $\eta=0,\ \CF$ is the irreducible skyscraper sheaf,
we have $\bc^\beta_\CQ(\CF)=\IC(\CQ^\beta)$, and we return to the
$\fg^L$-action on $\oplus_{\beta\in\BN[I]}H^\bullet(\CQ^\beta,\IC(\CQ^\beta))$
defined in ~\ref{!}.

We will formulate a few conjectures concerning this action,
following ~\cite{fk} ~\S6.

\subsection{}
\label{roman}
Let $\CN^L\subset\fg^L$ be the nilpotent cone. The Lie algebra $\fg^L$
acts on the cohomology $H^{\dim\bX}_{\fn_-^L}(\CN^L,\CO)$ of the structure
sheaf of $\CN^L$ with supports in $\fn_-^L$. The character of this module
is well known to be $\frac{|\CW_f|e^{2\check\rho}}
{\prod_{\ctheta\in{\check\CR}{}^+}(1-e^\ctheta)^2}$. Thus it coincides with
the character of $\fg^L$-module
$\oplus_{\beta\in\BN[I]}H^\bullet(\CQ^\beta,\IC(\CQ^\beta))$.

{\bf Conjecture.} $\fg^L$-modules
$\oplus_{\beta\in\BN[I]}H^\bullet(\CQ^\beta,\IC(\CQ^\beta))$ and
$H^{\dim\bX}_{\fn_-^L}(\CN^L,\CO)$ are isomorphic.

As explained in ~\cite{fk} ~\S6 the module
$H^{\dim\bX}_{\fn_-^L}(\CN^L,\CO)$ is {\em tilting}, and to prove the
Conjecture it suffices to check that
$\oplus_{\beta\in\BN[I]}H^\bullet(\CQ^\beta,\IC(\CQ^\beta))$ is tilting as
well. To this end it is enough to check that the action of the algebra
$\CA=U(\fn_+^L)$ on
$\oplus_{\beta\in\BN[I]}H^\bullet(\CQ^\beta,\IC(\CQ^\beta))$ is free.

\subsection{}
\label{roma}
For $\eta\in Y^+$ the irreducible $\fg^L$-module with the highest weight
$\eta$ is denoted by $W_\eta$. The Theorem ~13.2 of ~\cite{fm} implies that
the character of $\fg^L$-module
$\oplus_{\beta\in Y}H^\bullet(\CQ^\beta,\bc^\beta_\CQ(\IC(\oCG_\eta)))$
equals char($W_\eta)\times\frac{|\CW_f|e^{2\check\rho}}
{\prod_{\ctheta\in{\check\CR}{}^+}(1-e^\ctheta)^2}$.

{\bf Conjecture.} There is an isomorphism of $\fg^L$-modules
$\oplus_{\beta\in Y}H^\bullet(\CQ^\beta,\bc^\beta_\CQ(\IC(\oCG_\eta)))$ and
$W_\eta\otimes H^{\dim\bX}_{\fn_-^L}(\CN^L,\CO)$.

This Conjecture is again equivalent to the statement that the $\fg^L$-module
$\oplus_{\beta\in Y}H^\bullet(\CQ^\beta,\bc^\beta_\CQ(\IC(\oCG_\eta)))$ is
tilting. To check this it suffices to prove that the action of $\CA=U(\fn_+^L)$
on $\oplus_{\beta\in Y}H^\bullet(\CQ^\beta,\bc^\beta_\CQ(\IC(\oCG_\eta)))$ is
free.

\subsection{}
\label{romka}
Let $\CF$ be a $\bG[[z]]$-equivariant perverse sheaf on the affine Grassmannian
$\CG$. According to ~\cite{mv}, $H^\bullet(\CG,\CF)$ is equipped
with a canonical $\fg^L$-action. While the action of ${\fh}^L\subset\fg^L$
is defined pretty explicitly in geometric terms in ~\cite{mv}, the actions
of $\fn^L_-,\fn^L_+\subset\fg^L$ were not constructed as explicitly so far.
Motivated by the Conjecture ~\ref{roma} above we propose the following
conjectural constructions of the $\fn^L_\pm$-actions on $H^\bullet(\CG,\CF)$.

\subsubsection{}
Let $\CQ^{\eta+\beta}=\sqcup_{0\leq\gamma\leq\eta+\beta}\dQ^\gamma$
be the partition of $\CQ^{\eta+\beta}$ according to the defect at $0\in C$.
Let $\CQ^{\eta+\beta}=\cup_{0\leq\gamma\leq\eta+\beta}\CQ^\gamma$
be the corresponding filtration by closed subspaces.
In particular, we have $\bX=\dQ^0=\CQ^0\subset\CQ^{\eta+\beta}$. The perverse
sheaf $\bc_\CQ^\beta(\CF)$ is isomorphic to a direct sum of a few copies
of the sheaves $\IC(\CQ^\gamma)$. More precisely, the multiplicity of
$\IC(\CQ^\gamma)$ in $\bc_\CQ^\beta(\CF)$ equals
$H^{2|\beta-\gamma|}_c(T_{\gamma-\beta},\CF)$
(see the Theorem ~13.2 of ~\cite{fm}).
In particular, $H^{2|\beta|}_c(T_{-\beta},\CF)\otimes\IC(\CQ^0)=
H^{2|\beta|}_c(T_{-\beta},\CF)\otimes\ul\BC_\bX[n]$
is canonically a direct summand
of $\bc_\CQ^\beta(\CF)$. Hence
$H^{2|\beta|}_c(T_{-\beta},\CF)=
H^{2|\beta|}_c(T_{-\beta},\CF)\otimes H^{-n}(\bX,\ul\BC[n])$
is canonically a direct summand of
$H^\bullet(\CQ^{\eta+\beta},\bc_\CQ^\beta(\CF))=
H^\bullet(\CG\CQ^\beta_\eta,\IC(\fQ^\beta)\otimes\bp^*\CF)$.

We conjecture that
$\oplus_{\beta\in Y}
H^{2|\beta|}_c(T_{-\beta},\CF)=\oplus_{\beta\in Y}
H^{2|\beta|}_c(T_{-\beta},\CF)\otimes H^{-n}(\bX,\ul\BC[n])$
is an $\CA^{opp}=U(\fn^L_-)$-submodule of
$\oplus_{\beta\in Y}
H^\bullet(\CG\CQ^\beta_\eta,\IC(\fQ^\beta)\otimes\bp^*\CF)$
(with respect to the action defined in ~\ref{Roma}), and
the resulting action of $U(\fn^L_-)$ on
$\oplus_{\beta\in Y}H^{2|\beta|}_c(T_{-\beta},\CF)=H^\bullet(\CG,\CF)$
coincides with the action of ~\cite{mv}.

\subsubsection{}
Similarly,
$H^{2|\beta|}_c(T_{-\beta},\CF)=
H^{2|\beta|}_c(T_{-\beta},\CF)\otimes H^{n}(\bX,\ul\BC[n])$
is canonically a direct summand of
$H^\bullet(\CQ^{\eta+\beta},\bc_\CQ^\beta(\CF))=
H^\bullet(\CG\CQ^\beta_\eta,\IC(\fQ^\beta)\otimes\bp^*\CF)$.
We conjecture that the kernel of the natural projection
$\oplus_{\beta\in Y}
H^\bullet(\CG\CQ^\beta_\eta,\IC(\fQ^\beta)\otimes\bp^*\CF)\lra
\oplus_{\beta\in Y}
H^{2|\beta|}_c(T_{-\beta},\CF)=\oplus_{\beta\in Y}
H^{2|\beta|}_c(T_{-\beta},\CF)\otimes H^{n}(\bX,\ul\BC[n])$
is invariant under the action of $\CA=U(\fn^L_+)$ defined in
~\ref{FF},
and the induced action of $U(\fn^L_+)$ on
$\oplus_{\beta\in Y}H^{2|\beta|}_c(T_{-\beta},\CF)=H^\bullet(\CG,\CF)$
coincides with the action of ~\cite{mv}.

\subsection{}
\label{denis}
Recall the moduli scheme $\fM$ of $\bG$-torsors on $\BP^1$ equipped with
the formal trivialization at $\infty\in\BP^1$. Its stratification according
to the isomorphism classes of $\bG$-torsors $\fM=\sqcup_{\eta\in Y^+}\fM_\eta$
was described in ~\cite{fm} ~10.6. Let $\ol\fM_\eta$ denote the closure of
the stratum $\fM_\eta$, and let $\IC(\ol\fM_\eta)$ denote the corresponding
$\IC$-sheaf. According to the Decomposition Theorem, for $\alpha\in Y$,
the direct image $\bp_*\IC(\fQ^\xi)$ is a semisimple complex on $\fM$,
isomorphic to a direct sum of various $\IC(\ol\fM_\eta)$ with shifts and
multiplicities. We have $\bigoplus_{\xi\in Y}\bp_*\IC(\fQ^\xi)=
\bigoplus_{\eta\in Y^+}K^\bullet_\eta\otimes\IC(\ol\fM_\eta)$ 
for (infinite dimensional
in general) graded multiplicities spaces $K_\eta$.

One can introduce the stalkwise action of $\fg^L$ on
$\bigoplus_{\xi\in Y}\bp_*\IC(\fQ^\xi)$ entirely similar to ~\ref{hh}.
In other words, one obtains the action of $\fg^L$ on the multiplicities
spaces $K^\bullet_\eta$ for any $\eta\in Y^+$. For instance, 
$K^\bullet_0=\oplus_{\beta\in\BN[I]}H^\bullet(\CQ^\beta,\IC(\CQ^\beta))$
by definition. The action of $\fg^L$ on $K^\bullet_0$ was described in
~\ref{roman}. We will extend the conjecture ~\ref{roman} to the case of
arbitrary $\eta\in Y^+$.

Recall the Lusztig's quantum groups $\fu\subset\fU$, and the representations'
category $\fC$, introduced in
~\cite{fm} ~1.3. The highest weights of the regular block $\fC^0$ of the
category of $\fU$-modules are contained in the $\CW$-orbit 
$\CW\cdot0\subset X$ (the strong linkage principle). Let $m_\eta\in\CW$
be the shortest representative of the double coset $\CW_f\eta\CW_f$ in $\CW$.
Let $T(m_\eta\cdot0)\in\fC^0$ be the indecomposable {\em tilting} $\fU$-module
with the highest weight $m_\eta\cdot0$ (see ~\cite{s}). 
For example, $T(m_0\cdot0)$ is the
trivial module. Finally, recall the notion of semiinfinite cohomology 
$H^{\frac{\infty}{2}+\bullet}(\fu,?)$ introduced in ~\cite{a}. It is known
that for $T\in\fC^0$ the semiinfinite cohomology
$H^{\frac{\infty}{2}+\bullet}(\fu,T)$ carries a natural structure of
$\fn_-^L$-integrable $\fg^L$-module.

{\bf Conjecture.} For $\eta\in Y^+$ there is an isomorphism of graded
$\fg^L$-modules $K^\bullet_\eta\simeq
H^{\frac{\infty}{2}+\bullet}(\fu,T(m_\eta\cdot0))$. Under this isomorphism,
the weight $\xi\in Y$ part of the RHS corresponds to the multiplicity
$_{(\xi)}K^\bullet_\eta\subset K^\bullet_\eta$ of $\IC(\ol\fM_\eta)$
in $\bp_*\IC(\fQ^\xi)$.

{\bf Corollary of the Conjecture.} Suppose $\eta-2\check{\rho}\in Y^+$.
Then $K^\bullet_\eta$ is concentrated in degree 0, and is isomorphic
to the irreducible $\bG^L$-module $W_{\eta-2\check\rho}$.

{\em Remark.} S.Arkhipov has recently proved that 
$H^{\frac{\infty}{2}+\bullet}(\fu,T(0))\simeq 
H^{\dim\bX}_{\fn^L_-}(\CN^L,\CO)$. Together with the results of ~\cite{fkm}
this establishes the $\eta=0$ case of the above conjecture.

\section{Appendix. Kontsevich resolution of Quasimaps' space}

\subsection{}
For $\alpha=\sum_{i\in I}a_ii\in\BN[I]$ let $\CQ^\alpha_K=
\ol{M}_{0,0}(\BP^1\times\bX),(1,\alpha))$ be the Kontsevich space of stable
maps from the genus zero curves without marked points to $\BP^1\times\bX$
of bidegree $(1,\alpha)$ (see ~\cite{ko}). In this Appendix we construct
a regular birational map $\pi:\ \CQ^\alpha_K\lra\CQ^\alpha$.

To this end recall that $\bX$ is a closed subvariety in
$\prod_{i\in I}\BP(V_{\omega_i})$ given by Pl\"ucker relations.
Drinfeld's space $\CQ^\alpha$ is a closed subvariety in
$\prod_{i\in I}\BP H^0(C,V_{\omega_i}\otimes\CO(a_i))$ also given by
Pl\"ucker relations.

So to construct the desired map it suffices to solve the following problem.
For a vector space $V$ one has 2 compactifications of the space $M^a$
of algebraic maps $f:\ C\lra\BP(V)$ of degree $a\in H_2(\BP(V),\BZ)$.
Kontsevich compactification $M^a_K=\overline{M}_{0,0}(C\times\BP(V),(1,a))$
is the space of all stable maps from the genus 0 curves into the product
$C\times\BP(V)$ of bidegree $(1,a)$; the embedding of $M^a$ into $M^a_K$
takes $f$ to its graph $\Gamma_f$.
Drinfeld compactification $M^a_D$ is the space
of invertible subsheaves of degree $-a$ in $V\otimes\CO_C$. Twisting by
$\CO(a)$ we identify $M^a_D$ with $\BP H^0(C,V\otimes\CO(a))$. The embedding
of $M^a$ into $M^a_D$ takes $f$ to the line subbundle
$\CL_f:=p_*q^*\CO_{\BP(V)}(-1)\subset V\otimes\CO_C$ where
$C\stackrel{p}{\longleftarrow}\Gamma_f\stackrel{q}{\lra}\BP(V)$ are the
natural projections. We want to prove that the identification
$M^a_K\supset M^a=M^a\subset M^a_D$ extends to the regular map
$\pi:\ M^a_K\lra M^a_D$.

This problem was addressed by A.Givental in ~\cite{g}, ~``Main Lemma''.
Unfortunately, the validity of his proof is controversial. In the rest of the
Appendix we prove that Id: $M^a\iso M^a$ extends to the regular map
$\pi:\ M^a_K\lra M^a_D$.

\subsection{}
\label{kuznec}
One of the most important properties of the Drinfeld space $M^a_D$
is the following.

{\em Lemma.}
Let $S$ be a scheme and $\CL\subset V\otimes\CO_{C\times S}$ be
an invertible subsheaf in the trivial vector bundle $V\otimes\CO_{C\times S}$
on $C\times S$ of degree $-a$ over $S$ (i.e.\ for any $s\in S$
the restriction $\CL_s$ of $\CL$ to the fiber $C\times s\subset C\times S$
is equal to $\CO(-a)$). If the open subset $U=\{(x,s)\in C\times S\ |\
\CL_{(x,s)}\to V\otimes\CO_{x,s} \text{ is embedding}\}$ contains generic
points of all fibers of $C\times S$ over $S$ then the map
$f:S\to M^a_D$ sending each point $s\in S$ to the subsheaf
$\CL_s\subset V\otimes\CO_{C\times s}\cong C\otimes\CO_C$ is regular.

{\em Proof.}
Let $p:C\times S\to S$ be the natural projection and
$L=p_*(\CL\otimes\CO_C(a))$. It is evident that for any $s\in S$
we have
$$
\begin{array}{rcrcl}
H^0(p^{-1}(s),\CL\otimes\CO_C(a))&=&H^0(C,\CO_C)&=&\BC,\\
H^{>0}(p^{-1}(s),\CL\otimes\CO_C(a))&=&H^{>0}(C,\CO_C)&=&0
\end{array}
$$
and the map
$$
H^0(p^{-1}(s),\CL\otimes\CO_C(a))\to H^0(p^{-1}(s),V\otimes\CO_C(a))
$$
is injection, hence $L$ is a line subbundle in
$p_*(V\otimes\CO_C(a))=H^0(C,V\otimes\CO_C(a))\otimes\CO_S$.
The map $f$ is just the regular map associated with the subbundle $L$.
$\Box$

\subsection{}
Recall the statement of Givental's ``Main Lemma''.

{\bf ``Main Lemma''}
a) The map Id: $M^a\iso M^a$
extends to the regular map $\pi:M^a_K\to M^a_D$;

b) Let $(\vphi:\CC\to C\times \BP(V))\in M^a_K$ be a stable map and let
$\vphi':\CC\to C$, $\vphi'':\CC\to \BP(V)$ denote the induced maps. Let $\CC_0$
denote the irreducible component of $\CC$ such that
$\vphi':\CC_0\to C$ is dominant,
and let $\CC_1,\dots,\CC_m$ be the connected components of $\CC\setminus\CC_0$.
Finally, let $x_1=\vphi'(\CC_1),\dots,x_m=\vphi'(\CC_m)$. Then
$\pi(\CC)$ is the subsheaf with normalisation equal to
$\vphi'_*(\vphi''_{|\CC_0})^*\CO_{\BP(V)}(-1)$ 
and with defect at points $x_1$, \dots
$x_m$ of degree $\vphi''_*[\CC_1]$, \dots, $\vphi''_*[\CC_m]$ respectively.

{\em Proof.}
a) Consider the space $\tmk=\overline M_{0,1}(C\times \BP(V);(1,a))$ of
stable maps with one marked point. There are natural maps
${\sf f}:\tmk\to M^a_K$ forgeting the marked point, and
$\ev:\tmk\to C\times \BP(V)$ evaluating at the marked point.
Let $p:\tmk\to C\times M^a_K$ denote the product of the
map $p_1\circ\ev:\tmk @>\ev>> C\times \BP(V) @>p_1>> C$ and of
the map ${\sf f}$. Let $q:\tmk\to \BP(V)$ denote the map
$p_2\circ\ev:\tmk @>\ev>> C\times \BP(V) @>p_2>> \BP(V)$.

The map $p$ is an isomorphism over the open subspace
$U\subset C\times M^a_K$ of pairs $(x,\CC)$ such that
$x\in\vphi'(\CC_0-\CC^{\operatorname{sing}})$.
%, where $\vphi:\CC\to C\times \BP(V)$
%is a stable map, $\CC_0$ is irreducible component of $\CC$
%which maps dominant onto $C$ by the map $\vphi_1=p_1\cdot\vphi$.
Consider the sheaf $\CB=p_*q^*\CO_{\BP(V)}(-1)\subset
p_*q^*(V\otimes\CO_{\BP(V)})=
V\otimes\CO_{C\times M^a_K}$. Its reflexive hull $\CB^{**}$ is
a reflexive rank 1 sheaf on the smooth orbifold $C\times M^a_K$,
hence $\CB^{**}$ is an invertible sheaf. Note that since $\CB$ is
invertible over $U$, we have $\CB^{**}_{|U}=\CB_{|U}$.
Since $U$ contains $C\times M^a$, the restriction of $\CB^{**}$ to
the generic fiber of $C\times M^a_K$ over $M^a_K$ is $\CO(-a)$,
hence the restriction of $\CB^{**}$ to any fiber is $\CO(-a)$ as well.
Over the set $U$, the embedding $\CB^{**}_{|U}\to V\otimes\CO_U$ is an
embedding of vector bundles (because over $U$ the map $p$ is an
isomorphism). On the other hand, $U$ contains generic points of all
fibers of $C\times M^a_K$ over $M^a_K$, therefore we can apply
the Lemma ~\ref{kuznec} and the first part of the ``Main Lemma'' follows.

b) Let $s=(\vphi:\CC\to C\times \BP(V))\in M^a_K$ be a stable map, and
let $L\subset V\otimes\CO_C$ denote a subsheaf $\pi(s)\in M^a_D$;
let $L'=\vphi'_*({\vphi''_{|\CC_0}}^*)\CO_{\BP(V)}(-1)$. The Lemma
~\ref{kuznec} implies
$L=\CB^{**}_s$. We want to prove that the normalization $\widetilde L$ of $L$
equals $L'$.

Denote $C-\{x_1,\dots,x_m\}$ by $C^0$.
Then $U\bigcap(C\times s)=C^0\times s$, hence
$$%\begin{multline*}
L_{|C^0\times s}={\CB_s}_{|C^0\times s}=
(\vphi'_*{\vphi''}^*\CO_{\BP(V)}(-1))_{|C^0\times s}=
(\vphi'_*(\vphi''_{|\CC_0})^*\CO_{\BP(V)}(-1))_{|C^0\times s}=
L'_{|C^0\times s},
$$%\end{multline*}
i.e. $L$ and $L'$ coincide on the open subset of $C\times s$.
Since $L'$ is a line subbundle in $V\otimes\CO_{C\times s}$ this means that
$L$ is a subsheaf in $L'$ and the quotient $L'/L$ is concentrated
at the points $x_1,\dots,x_m$. Therefore $\widetilde L=L'$ and the defect
of $L$ is concentrated at the points $x_1,\dots,x_m$.

Note that in the case $m=1$ the ``Main Lemma'' follows. In the general case
we proceed as follows. Let $c_k=\vphi''_*[\CC_k]$, $(k=1,\dots,m)$.
Let ${\bf M}_K^k=\{(\psi:{\sf C}\to C\times \BP(V))\}\subset M^a_K$ denote the
subspace of all stable maps such that a curve ${\sf C}$ has two
irreducible components ${\sf C}={\sf C}_0\cup{\sf C}_1$,
$\psi_*[{\sf C}_0]=(1,a-c_k)$,
$\psi_*[{\sf C}_1]=(0,c_k)$ and $\psi'({\sf C}_1)=x_k$. Then we have
$s\in\overline{\bf M}_K^1\bigcap\overline{\bf M}_K^2\bigcap\dots
\bigcap\overline{\bf M}_K^m.$

Let ${\bf M}_D^k\subset M^a_D$ denote the subspace of all invertible
subsheaves with defect of degree $c_k$ concentrated at the point $x_k$.
According to the case $m=1$ we have $\pi({\bf M}_K^k)\subset{\bf M}_D^k$, hence
$\pi(s)\in\overline{\bf M}_D^1\bigcap\overline{\bf M}_D^2\bigcap\dots
\bigcap\overline{\bf M}_D^m.$
Since the degree of defect at a given point increases under specialization,
the degrees of the defect of $\pi(s)$ at the points $x_k$ are greater or
equal than $c_k$, hence they are equal to $c_k$,
 and the ``Main Lemma'' follows. $\Box$

\end{document}